\documentclass[10pt,a4paper]{article}

\usepackage{amsmath,amssymb, bbm}
\usepackage{color}
\usepackage{float}
\newcommand{\R}{\mathbb{R}}
\newcommand{\N}{\mathbb{N}}
\newcommand{\Z}{\mathbb{Z}}

\newcommand{\Co}{\mathbb{C}_0}

\newcommand{\bP}{\mathbb{P}}
\def\cA{{\mathcal A}}
\def\cC{{\mathcal C}}
\def\cB{{\mathcal B}}

\def\cF{{\mathcal F}}
\def\cP{{\mathcal P}}
\def\cQ{{\mathcal Q}}
\newcommand{\ee}{\varepsilon}
\renewcommand{\aa}{\alpha}
\newcommand{\bb}{\beta}

\renewcommand{\div}{{\rm div}\,}

\newcommand{\Hs}{\dot{H^s}}

\def\d{\partial}

\def\ddj{\dot \Delta_j}

\def\tilde{\widetilde}
\def\hat{\widehat}
\newcommand{\D}{\Delta}

\newcommand{\n}{\nabla}

\newcommand{\Fe}{F_{ext}}

\newcommand{\tOm}{\tilde{\Omega}_{QG}}

\newcommand{\ome}{\omega_\varepsilon}
\newcommand{\tom}{\tilde{\omega}}

\newcommand{\voe}{v_{0,\ee}}
\newcommand{\ub}{\bar{u}}

\newcommand{\Ue}{U_\ee}
\newcommand{\ve}{v_\ee}
\newcommand{\Ve}{V_\ee}
\newcommand{\He}{H_\ee}
\newcommand{\Ven}{V_\ee^n}
\newcommand{\Hen}{H_\ee^n}

\newcommand{\Uoe}{U_{0,\ee}}
\newcommand{\UoeS}{U_{0,\ee,S}}
\newcommand{\Uoeosc}{U_{0,\ee, osc}}
\newcommand{\tU}{\tilde{U}}

\newcommand{\De}{D_\ee}
\newcommand{\Den}{D_\ee^n}
\newcommand{\DeS}{D_{\ee,S}}
\newcommand{\Dosc}{D_{\ee, osc}}

\newcommand{\tD}{\tilde{D}}

\newcommand{\Pe}{P_\ee}
\newcommand{\tP}{\tilde{P}}

\newcommand{\Phie}{\Phi_\ee}
\newcommand{\Thee}{\theta_\ee}
\newcommand{\tThe}{\tilde{\theta}}
\newcommand{\tThee}{\tilde{\theta}_\ee}
\newcommand{\Theeo}{\theta_{0,\ee}}
\newcommand{\tTheo}{\tilde{\theta}_0}
\newcommand{\tTheeo}{\tilde{\theta}_{0,\ee}}

\newcommand{\tv}{\tilde{v}}

\newcommand{\tvo}{\tilde{v}_{0}}

\newcommand{\tX}{\tilde{X}}
\newcommand{\tY}{\tilde{Y}}
\newcommand{\tZ}{\tilde{Z}}
\newcommand{\tT}{\tilde{T}}
\newcommand{\tpi}{\tilde{\pi}}
\newcommand{\tq}{\tilde{q}}
\newcommand{\tg}{\tilde{g}}
\newcommand{\tG}{\tilde{G}}
\newcommand{\tZe}{\tilde{Z}_{\ee}}
\newcommand{\tKe}{\tilde{K}_{\ee}}

\newcommand{\qe}{q_{\ee}}
\newcommand{\qen}{q_{\ee}^n}

\newcommand{\We}{W_\ee}

\renewcommand{\Re}{R_\ee}
\newcommand{\re}{r_\ee}

\newcommand{\cPrR}{\cP_{\re, \Re}}
\newcommand{\cPrRb}{\cP_{\frac{\re}2, 2\Re}}

\newcommand{\exi}{(\ee,\xi)}

\newcommand{\IbR}{I_{1, \bb}^R(\sigma)}
\newcommand{\IbRd}{I_{1,\frac2{3\sqrt{3}}}^R(\sigma)}
\newcommand{\IbRo}{I_{1, \bb_0}^R(\sigma)}
\newcommand{\Jb}{J_{1,\bb}(\sigma)}
\newcommand{\Kb}{K_{1,\bb}(\sigma)}
\newcommand{\Jbb}{\overline{J_2}}
\newcommand{\Jub}{\overline{J_1}}

\newcommand{\zib}{z_1(\bb)}
\newcommand{\fdo}{f_1'(z_1(2\delta_0))}
\newcommand{\Kbb}{\overline{K_2}}
\newcommand{\ziib}{z_2(\bb)}
\newcommand{\bJ}{\overline{J}}

\newtheorem{thm}{Theorem}[section]
\newtheorem{lem}{Lemma}[section]

\newtheorem{prop}{Proposition}[section]
\newtheorem{defi}{Definition}[section]

\newtheorem{rem}{Remark}[section]

\usepackage{graphicx}
\usepackage{tikz}
\usepackage{scalerel}
\usepackage{pict2e}
\usepackage{tkz-euclide}
\usetikzlibrary{calc}
\usetikzlibrary{patterns,arrows.meta}
\usetikzlibrary{shadows}
\usetikzlibrary{external}
\usepackage{pgfplots}
\pgfplotsset{compat=newest}
\usepgfplotslibrary{statistics}
\usepgfplotslibrary{fillbetween}
\usepackage{xcolor}
\usepackage{nicefrac}
\usepackage{pgfplots}

\title{Hidden asymptotics for the weak solutions of the strongly stratified Boussinesq system without rotation}

\author{Fr\'ed\'eric Charve\footnote{Univ Paris Est Creteil, Univ Gustave Eiffel, CNRS, LAMA UMR8050, F-94010 Creteil, France. E-mail: frederic.charve@u-pec.fr}}

\date{}
\begin{document}

\maketitle

\begin{abstract} The asymptotics of the strongly stratified Boussinesq system when the Froude number goes to zero have been previously investigated, but the resulting limit system surprisingly did not depend on the thermal diffusivity $\nu'$. In this article we obtain richer asymptotics (depending on $\nu'$) for more general initial data.

As for the rotating fluids system, the only way to reach this limit consists in considering non-conventional initial data: to a function classically depending on the full space variable, we add a second one only depending on the vertical coordinate.

Thanks to a refined study of the structure of the limit system and to new adapted Strichartz estimates, we obtain convergence in the context of weak Leray-type solutions providing explicit convergence rates when possible. In the usually simpler case $\nu=\nu'$ we are able to improve the Strichartz estimates and the convergence rates. The last part of the appendix is devoted to the proof (by elementary techniques) of a new and crucial dispersion estimate, as classical methods fail.

Finally, our theorems can also be rewritten as a global existence result and asymptotic expansion for the classical Boussinesq system near an explicit stationary solution and for non-conventional vertically stratified initial data.
\end{abstract}

\section{Introduction}

\subsection{Geophysical fluids: Primitive and Rotating fluids systems}

It is commonly known that Geophysical fluids dynamics are greatly influenced by two concurrent phenomena: the Coriolis and centrifugal forces created by the rotation of the Earth around its axis, and the vertical stratification of the density induced by gravity. Both of them create rigidities whose influence on the dynamics is quantified by two non-dimensional numbers introduced by Physicists: the Rossby $Ro$ and Froude $Fr$ numbers. The smaller they are, the greater are the actions of these rigidities on the fluid dynamics.

Let us present a few models taking these forces into account. On one hand, the Primitive System (sometimes also called Primitive Equations) that we write here in the whole space and in the particular regime where both phenomena are of the same scale (that is we choose $Ro=\ee$ and $Fr=\ee F$ with $F>0$):
\begin{equation}
\begin{cases}
\d_t \Ue +\ve\cdot \n \Ue -L \Ue +\frac{1}{\ee} \cA \Ue=\frac{1}{\ee} (-\n \Phie, 0),\\
\div \ve=0,\\
{\Ue}_{|t=0}=U_{0,\ee}.
\end{cases}
\label{PE}
\tag{$PE_\ee$}
\end{equation}
The unknowns are $\Ue =(\ve, \Thee)=(\ve^1, \ve^2, \ve^3, \Thee)$, where $\ve$ denotes the velocity of the fluid and $\Thee$ the scalar potential temperature (linked to the density, temperature and salinity), and $\Phie$, which is called the geopotential and contains in particular the pressure term and the centrifugal force. The diffusion operator $L$ takes into account two heat regularization effects  and is defined by
$$
L\Ue \overset{\mbox{def}}{=} (\nu \D \ve, \nu' \D \Thee),
$$
where $\nu, \nu'>0$ respectively denote the kinematic viscosity and thermal diffusivity (both will be considered as viscosities in what follows). The last term $\ee^{-1}\cA$ gathers the rotation and stratification effects and the matrix $\cA$ is defined by
$$
\cA \overset{\mbox{def}}{=}\left(
\begin{array}{llll}
0 & -1 & 0 & 0\\
1 & 0 & 0 & 0\\
0 & 0 & 0 & F^{-1}\\
0 & 0 & -F^{-1} & 0
\end{array}
\right).
$$
On the other hand, the Rotating fluids system only focusses on rotational effects and can be seen as deduced from the previous system only considering the velocity (so we get rid of the last line and column of $\cA$), it is written as follows:

\begin{equation}
 \begin{cases}
  \d_t \ve +\ve\cdot \n \ve-\nu \D \ve +\frac{e_3\wedge \ve}{\ee} =-\n p_\ee,\\
  \div \ve=0,\\
  {\ve}_{|t=0}= v_0.
 \end{cases}
\label{RF} \tag{$RF_\ee$}
\end{equation}

Both of these systems are variations of the famous Navier-Stokes system, but each of them shows better behaviour induced by the special structure brought by their respective limit systems as $\ee$ goes to zero: the QG/oscillating structure for \eqref{PE}, and the 2D-3D structure for \eqref{RF}. We refer to \cite{CDGG, CDGGbook, BMN5, Dutrifoy2, FC1, FCPAA, FCRF} for more details. For a physical presentation of the models we refer to \cite{Cushman, EmMa, EmMa2} as well as the introduction of \cite{FCthese}. We refer to \cite{FCcompl, FCRF} for a small survey presenting recent articles devoted to Systems \eqref{RF} and \eqref{PE}.
\\

We will use the same notations as in \cite{FCPAA, FCcompl}: for $s\in \R$ and $T>0$ we define the spaces:
$$
\begin{cases}
\vspace{0.2cm}
 \dot{E}_T^s=\mathcal{C}_T(\Hs (\R^3)) \cap L_T^2(\dot{H}^{s+1}(\R^3)),\\
 \dot{B}_T^s=\mathcal{C}_T(\dot{B}_{2,1}^s (\R)) \cap L_T^1(\dot{B}_{2,1}^{s+2}(\R)),
\end{cases}
$$
endowed with the following norms (For \eqref{PE}, $\nu_0=\min(\nu,\nu')$):
$$
\begin{cases}
\vspace{0.2cm}
 \|f\|_{\dot{E}_T^s}^2 \overset{def}{=}\|f\|_{L_T^\infty \Hs }^2+\nu_0 \int_0^T \|f(\tau)\|_{\dot{H}^{s+1}}^2 d\tau,\\
 \|f\|_{\dot{B}_T^s} \overset{def}{=}\|f\|_{L_T^\infty \dot{B}_{2,1}^s}+\nu' \int_0^T \|f(\tau)\|_{\dot{B}_{2,1}^{s+2}} d\tau,
\end{cases}
$$
where $H^s(\R^3)$, $\dot{H}^s(\R^3)$ and $\dot{B}_{2,1}^s(\R)$ respectively denote the inhomogeneous and homogeneous Sobolev spaces of index $s\in \R$ and the homogeneous Besov space of indices $(s,2,1)$.\\
When $T=\infty$ we simply write $\dot{E}^s$ or $\dot{B}^s$ and the corresponding norms are understood as taken over $\R_+$ in time.

\textbf{Notation:} For any $\R^3$ or $\R^4$-valued vector field, we will write $f^h=(f^1,f^2)$ and will define $f\cdot \n f=\sum_{i=1}^3 f^i \d_i f$. So that for instance we will indifferently write $\ve\cdot \n\Ue=\Ue \cdot \n\Ue$. We also denote $\nabla_h f =(\d_1 f, \d_2 f)$ and $\D_h=\d_1^2+\d_2^2$.

\subsection{Strongly stratified Boussinesq system without rotation}

\subsubsection{Presentation of the model and quick survey of recent results}

In this article we will focus on the following system, that only takes into account the stratification aspects. Here we will denote as $Fr=\ee$ the Froude number (this corresponds to $\ee F$ for the Primitive system $(PE_\ee)$) and the model is written as follows:
\begin{equation}
\begin{cases}
\d_t \Ue +\ve\cdot \n \Ue -L \Ue +\frac{1}{\ee} \cB \Ue=\frac{1}{\ee} (-\n \Phie, 0),\\
\div \ve=0,\\
{\Ue}_{|t=0}=U_{0,\ee}.
\end{cases}
\label{Stratif}
\tag{$S_\ee$}
\end{equation}
The unknowns are the same as in the previous section $\Ue =(\ve, \Thee)=(\ve^1, \ve^2, \ve^3, \Thee)$ and $\Phie$. In what follows we will explicitely decompose it as the sum of the pressure term and another penalized pressure term that could be seen as an analoguous of the centrifugal force. The diffusion operator $L$ is also defined as previously with $\nu, \nu'>0$. In the last term, $\ee^{-1}\cB$ only takes into account stratification effects and the matrix $\cB$ is still skewsymmetric and defined by:
$$
\cB \overset{\mbox{def}}{=}\left(
\begin{array}{llll}
0 & 0 & 0 & 0\\
0 & 0 & 0 & 0\\
0 & 0 & 0 & 1\\
0 & 0 & -1 & 0
\end{array}
\right).
$$

Let us begin with a brief review about recent works dedicated to System \eqref{Stratif} or its inviscid version. Let us first mention that in \cite{Wid}, K. Widmayer considers a solution $\Ue$ of \eqref{Stratif}, in the inviscid setting $\nu=\nu'=0$, which belongs to $\cC ([0,T], H^N)$ (for some $T>0$ and $N\geq 6$) with $\|\Ue\|_{L_T^\infty H^N}\leq C$ (uniformly in $\ee$). If, in addition, the initial data is independent of $\ee$ and satisfies $\|U_0\|_{W^{5,1}}<\infty$, then the solution can be decomposed into two parts $\Ue=(\Ue^1,0,0)+\Ue^2$ such that for all $t>0$:
$$
\begin{cases}
 \Ue^2(t) \underset{\ee \rightarrow 0}{\longrightarrow} 0 \mbox{ in }W^{1,\infty}(\R^3),\\
 \Ue^1(t)\underset{\ee \rightarrow 0}{\longrightarrow} \ub(t)=(\ub^1(t), \ub^2(t)) \mbox{ in }L^2(\R^3),
\end{cases}
$$
where $\ub:\R_+\times \R^3 \rightarrow \R^2$ solves the following two-dimensional incompressible Euler system (defining $\overline{\bP_0}$ as the projector onto four-component functions with zero horizontal divergence, the last two components being zero, see below for details):
\begin{equation}
 \begin{cases}
  \d_t \ub +\ub \cdot \n_h \ub = -\n_h \bar{p},\\
  \div_h \ub=0,\\
  \ub_{|t=0}= (\overline{\bP_0}U_0)^h,
 \end{cases}
\label{Eulerh}
\end{equation}
This result was later improved by R. Takada in \cite{T2} for a divergence-free initial data in $H^{s+4}(\R^3)$ ($s\geq 3$) with precise statement of the global existence of the solution $\Ue$ and an explicit convergence rate ($q\in[4,\infty[$):
$$
\|\Ue-(\ub,0,0)\|_{L_T^q W^{1,\infty}(\R^3)} \leq C_{T,q,s,\|U_0\|_{H^{s+4}(\R^3)}} \ee^{\frac1{q}}.
$$
Concerning the viscous system, let us begin with \cite{LT}: S. Lee and R. Takada study, in the particular setting $\nu=\nu'$, the global existence of strong solutions when $\ee$ is small enough and for more regular initial data (that may seem large but satisfy some smallness condition). More precisely, rewording their result with the previous notations, they prove that if $s\in ]\frac12, \frac58]$ there exist $\delta_1,\delta_2>0$ such that for any initial data $U_0$ such that $\overline{\bP_0} U_0\in \dot{H}^\frac12$, $U_{0,osc}\overset{def}{=}(I_d-\overline{\bP_0}) U_0 \in \dot{H}^s$ with:
$$
\|U_{0,osc}\|_{\dot{H}^s} \leq \delta_1 \ee^{-\frac12(s-\frac12)}, \mbox{ and } \|\overline{\bP_0} U_0\|_{\dot{H}^\frac12} \leq \delta_2,
$$
there exists a unique global mild solution $\Ue\in L^4(\dot{W}^{\frac12,3})$.
They also provide another result: if $\|\overline{\bP_0} U_0\|_{\dot{H}^\frac12}$ is sufficiently small, there exists a global solution for small enough $\ee$ (in some sense it is an improvement of the Fujita-Kato theorem as the smallness is not required for the complete initial data but only for what we will call its "stratified part": $\overline{\bP_0} U_0$). As in the other works of their series (see \cite{KLT, IMT}) the main tool are Strichartz estimates obtained through the Littman theorem (see \cite{Litt}). In \cite{LT} the difficulty is that, in contrast to the cases of Systems \eqref{RF} and \eqref{PE}, the Littman theorem cannot be easily applied as the phase function $|\xi_h|/|\xi|$ presents singularities and frequency truncations are necessary to obtain the result.

We point out that in \cite{KLT, IMT, LT} the authors do not study the limit system (these works only focus on global existence), but this is the object of \cite{Scro3} in which S. Scrobogna adapts the ideas of \cite{CDGG, FC2} to state, in the general case ($\nu=\nu'$ is not assumed) that $\Ue$ converges to $(\tv^h,0,0)$ where $\tv^h$ is the unique global solution of the two-component Navier-Stokes system:
\begin{equation}
 \begin{cases}
  \d_t \tv^h +\tv^h \cdot \n_h \tv^h -\nu \D \tv^h & = -\n_h \tpi^0,\\
  \div_h \tv^h=0,\\
  \tv^h_{|t=0}= (\overline{\bP_0} U_0)^h,
\end{cases}
\label{SNS0}
\end{equation}
We emphasize that $\tv^h:\R_+\times \R^3 \rightarrow \R^2$ depends on the full space variable and has two components. More precisely, still rewording the result from \cite{Scro3} using the previous notations, if $U_{0,osc}=(I_d-\overline{\bP_0})U_0\in H^\frac12$ and if $\overline{\bP_0} U_0\in H^1$ (inhomogeneous spaces) then System \eqref{SNS0} admits a unique global strong solution $\tv^h\in \dot{E}^0\cap \dot{E}^1$ and there exists $\ee_0>0$ such that for any $\ee<\ee_0$, System \eqref{Stratif} (with $U_0$ as initial data) admits a unique global strong solution $\Ue \in \dot{E}^\frac12$ (in contrast to \cite{LT}, and as in \cite{FC2}, no smallness condition is required neither from $\overline{\bP_0} U_0$ nor from $U_{0,osc}$). Moreover $\Ue$ converges towards $(\tv^h,0,0)$ in the following sense: if we define $\We$ as the unique global solution of the following linear system (where $\bP$ denotes the classical Leray projector over divergence-free vectorfields),
\begin{equation}
\begin{cases}
 \d_t \We -L\We +\frac1{\ee} \bP \cB \We=\bP\left(\begin{array}{c}0\\0\\-\d_3(-\D_h)^{-1}\div_h(\tv^h\cdot \n_h \tv^h)\\0\end{array}\right),\\
 W_{\ee|t=0}=(I_d-\overline{\bP_0}) U_0,
\end{cases}
 \label{LinScro}
\end{equation}
then $\Ue-(\tv^h,0,0)-\We$ goes to zero in $\dot{E}^\frac12$. Let us precise that in the previous system, the external force term comes from the ideas in \cite{FC2} (also used in \cite{T2}), which can be simply stated as follows: the system satisfied by $\Ue-(\tv^h,0,0)$ features an external force term $G$ which is independent of $\ee$ and prevents any convergence via simple energy methods. But if we "make it oscillate" by putting it (or at least a part of it) as an external force in the previous dispersive linear system then we can absorb a sufficiently large part of $G$ to make the convergence reachable (considering $\Ue-(\tv^h,0,0)-\We$). In the case of \eqref{PE} we pushed further this idea and obtained explicit convergence rates in terms of $\ee$ even for large ill-prepared initial data depending on $\ee$ and, as usual, the results are better when $\nu=\nu'$ because the linearized system is nicer (see \cite{FC3, FCPAA, FCcompl, FCRF}).

The previous result from \cite{Scro3} could be generalized asking less assumptions on the initial data: $\tv_0$ less regular, $\Uoe$ dependent on $\ee$, no low frequency assumption on the initial oscillating part, lower regularity assumptions ($\Uoe\in\dot{H}^\frac12 \cap \dot{H}^{\frac12+\delta}$) and for very large ill-prepared intial data (that is $\Uoeosc$ has size $\ee^{-\alpha}$). This would give a result close to \cite{FCPAA,FCcompl, FCRF} (with explicit convergence rates), but in the present paper we would like to propose \emph{a different kind of generalization}, and answer the following surprising question:
\\

\emph{Why does the limit system \eqref{SNS0} in the result from \cite{Scro3} not depend on $\nu'$ ?}
\\

Before focusing on this question, we would like to end this overview with some very interesting recent results devoted to an intermediate model, the three-scale limit. The models we have presented so far focus on regimes with either only rotation in \eqref{RF}, or only stratification in \eqref{Stratif}, or both of them \emph{but with comparable size} in \eqref{PE} (the non-dispersive case when $F=1$ features tools and results completely different from what is done when $F\neq 1$, see for example \cite{Chemin2, FCF1}).
\\
In \cite{MuScho, MuWei}, the authors completely disconnect the two parameters (namely the Rossby number $\ee$ and the Froude number $\delta$) and study the following \emph{inviscid} system:
$$
\begin{cases}
 \d_t v+v\cdot \n v+\frac1{\ee} e_3\times v+\frac1{\delta} \rho e_3 & =-\n \Phi,\\
 \d_t \rho+v\cdot \n \rho+\frac1{\delta} v_3 & =0,\\
 \div v=0.
\end{cases}
$$
There are now obviously two kinds of penalization ($\delta,\ee \rightarrow 0$) and, in \cite{MuScho}, P. Mu and S. Schochet study the case $\mu=\delta/\ee \rightarrow 0$ (stratification-dominant), whereas in \cite{MuWei}, P. Mu and Z. Wei focus on the alternative case $\nu=\ee/\delta \rightarrow 0$ (rotation-dominant). In both cases the authors need suitable dispersion and Strichartz-type estimates, let us state the dispersive estimates from \cite{MuScho}: if $p_\mu(\xi) =\sqrt{\xi_1^2+\xi_2^2+\mu^2\xi_3^2}/|\xi|$ and $\psi$ is some frequency truncation function, there exists $C=C(\psi)$ such that if $\ee$ and $\delta$ are sufficiently small, for all $f\in L^1(\R^3)$ and all $(t,x)\in \R^4$, we have (denoting $\kappa\overset{def}{=} \ee^{\frac15}\delta^\frac45$)
$$
\left| \int_{\R^3} e^{ix\cdot \xi \pm \frac{t}{\delta} p_\mu(\xi)} \psi(\xi) \hat{f}(\xi) d\xi\right| \leq \frac{C}{(1+\frac{|t|}{\kappa})^\frac14} \|f\|_{L^1(\R^3)}.
$$
The authors also prove that the limit is $U^S=(U_H^S, 0,0)$, where $U_H^S$ solves \eqref{Eulerh}, and they also provide a convergence rate: if $k\geq 6$ is an integer, if $q\in[8,\infty[$, then for any divergence-free initial data $U_0\in H^k(\R^3)\cap L^1(\R^3)$ (inhomogeneous Sobolev space) and for any $T>0$ there exists $\kappa_0$ and $\mu_0$ such that for any $\ee,\delta>0$ for which $\kappa\leq \kappa_0$ and $\mu\leq \mu_0$, there exists a unique solution in $\cC_T([0,T],H^k(\R^3)) \cap \cC^1 ([0,T],H^{k-1}(\R^3))$ and for any $p\in \N\cap[0,k-5]$,
$$
\|U-U^S\|_{L^q([0,T],W^{k-5-p,\infty})} \leq C(\kappa^\frac1{q}+\mu^\frac{3+p}{k+3/2} \sqrt{-\ln \mu}).
$$
In \cite{MuWei}, similar dispersive estimates are obtained, for the phase $p_\nu(\xi)=\sqrt{\nu^2(\xi_1^2+\xi_2^2)+\xi_3^2}/|\xi|$ and now with $\kappa=\ee^{\frac23}\delta^\frac13$. The authors show that when the initial data is taken as follows  $U_0(x)=(u_0^R(x_h),0,0)+w_0(x)$, the limit is $U^R(x)=(U_h^R(x_h),0,\rho^R(x))$, where $U_h^R$ solves the 2D-Euler system with 3 components (which is to be related to the case of \eqref{RF}, see \cite{CDGG, FCRF} and not to be confused with the 3D-Euler system with two components \eqref{Eulerh}) and $\rho^R$ solves the following 2D-transport equation (in the full space variable $x\in \R^3$):
$$
\begin{cases}
 \d_t \rho^R+U_h^R\cdot \n_h\rho^R=0,\\
 \rho^R_{|t=0}=w_0^4.
\end{cases}
$$
If the initial data is in $H^{s+2}(\R^2)+H^s(\R^3)$ (with $s\geq 3$), the authors also provide a convergence rate involving powers of $\kappa$ and $\nu$ (and some radius $\sigma$).

Let us end this survey with \cite{JuMu} where the authors focus on the anisotropic case ($\nu>0$ and $\nu'=0$) in the torus. In this stratification-dominant case they obtain that, for well-prepared initial data (i.-e. with small initial oscillating part) the limit of the global weak solution is $U^S(x)=(U_h^S(x_h),0,\rho^S(x))$ where $U^S$ solves the 2D-Navier-Stokes system and $\rho^S$ is independent of $(t,x_h)$ (but is not explicitely specifed). The authors also provide a convergence rate in terms of $\max(\ee, \delta/\ee)$.

\subsubsection{Statement of the results}

We can now precise the aim of the present article: in both cases $\nu\neq \nu'$ and $\nu=\nu'$, we wish to go in a different direction and \emph{look for ill-prepared initial data that allow a limit really depending on $\nu'$}. We will specify this limit and provide, when possible, global-in-time estimates. As in \cite{CDGG, MuWei}, it will obviously be necessary to get out of the setting of initial data depending on the full space variable and in the next section, which is devoted to the formal obtention of the limit system, we will explain in details how we are lead to modify the initial data.\\

\textbf{Notation}: from now on, the operator $\overline{\bP_0}$ will be denoted as $\bP_2$, and for a given four-component function $f$, we will say that $f_S=\bP_2 \Uoe$ and $f_{osc}=(I_d-\bP_2)\Uoe$ are respectively the "stratified" and "oscillating" parts of $\Uoe$  (see precise definitions below, we recall that $f^h=(f^1, f^2)$).\\

We can now state a simplified version of the results that we prove in this article (the precise results are respectively Theorems \ref{ThLeray} and \ref{ThCV}):

\begin{thm}(Existence)
 \sl{Let $\tTheo\in \dot{B}_{2,1}^{-\frac12}(\R)$ and for all $\ee>0$, let $\Uoe\in L^2(\R^3)$ (divergence-free). Then for all $\ee>0$, System \eqref{Stratif} admits a global weak solution $\Ue\in\dot{E}^0(\R^3)^3\times(\dot{E}^0(\R^3)+\dot{B}^{-\frac12}(\R))$ corresponding to the following initial data:
 $$
 \Uoe(x_1,x_2,x_3)+\left(\begin{array}{c}0\\0\\0\\ \tTheo(x_3)\end{array}\right).
 $$
 Moreover, there exists $C=C_{\nu,\nu',\tTheo}$ such that $\|\Ue\|_{\dot{E}^0(\R^3)+\dot{B}^{-\frac12}(\R)} \leq C(\|\Uoe\|_{L^2}+1)$
 }
\label{Th0exist}
\end{thm}

\begin{thm}(Convergence)
 \sl{For any $\delta>0$, $\Co\geq 1$, $\tTheo\in \dot{B}_{2,1}^{-\frac12}(\R)$, any horizontal divergence-free two-component vectorfield $\tv_0^h\in H^{\frac12+\delta}(\R^3)$ (that is $\div_h \tv_0^h=0$) and, for all $\ee>0$, any divergence-free $\Uoe\in L^2(\R^3)$ ($\Uoe=\UoeS+\Uoeosc$), with:
\begin{equation}
  \begin{cases}
  \vspace{0.1cm}
  \|\tv_0^h\|_{H^{\frac12+\delta}(\R^3)}\leq \Co, \quad \mbox{and} \quad \|\tTheo\|_{\dot{B}_{2,1}^{-\frac12}(\R)}\leq \Co,\\
  \vspace{0.2cm}
  \sup_{\ee>0} \|\Uoe\|_{L^2(\R^3)} \leq \Co, \quad \mbox{and} \quad \|\UoeS-(\tv_0^h,0,0)\|_{L^2(\R^3)} \underset{\ee\rightarrow 0}{\longrightarrow} 0,
 \end{cases}
\end{equation}
then $\Ue$ converges (as $\ee$ goes to zero) to $(\tv^h,0,\tThe)$, where $\tv^h$ and $\tThe$ are the unique global solutions, respectively, of \eqref{SNS0} (with initial data $\tv_0^h$) and the following system:
$$
\begin{cases}
  \d_t \tThe-\nu' \d_3^2 \tThe = 0,\\
  \tThe_{|t=0}= \tTheo,
 \end{cases}
$$
in the following sense:
\begin{itemize}
 \item If $\De=\Ue-(\tv^h,0,\tThe)$ then for any $q\in]2,6[$,
 $$
 \|\De\|_{L_{loc}^2(\R_+, L_{loc}^q(\R^3))}\underset{\ee \rightarrow 0}{\longrightarrow}0.
 $$
 \item For all $q\in]2,6[$, there exist $\ee_1=\ee_1(\nu,\nu',q)>0$, $k_q>0$ and, for all $t\geq 0$, a constant $\mathbb{D}_t=\mathbb{D}_{t,\delta,\nu,\nu',q,\Co}$ such that for all $\ee\in]0, \ee_1]$:
\begin{equation}
\|\Dosc\|_{L_t^2 L^q} =\|(I_d-\bP_2)\De\|_{L_t^2 L^q} \leq \mathbb{D}_t \ee^{k_q}.
 \end{equation}
\item Moreover, when $\nu=\nu'$, the previous estimates can be upgraded into the following global-in-time estimates: there exists a constant $C=C_{\nu,\delta,\Co}>0$ such that, for any $\ee>0$,
$$
\|\Dosc\|_{\tilde{L}^\frac43 \dot{B}_{8, 2}^0 +\tilde{L}^1 \dot{B}_{8, 2}^0} =\|(I_d-\bP_2)\De\|_{\tilde{L}^\frac43 \dot{B}_{8, 2}^0 +\tilde{L}^1 \dot{B}_{8, 2}^0} \leq C \ee^\frac3{16}.
$$
 \end{itemize}
}
\label{Th0CV}
\end{thm}

\subsubsection{Relation with the classical Boussinesq system}

We emphasize that System \eqref{Stratif} is related to the following well-known Boussinesq system:
\begin{equation}
 \begin{cases}
 \d_t v +v\cdot \n v -\nu \D v+\kappa^2 \rho e_3= -\n P,\\
 \d_t \rho + v\cdot \n \rho -\nu' \D \rho =0,\\
 \mbox{div }v=0.
\end{cases}
\label{Bo}
\end{equation}
First, it is easy to see that for any functions $\bar{\rho},\bar{p}:\R_+\times \R \rightarrow \R$ satisfying:
$$
(\d_t-\nu' \d_3^2) \bar{\rho}(t,x_3)=0, \mbox{ and}\quad \d_3 \bar{p}(t,x_3)=-\kappa^2 \bar{\rho}(t,x_3),
$$
then $(\bar{V},\bar{P})$ is an explicit solution of the Boussinesq system \eqref{Bo}, where we have posed:
$$
\bar{V}(t,x)=\left(\begin{array}{c}0\\0\\0\\ \bar{\rho}(t,x_3)\end{array}\right), \quad \bar{P}(t,x)=\bar{p}(t,x_3).
$$
This solution is obviously vertically stratified, and focusing on solutions that are perturbations of $\bar{V}$, if we pose:
\begin{equation}
 V(t,x)=\left(\begin{array}{c}v(t,x) \\ \rho(t,x)\end{array}\right) =\left(\begin{array}{c} v(t,x)\\ \bar{\rho}(t,x_3)+\frac{\theta(t,x)}{\mu} \end{array}\right), \quad U(t,x)=\left(\begin{array}{c}v(t,x) \\ \theta(t,x)\end{array}\right),
\end{equation}
then $(V,P)$ solves \eqref{Bo} if, and only if, $(U, P-\bar{P})$ solves the following system:
\begin{equation}
 \begin{cases}
 \d_t v +v\cdot \n v -\nu \D v+\frac{\kappa^2}{\mu} \theta e_3= -\n \left(P-\bar{P}\right),\\
 \d_t \theta + v\cdot \n \theta -\nu' \D \theta+\mu \d_3 \bar{\rho}(t,x_3) v_3 =0,\\
 \mbox{div }v=0.
\end{cases}
\end{equation}
With a view to connect Systems \eqref{Bo} and \eqref{Stratif}, if we ask in addition that:
$$
\mu \d_3 \bar{\rho}(t,x_3) = -\frac{\kappa^2}{\mu} =-\frac1{\ee},
$$
we end-up with $\mu=\ee \kappa^2$ and
$$
\begin{cases}
 \bar{\rho}_\ee(x)=\bar{\rho}_{0,\ee}-\frac{x_3}{\ee^2 \kappa^2},\\
 \bar{P}_\ee(x)=\bar{P}_{0,\ee}-\kappa^2 \bar{\rho}_{0,\ee} x_3+\frac{x_3^2}{2\ee^2},
\end{cases}
$$
which corresponds (when $\bar{\rho}_{0,\ee}=0$) to the explicit solution $V_\lambda=(0,0,0,\lambda^2 x_3)$ mentioned in \cite{Wid} with  $\lambda=\frac1{\ee\kappa}$ (in \cite{Wid} $\cB$ is replaced by $-\cB$, see also \cite{LT, GenTa}). So that, introducing:
\begin{equation}
  \bar{V}_\ee(x)=\left(\begin{array}{c}0 \\ \bar{\rho}_\ee(x)\end{array}\right) =\left(\begin{array}{c}0\\0\\0\\ \bar{\rho}_{0,\ee}-\frac{x_3}{\ee^2\kappa^2}\end{array}\right),
\label{Explicitsol}
  \end{equation}
then $(\bar{V}_\ee, \bar{P}_\ee)$ is a (stationary) vertically stratified solution of \eqref{Bo} and $(V_\ee, P_\ee)$ solves \eqref{Bo} if, and only if, $(\Ue, \Phie)$ solves \eqref{Stratif}, where we have denoted:
\begin{equation}
 V_\ee=\left(\begin{array}{c}\ve \\ \rho_\ee\end{array}\right) =\left(\begin{array}{c}\ve\\ \bar{\rho}_\ee +\frac{\Thee}{\ee \kappa^2}\end{array}\right), \quad \Ue=\left(\begin{array}{c}\ve \\ \Thee\end{array}\right), \quad \frac1{\ee}\Phie= P_\ee-\bar{P}_\ee.
 \label{ChgtvarBouss}
\end{equation}
Put differently, aside from its own geophysical interest, studying \eqref{Stratif} provides solutions for the Boussinesq system \eqref{Bo} which are perturbations of the explicit solution $(\bar{V}_\ee, \bar{p}_\ee)$.
\\

Thanks to the change of variables described in \eqref{ChgtvarBouss} the previous theorems can be rewritten as a global existence result and asymptotic expansion for the Boussinesq system near the previous explicit stationary solution (see \eqref{Explicitsol}) and for non-conventional vertically stratified initial data:
\begin{thm}
 \sl{With the previous assumptions, for any $\ee>0$, there exists a weak global solution $V_\ee=(v_\ee,\rho_\ee)$ of the Boussinesq system corresponding to the following initial data (the last term is $\Uoe$ with a scaling on its last component):
 $$
 V_\ee|_{t=0}=\left(\begin{array}{c}0\\0\\0\\ \bar{\rho}_{0,\ee}-\frac{x_3}{\ee^2 \kappa^2} \end{array}\right) +\left(\begin{array}{c}0\\0\\0\\ \frac{\tTheo(x_3)}{\ee \kappa^2}\end{array}\right) +\left(\begin{array}{c} v_{0,\ee}(x) \vspace{0.2cm}\\ \frac{\theta_{0,\ee}(x)}{\ee \kappa^2}\end{array}\right).
 $$
Moreover, we have an asymptotic expansion of this solution $V_\ee$ when $\ee$ goes to zero: there exists a four-component function $\De$ such that for any $q\in]2,6[$,
 $$
 \|\De\|_{L_{loc}^2(\R_+, L_{loc}^q(\R^3))}\underset{\ee \rightarrow 0}{\longrightarrow}0,
 $$
and
$$
V_\ee(t,x)=\left(\begin{array}{c}\De^h(t,x) +\tv^h(t,x) \vspace{0.1cm}\\ \De^3(t,x)\\  \bar{\rho}_{0,\ee} -\frac{x_3}{\ee^2 \kappa^2} +\frac{\tThe(t,x_3) +\De^4(t,x)}{\ee \kappa^2}\end{array}\right),
$$
which means that in some sense,
$$
V_\ee(t,x) \underset{\ee \rightarrow 0}{\sim} \left(\begin{array}{c}0\\0\\0\\ \bar{\rho}_{0,\ee}-\frac{x_3}{\ee^2 \kappa^2} \end{array}\right) +\left(\begin{array}{c}0\\0\\0\\ \frac{\tThe(t,x_3)}{\ee \kappa^2}\end{array}\right) +\left(\begin{array}{c}\tv^h(t,x)\\ 0\\ 0\end{array}\right).
$$
}
 \label{ThBouss}
\end{thm}
\begin{rem}
 \sl{The parameters $\bar{\rho}_{0,\ee}, \bar{P}_{0,\ee}$ and $\kappa$ are free (as they do not appear in \eqref{Stratif}): we can make them depend on $\ee$ the way we wish, and choose for instance $\bar{\rho}_{0,\ee}=\bar{\rho}_0 \ee^{-2}$, $\kappa=\ee^{-1}$ or $\ee^{-\frac12}$.}
\end{rem}

\subsubsection{Outline of the article}

The present article is structured as follows: in the next section, we will formally obtain the limit system and the natural decomposition it induces (namely the stratified part and the oscillating part), and study their properties. This decomposition determines what means to be well- or ill-prepared. The limit systems are studied in the end of Section 2.

Then we will decompose the solutions of \eqref{Stratif} and obtain a more classical system for which we have to adapt the classical Leray theorem in the same spirit as what was done in \cite{CDGG, CDGG2, CDGGbook} for the rotating fluids system. A precise existence theorem for weak solutions (for more general initial data depending on $\ee$) is given in the beginning of Section 3.

Using the Strichartz estimates that we put in Appendix 1, we can show the oscillating part $(I_d-\bP_2)\De$ goes to zero. From this, we obtain the convergence of the stratified part $\bP_2 \De$. The precise convergence theorem is stated in Section 4. This theorem also features far better convergence rate in the case when $\nu=\nu'$ (using improved Strichartz estimates also proved in Appendix 1). The proof of the Strichartz estimates relies on a crucial technical result which is proved in Appendix 2.

The study of these asymptotics for strong solutions (i.-e. obtained with the Fujita-Kato theorem) with large regular ill-prepared initial data is dealt in the companion work \cite{FCStratif2}.

\section{The limit system}

In this section we will obtain the formal limit system when $\ee$ goes to zero. As pointed out in the introduction, a first attempt was made in \cite{Scro3} only with initial data depending on the full space variables. In our article thanks to a decomposition of the geopotential and to formal arguments in the spirit of \cite{FC1}, we will obtain a much more general limit system (see \eqref{SNS2}). Studying its solutions and structure will help reformulate \eqref{Stratif} in a more suitable way.

\subsection{Formal argument}

Taking the divergence of the velocity part of \eqref{Stratif} leads to:
$$
\D \Phie =-\d_3 \Thee -\ee\sum_{i=1}^3 \d_i \ve\cdot \nabla \ve^i.
$$
As $\ve$ is divergence-free we have:
$$
\sum_{i=1}^3 \d_i \ve\cdot \nabla \ve^i= \div (\ve \cdot \n \ve) =\sum_{i=1}^3 \d_i (\ve\cdot \nabla \ve^i)= \sum_{i,j=1}^3 \d_i \d_j (\ve^i \ve^j),
$$
from which we decompose the geopotential into $\Phie\overset{def}{=} \Pe^1 +\ee \Pe^0$, where:

\begin{equation}
\begin{cases}
 \Pe^1= -\D^{-1} \d_3 \Thee,\\
 \Pe^0= -\sum_{i,j=1}^3 \d_i \d_j \D^{-1} (\ve^i \ve^j).
\end{cases}
 \label{defpres}
\end{equation}
We can then write in extension System \eqref{Stratif} as follows:
\begin{equation}
\begin{cases}
\d_t \ve^1 +\ve\cdot \n \ve^1 -\nu \D \ve^1 & = -\d_1 \Pe^0-\frac{1}{\ee} \d_1 \Pe^1,\\
\d_t \ve^2 +\ve\cdot \n \ve^2 -\nu \D \ve^2 & = -\d_2 \Pe^0-\frac{1}{\ee} \d_2 \Pe^1,\\
\d_t \ve^3 +\ve\cdot \n \ve^3 -\nu \D \ve^3 & = -\d_3 \Pe^0-\frac{1}{\ee} (\d_3 \Pe^1+\Thee),\\
\d_t \Thee +\ve\cdot \n \Thee -\nu' \D \Thee & = \frac{1}{\ee} \ve^3,\\
\div \ve & =0.
\end{cases}
\label{Strext}
\end{equation}
If we assume that $(\ve,\Thee, \Pe^0, \Pe^1) \underset{\ee \rightarrow 0}{\longrightarrow} (\tv, \tThe, \tP^0, \tP^1)$ in a sufficiently strong way such that the convergence also occurs for the derivatives and nonlinear terms, then taking into account the penalized terms from the right-hand-side, necessarily:
\begin{equation}
\begin{cases}
  \d_1 \tP^1=\d_2 \tP^1=0,\\
 \d_3 \tP^1+\tThe=0,\\
 \tv^3=0,
\end{cases}
\mbox{ which means that }
\begin{cases}
\tP^1 \mbox{ and }\tThe=-\d_3 \tP^1 \mbox{ only depend on } x_3,\\
 \tv^3=0.
\end{cases}
\end{equation}
Moreover, as $\ve$ is divergence-free and as $\tv^3=0$, we obtain:
\begin{equation}
 \div_h \tv^h \overset{def}{=} \d_1 \tv^1+\d_2 \tv^2 =0 \quad \mbox{where we define } \tv^h \overset{def}{=} (\tv^1, \tv^2).
\end{equation}
Next, if we ask in addition that:
\begin{equation}
 \begin{cases}
-\frac{1}{\ee} \d_1 \Pe^1 & \underset{\ee \rightarrow 0}{\longrightarrow} \tX,\\
-\frac{1}{\ee} \d_2 \Pe^1 & \underset{\ee \rightarrow 0}{\longrightarrow} \tY,\\
-\frac{1}{\ee} (\d_3 \Pe^1+\Thee) & \underset{\ee \rightarrow 0}{\longrightarrow} \tZ,\\
\frac{1}{\ee} \ve^3 & \underset{\ee \rightarrow 0}{\longrightarrow} \tT,
\end{cases}
\end{equation}
then we end up with the following limit system:
\begin{equation}
\begin{cases}
\d_t \tv^1 +\tv^h \cdot \n_h \tv^1 -\nu \D \tv^1 & = -\d_1 \tP^0+\tX,\\
\d_t \tv^2 +\tv^h \cdot \n_h \tv^2 -\nu \D \tv^2 & = -\d_2 \tP^0+\tY,\\
0 & = -\d_3 \tP^0+\tZ,\\
\d_t \tThe -\nu' \d_3^2 \tThe & = \tT,\\
\div_h \tv^h & =0.
\end{cases}
\label{SNS1}
\end{equation}
Taking the limit in the expressions of $\Pe^0$ and $\Pe^1$ we obtain that:
\begin{equation}
\begin{cases}
\vspace{0.2cm}
  \D \tP^1=\d_3^2 \tP^1 =-\d_3 \tThe \mbox{, which brings nothing new},\\
 \D \tP^0=-\sum_{i=1}^2 \d_i \tv^h \cdot \nabla_h \tv^i= -\sum_{i=1}^2 \d_i (\tv^h \cdot \nabla_h \tv^i)= -\sum_{i,j=1}^2 \d_i \d_j (\tv^i \tv^j),
\end{cases}
\label{P0P1}
\end{equation}
As the horizontal divergence $\div_h \tv^h$ is zero, similarly as for \eqref{Stratif}, computing for \eqref{SNS1} $\d_1 (\tilde{\mbox{line }1})+\d_2(\tilde{\mbox{line }2})+\d_3 (\tilde{\mbox{line }3})$, we get:
\begin{equation}
 \d_1 \tX+ \d_2 \tY+ \d_3 \tZ=0.
\end{equation}
Moreover, computing on \eqref{Stratif} $\d_1 (\mbox{line }2)-\d_2 (\mbox{line }1)$, we get that, introducing $\ome \overset{def}{=} \d_1 \ve^2-\d_2 \ve^1$:
\begin{equation}
 \d_t \ome + \ve\cdot \n \ome + \div_h \ve \cdot \ome +\d_1 \ve^3\cdot \d_3 \ve^2 -\d_2 \ve^3\cdot \d_3 \ve^1 -\nu \D \ome =0.
\end{equation}
Performing $\ee \rightarrow 0$ in the previous equation, using that $\tv^3=0=\div_h \tv^h$, and  defining $\tom \overset{def}{=} \d_1 \tv^2-\d_2 \tv^1$, we get that:
\begin{equation}
 \d_t \tom + \tv^h \cdot \n_h \tom -\nu \D \tom =0.
 \label{eqomtilde1}
\end{equation}
On the other hand, computing on \eqref{SNS1} $\d_1 (\tilde{\mbox{line }2})-\d_2 (\tilde{\mbox{line }1})$, we get that:
\begin{equation}
 \d_t \tom + \tv^h \cdot \n_h \tom -\nu \D \tom =\d_1 \tY- \d_2 \tX,
 \label{eqomtilde2}
\end{equation}
so that, identifying the right-hand side from \eqref{eqomtilde1} and \eqref{eqomtilde2}, we obtain that:
\begin{equation}
 \d_1 \tY- \d_2 \tX=0.
\end{equation}
Next, using the third line from \eqref{SNS1} and formally solving in the Fourier variable the system:
$$
\begin{cases}
\d_1 \tX+ \d_2 \tY & =-\d_3 \tZ =-\d_3^2 \tP^0,\\
- \d_2 \tX +\d_1 \tY & =0,
\end{cases}
$$
we end up with:
\begin{equation}
 (\tX, \tY)=- \n_h\d_3^2 \D_h^{-1} \tP^0,
\end{equation}
so that:
\begin{equation}
 -\n_h \tP^0+(\tX, \tY)=-\n_h \D_h^{-1} (\D_h \tP^0 +\d_3^2 \tP^0) =-\n_h \D_h^{-1} \D \tP^0 =-\n_h \tpi^0,
\end{equation}
where from \eqref{P0P1} we have that:
\begin{equation}
 \tpi^0 =\D_h^{-1} \D \tP^0 =-\sum_{i,j=1}^2 \D_h^{-1} \d_i \d_j (\tv^i \tv^j).
\label{defq}
 \end{equation}
Gathering the previous informations, we can now write the formal limit system \eqref{SNS1} as follows, still denoting $\tU=(\tv^1, \tv^2, 0, \tOm)=(\tv^h, 0, \tOm)$:
\begin{equation}
    \begin{cases}
  \d_t \tv^h +\tv^h \cdot \n_h \tv^h -\nu \D \tv^h & = -\n_h \tpi^0,\\
  \div_h \tv^h=0,\\
  \end{cases}
 \quad \mbox{and} \quad \d_t \tThe-\nu' \d_3^2 \tThe = \tT,
 \label{SNS2}
\end{equation}
where $\tv^3=0$, $\tThe=-\d_3 \tP^1$, $\tZ=\d_3 \tP^0$, and with $\tP^1, \tThe, \tT$ only depending on $(t,x_3)$.
\begin{rem}
 \sl{
 \begin{enumerate} \item Solving in the Fourier variable the following system:
$$
\begin{cases}
 \d_1 \tv^1 +\d_2 \tv^2 & =0,\\
 -\d_2 \tv^1 +\d_1 \tv^2  & =\tom,
\end{cases}
$$
we get that $\tv^h= (-\d_2 \D_h^{-1} \tom, \d_1 \D_h^{-1} \tom) =\n_h^\perp \D_h^{-1} \tom$. Therefore the first system can be rewritten in terms of the vorticity $\tom=\omega(\tv)=\d_1 \tv^2-\d_2 \tv^1$ as follows:
\begin{equation}
\begin{cases}
 \d_t \tom + \tv^h \cdot \n_h \tom -\nu \D \tom =0,\\
 \tv^h=\n_h^\perp \D_h^{-1} \tom.
\end{cases}
 \label{eqomtilde3}
\end{equation}
\item We emphasize that there is no stretching term $\tom \cdot \n_h \tv^h$ in \eqref{eqomtilde3}, like in the $2D$-Navier-Stokes case (full limit of the rotating fluids system, see \cite{CDGG}), or the Quasi-geotrophic system $(\tilde{QG})$ (limit of the Primitive system, see \cite{FC1, FCPAA}) and unlike in the $3D$-Navier-Stokes case.
\end{enumerate}
\label{Remomtilde}
}
\end{rem}

\begin{rem}
\sl{
\begin{enumerate}
 \item At this stage, nothing can help us precise the term $\tT$, we only know that $\tT=\underset{\ee \rightarrow 0}{lim} \frac{\ve^3}{\ee}$.
 \item Moreover, nothing allows us to make a direct link between $\tT(0, .)$ and the initial velocity such as for instance:
 $$
 \tT(0,x_3)=\underset{\ee \rightarrow 0}{lim} \frac{\ve^3(0,x)}{\ee}.
 $$
 Indeed, we easily obtain in the case of the Primitive system $(PE_\ee)$ (see \cite{FC1}) that $\frac{\ve^3}{\ee}$ has a limit when $\ee$ goes to zero:
 $$
 \frac{\ve^3}{\ee F} \underset{\ee \rightarrow 0}{\longrightarrow} (\nu-\nu') \D \D_F^{-1} \tilde{\theta_{QG}}+ F[\d_3 \D_F^{-1}, \tilde{v_{QG}}^h\cdot \n_h] \tOm,
 $$
 whereas if, for instance, the initial data does not depend on $\ee$, then for all $x\in \R^3$ such that $\ve^3(x)\neq 0$,
 $$
 \frac{\ve^3(0,x)}{\ee F} =\frac{v_0^3(x)}{\ee F} \underset{\ee \rightarrow 0}{\longrightarrow} \pm \infty.
 $$
 \item It seems that we have better considering $\tT$ as a parameter, and in this article we will search for initial data corresponding to the case $\tT=0$.
 \item On the contrary, nothing forces $\tThe$ to be zero, and this term will be the key to obtain new asymptotics, in the spirit of \cite{CDGG, CDGGbook}: for the rotating fluids system, when the initial data is in $L^2(\R^3)$ or in $\dot{H}^s(\R^3)$ the only limit we can hope for is zero. But when Chemin, Desjardins, Gallagher and Grenier consider as an initial data the sum of a function of the previous type with some  $\bar{u}\in L^2(\R^2_h)$, then richer asymptotics are at hand and they manage to reach as limit system the $2D$-Navier-Stokes system with three components with $\bar{u}$ as an initial data (as described by physicists in the Taylor-Proudman theorem). We show that similar phenomena will occur for the strongly stratified Boussinesq system.
\end{enumerate}
}
\end{rem}

\subsection{The Stratified/osc structure, final form of the limit system}

Thanks to the previous observations we will consider initial data satisfying:
$$
{\Thee}_{|t=0}(x)=\Theeo (x) +\tTheeo (x_3),
$$
so that our complete initial data is:
$$
{\Ue}_{|t=0}(x) =\Uoe(x) +(0,0,0,\tTheeo (x_3)).
$$
Moreover, the structure of the formal limit system suggests to introduce the following operators:
\begin{defi}
 \sl{For a $R^4$-valued function, we introduce the following quantity, that we will call its vorticity:
 $$
 \omega(f)=\d_1 f^2 -\d_2 f^1.
 $$
 From this we define the \emph{stratified} and \emph{oscillating (or oscillatory) parts} of f, respectively denoted as $f_S$ and $f_{osc}$, according to:
 \begin{equation}
 f_S= \left(
 \begin{array}{c}
  \n_h^\perp \D_h^{-1} \omega(f)\\
  0\\
  0
 \end{array}\right)
=\left(
 \begin{array}{c}
  -\d_2 \D_h^{-1} \omega(f)\\
  \d_1 \D_h^{-1} \omega(f)\\
  0\\
  0
 \end{array}\right),
 \label{defS}
 \end{equation}
 and, denoting $\div_h f^h \overset{def}{=} \d_1 f^1 + \d_2 f^2$,
\begin{equation}
 f_{osc}=f-f_{S}= \left(
 \begin{array}{c}
  \n_h \D_h^{-1} \div_h f^h\\
  f^3\\
  f^4
 \end{array}\right)
 =\left(
 \begin{array}{c}
  \d_1 \D_h^{-1} \div_h f^h\\
  \d_2 \D_h^{-1} \div_h f^h\\
  f^3\\
  f^4
 \end{array}\right).
 \label{defosc}
 \end{equation}
}
\end{defi}
\begin{rem}
 \sl{For a $\R^2$-valued function $f=(f^1,f^2)=f^h$, we could introduce $\bP_2^h$ and $f_S= \bP_2^h f= \n_h^\perp \D_h^{-1} \omega(f)$ (this corresponds to $\overline{\bP_0}$ in \eqref{Eulerh}), but with a slight notational abuse, we may also denote $f_S=\bP_2 f$} and $f_{osc}=f-f_{S}=\n_h \D_h^{-1} \div_h f$.
\end{rem}
Now we can completely precise the initial data and limit system that we will consider in this article:
\begin{equation}
 {\Ue}_{|t=0}(x) =\Uoe(x)
 +\left(\begin{array}{c}
  0\\0\\0\\\tTheeo (x_3)
 \end{array}\right)
 =\UoeS(x) +\Uoeosc(x)
 +\left(\begin{array}{c}
  0\\0\\0\\\tTheeo (x_3)
 \end{array}\right).
\end{equation}
And we will denote:
$$
\Uoe=\left(\begin{array}{c}
  \voe^1\\ \voe^2\\ \voe^3\\ \Theeo
 \end{array}\right)
 \quad \mbox{and} \quad
 \Ue=\left(\begin{array}{c}
  \ve^1\\ \ve^2\\ \ve^3\\ \Thee
 \end{array}\right)
 =\left(\begin{array}{c}
  \ve \\ \Thee
 \end{array}\right)
=\left(\begin{array}{c}
  \ve^h\\ \ve^3\\ \Thee
 \end{array}\right).
 $$
The previous formal study suggests that the initial data in \eqref{SNS2} are related to $\UoeS$ and $\tTheeo$, in the following way:
\begin{equation}
 \begin{cases}
 \UoeS^h (x) \underset{\ee \rightarrow 0}{\longrightarrow} \tvo^h (x), \quad \mbox{or equivalently } \UoeS (x)=\bP_2 \Uoe (x) \underset{\ee \rightarrow 0}{\longrightarrow} (\tvo^h (x),0,0),\\
 \tTheeo (x_3) \underset{\ee \rightarrow 0}{\longrightarrow} \tTheo (x_3),
\end{cases}
\label{Condinit}
\end{equation}
so that our complete limit system is composed by the following two systems:
\begin{equation}
    \begin{cases}
  \d_t \tv^h +\tv^h \cdot \n_h \tv^h -\nu \D \tv^h & = -\n_h \tpi^0,\\
  \div_h \tv^h=0,\\
  \tv^h_{|t=0}= \tvo^h,
  \end{cases}
 \label{SNS3}
\end{equation}
and
\begin{equation}
 \begin{cases}
  \d_t \tThe-\nu' \d_3^2 \tThe = 0,\\
  \tThe_{|t=0}= \tTheo.
 \end{cases}
 \label{SNS4}
\end{equation}

\begin{rem}
 \sl{We could simply choose $\tTheeo = \tTheo$, but for more generality we made the previous choice. The same can be done to generalize a little the result in \cite{FCRF} concerning rotating fluids.}
\end{rem}
The following proposition states properties induced by the stratified/oscillating structure.
\begin{prop}
 \sl{With the notations from \eqref{defS} and \eqref{defosc}, there exist two pseudodifferential operators of order zero $\mathcal{P}$ and $\mathcal{Q}$ such that for any $f$,
$$
f_{S}=\mathcal{Q} f,\quad \mbox{and} \quad f_{osc}=\mathcal{P} f.
$$
These operators satisfy:
\begin{enumerate}
 \item $\mathcal{Q}=\bP_2$ and $\mathcal{P}=I_d-\bP_2$ (where the operators $\bP_k$ are defined in \eqref{defP2} and \eqref{defPk}).
 \item For any $s\in \R$, we have $((I_d-\bP_2) f| \bP_2 f)_{H^s/\dot{H}^s} =0=(\cB f| \bP_2 f)_{H^s/\dot{H}^s}$ (when defined).
 \item $(I_d-\bP_2) f=f \Longleftrightarrow \bP_2 f =0 \Longleftrightarrow \omega(f)=0$.
 \item $(I_d-\bP_2) f=0 \Longleftrightarrow \bP_2 f =f \Longleftrightarrow f^3=f^4=0 \mbox{ and } \div_h f=0 \Longleftrightarrow$ there exists a scalar function $\phi$ such that $f=(-\d_2 \phi, \d_1 \phi,0,0) =(\n_h^\perp \phi, 0,0)$. Such a vector field is obviously divergence free (and horizontal divergence-free) and we will say that it is \emph{stratified}. It also satisfies $f=(f^h,0,0)$.
\item If $f$ is divergence-free, so is $(I_d-\bP_2) f$.
 \item $\cB \bP_2 f=0$ (in $\R^4$).
 \item $\bP_2 \bP=\bP \bP_2=\bP_2$ and $\bP_2 (I_d-\bP)=(I_d-\bP) \bP_2=0$ (in particular $\bP_2 (\n q,0)=0$).
 \item If f is a divergence-free vector field, then (we recall that we denote $f\cdot \n f=\sum_{i=1}^3 f^i \d_i f$)
 $$
 \omega(f\cdot \n f)=-\d_3 f^3\cdot \omega(f)+\d_1 f^3 \cdot \d_3 f^2 -\d_2 f^3 \cdot \d_3 f^1 +f\cdot \omega (f).
 $$
 \item If f is a stratified vector field, then $\omega(f\cdot \n f)=f\cdot \omega (f)$.
\end{enumerate}}
\label{PropSosc}
\end{prop}
The proof of this proposition is similar to what is done in \cite{FC1, FC2, FCPAA}.
\begin{rem}
\sl{\begin{enumerate}
     \item Let us recall that, in the case of the Primitive system, the decomposition QG/osc (with its corresponding $\cQ/\cP$) and $\bP_2/\bP_{3,4}$ do not coincide except when $\nu=\nu'$, so that for instance the oscillating part still has a "$\bP_2$" part (but we show it is very small in terms of $\ee$). In the present article, we have $\cQ=\bP_2$ in the general case (but it is only when $\nu=\nu'$ that $\bP_3$ and $\bP_4$ are orthogonal projectors of norm $1$).
     \item We chose to use the denomination \emph{Stratified} rather than \emph{quasi-geostrophic}.
    \end{enumerate}
}
\end{rem}

\subsection{Reformulation of the different systems}

In order to properly study the solutions of \eqref{Stratif} with such unusual initial data, as in \cite{CDGG, CDGGbook} in the case of the rotating fluids, we need to rewrite System \eqref{Stratif} in a form where functions only depending on $x_3$ do not appear in the initial data anymore. Doing this will move these functions in the transport terms which explains why a little adaptation of the proofs of the classical existence results such as the Leray or Fujita-Kato theorems is needed.
\begin{rem}
 \sl{When studying $(PE_\ee)$ or in \cite{Scro3}, this was not necessary because the systems only feature functions of $x\in\R^3$ and adapting the classical existence theorems required no effort.}
\end{rem}
Let us first merge Systems \eqref{SNS3} and \eqref{SNS4}. Introducing $\tv^3=0$ and $\tP^1$ such that $\tThe(x_3)=-\d_3 \tP^1(x_3)$, we can rewrite the limit systems as follows:
\begin{equation}
\begin{cases}
\d_t \tv^1 +\tv \cdot \n \tv^1 -\nu \D \tv^1 & = -\d_1 \tpi^0,\\
\d_t \tv^2 +\tv \cdot \n \tv^2 -\nu \D \tv^2 & = -\d_2 \tpi^0,\\
\d_t \tv^3 +\tv \cdot \n \tv^3 -\nu \D \tv^3 +\frac1{\ee} \tThe & = -\frac1{\ee}\d_3 \tP^1,\\
\d_t \tThe  +\tv \cdot \n \tThe -\nu' \D \tThe -\frac1{\ee} \tv^3 & = 0,\\
\div \tv & =0.
\end{cases}
\label{SNS5}
\end{equation}
If we set $\tU\overset{def}{=}(\tv^h,0,\tThe)$, this system can be rewritten into:
$$
\d_t \tU +\tv \cdot \n \tU -L \tU + \frac1{\ee} \cB \tU = -\left( \begin{array}{c}\n_h \tpi^0 \\0 \\0
\end{array}\right)
-\frac1{\ee} \left(\begin{array}{c}\n \tP^1 \\0 \end{array} \right).
$$
Denoting as $\bP$ the orthogonal Leray projector onto divergence-free vectorfields, we introduce $\tG=\bP(\n_h \tpi^0,0,0)$ which satisfies the following result.
\begin{prop}
 \sl{With the previous notations, there exists some $\tg$ such that
 \begin{equation}
   \left( \begin{array}{c}\n_h \tpi^0 \\0 \\0
\end{array}\right)=\tG+\left(\begin{array}{c}
  \n \tg\\0 \end{array}\right) \quad \mbox{with} \quad
  \tG= \bP \left(\begin{array}{c}
  \d_1 \tpi^0\\ \d_2 \tpi^0 \\0\\0 \end{array}\right)
  =\left(\begin{array}{c}
  \d_1 \d_3^2 \D^{-1} \D_h^{-1} \tq_0\\\d_2 \d_3^2 \D^{-1} \D_h^{-1} \tq_0 \\-\d_3 \D^{-1} \tq_0\\0 \end{array}\right),
  \label{deftG}
 \end{equation}
 where $\tq_0$ is defined by ($\tpi^0$ has been introduced in \eqref{defq})
 $$
 \tq_0 =-\D_h \tpi^0 =\sum_{i,j=1}^2 \d_i \d_j (\tv^i \tv^j)=\sum_{i=1}^2 \d_i (\tv^h\cdot \n_h \tv^i).
 $$
 Moreover, it is obvious that $\omega(\tG)=0=\div \tG$ and $\bP_2 \tG=0$.
 }
 \label{PropG}
\end{prop}
\textbf{Proof:} Just compute $\tG$ using the fact that $\bP=I_d-\n \D^{-1} \div$. $\blacksquare$
\begin{rem}
 \sl{Adding some gradient, we can also write that:
 $$
 \tG=\bP\left(\begin{array}{c}
 0\\0\\ \d_3 \D_h^{-1} \tq_0\\0 \end{array}\right),
 $$
which appears in the external force term from the auxiliary system \eqref{LinScro} introduced in \cite{Scro3} when reproducing the method from \cite{FC2}, where the analoguous term features two parts with different behaviours, namely $G=G^b+G^l$.}
\end{rem}
Thanks to this, we can finally reformulate our limit system as follows:
\begin{equation}
\begin{cases}
 \d_t \tU +\tU \cdot \n \tU -L \tU + \frac1{\ee} \cB \tU = -\tG -\left(\begin{array}{c}\n \tg \\0\end{array}\right)-\frac1{\ee} \left(\begin{array}{c}\n \tP^1\\0\end{array}\right),\\
 \div \tv=0,\\
 \tU_{|t=0}=(\tvo^h,0, \tTheo).
\end{cases}
\label{SNS6}
\end{equation}
Our first idea would be to study the system satisfied by $\Ue-\tU$, but as already mentionned, its initial data:
$$
(\Ue-\tU)_{|t=0}(x)=\Uoeosc(x)+ (\bP_2 \Uoe^h(x) -\tvo^h(x), 0, \tTheeo(x_3)-\tTheo(x_3)),
$$
would feature functions depending only on $x_3$. As we wish to do minimal transformations in order to adapt the Leray theorem, we simply define the following function, that will help us to neutralize the $x_3$-only-dependent part:
\begin{equation}
 \tZe=\left(\begin{array}{c}
  0\\0\\0\\ \tKe \end{array}\right), \quad \mbox{where } \tKe \mbox{ solves }
  \begin{cases}
   \d_t \tKe-\nu' \d_3^2 \tKe = 0,\\
  \tThe_{|t=0}= \tTheeo-\tTheo.
  \end{cases}
\end{equation}
\begin{rem}
 Thanks to the estimates from Theorem \ref{ThHeat} and \eqref{Condinit}, $\tKe$ goes to zero as $\ee$ goes to zero.
\end{rem}
Next, introducing $\tP_{\ee}^2$ such that $\tKe=-\d_3 \tP_{\ee}^2$ (both of these functions only depend on $x_3$), $\tZe$ satisfies:
\begin{equation}
\begin{cases}
  \d_t \tZe -L \tZe + \frac1{\ee} \cB \tZe = -\frac1{\ee}(\n \tP_{\ee}^2,0),\\
  \tilde{Z}_{\ee|t=0} =(0,0,0, \tTheeo-\tTheo).
\end{cases}
\label{SystZe}
\end{equation}
Now we can properly rewrite the system we will study, let us finally set:
\begin{equation}
 \De \overset{def}{=} \Ue-\tU-\tZe =\left(\begin{array}{c}
  \ve^h-\tv^h\\\ve^3\\ \Thee-(\tThe+\tKe) \end{array}\right) =\left(\begin{array}{c}
  \ve^h-\tv^h\\\ve^3\\ \Thee-\tThee \end{array}\right)=\Ue-\left(\begin{array}{c}\tv^h\\ 0\\ \tThee \end{array}\right),
  \label{defD}
\end{equation}
where the function $\tThee\overset{def}{=}\tThe+\tKe$ is nothing but the solution of:
\begin{equation}
 \begin{cases}
   \d_t \tThee-\nu' \d_3^2 \tThee = 0,\\
  \tilde{\theta}_{\ee|t=0}= \tTheeo.
  \end{cases}
  \label{SysttThee}
\end{equation}
Rewriting System \eqref{Stratif} as follows,
\begin{equation}
\begin{cases}
\d_t \Ue +\Ue\cdot \n \Ue -L \Ue +\frac{1}{\ee} \cB \Ue=-\left(\begin{array}{c}\n \Pe^0 \\0\end{array}\right)-\frac1{\ee} \left(\begin{array}{c}\n \Pe^1\\0\end{array}\right),\\
{\Ue}_{|t=0}=U_{0,\ee}.
\end{cases}
\label{Stratif2}
\tag{$S_\ee$}
\end{equation}
and substracting Systems \eqref{SNS6} and \eqref{SystZe}, we obtain that $\De$ satisfies:
$$
\d_t \De -L \De + \frac1{\ee} \cB \De = -\left[\De\cdot \n \De +\De \cdot \n \left(\begin{array}{c}
  \tv^h \\0\\ \tThee \end{array}\right) +\tv^h\cdot \n_h \De \right] +\tG -\left(\begin{array}{c}\n \qe \\0 \end{array}\right),
$$
where
$$
\qe=-\tg+\frac1{\ee}(\phi_\ee-\tP^1-\tP_\ee^2) =\Pe^0-\tg +\frac1{\ee}(\Pe^1-\tP^1-\tP_\ee^2).
$$
Reformulating the second term from the right-hand side, we obtain the final form \eqref{StratifD} of the system satisfied by $\De=(\Ve,\He)$:
\begin{equation}
\begin{cases}
 \d_t \De -L \De + \frac1{\ee} \cB \De = -\left[\De\cdot \n \De
 +\left(\begin{array}{c}
  \De \cdot \n \tv^h \\0\\ \De^3\cdot \d_3 \tThee \end{array}\right) +\tv^h\cdot \n_h \De \right] +\tG -\left(\begin{array}{c}\n \qe \\0 \end{array}\right),\\
 \div \Ve=0,\\
 D_{\ee|t=0}=\Uoeosc +\UoeS-(\tvo^h,0,0) =\Uoeosc +(\UoeS^h-\tvo^h,0,0).
\end{cases}
\label{StratifD}
\end{equation}

\subsection{Study of the limit system}

We show in this section that the limit system has global regular solutions which specifies the necessary informations about $\tv^h$ and $\tThee$ required to study System \eqref{StratifD}.

Let us begin with System \eqref{SNS4}, which is only a one-dimensional heat equation (we refer for example to \cite{Dbook}, Section 3.4.1, Lemma 5.10 and Proposition 10.3, see also Definition \ref{deftilde}).

\begin{thm}
 \sl{Let $s\in \R$. For any $\tTheo\in \dot{H}^s(\R)$ (respectively $\tTheo\in \dot{B}_{2,1}^s(\R)$) there exists a unique global solution $\tThe$ of \eqref{SNS4} and for all $t\geq 0$, we have:
 \begin{equation}
  \|\tThe\|_{\tilde{L}_t^\infty \dot{H}^s}^2 +\nu' \|\tThe\|_{L_t^2 \dot{H}^{s+1}}^2 \leq 2 \|\tTheo\|_{\dot{H}^s}^2.
 \end{equation}
\begin{equation}
(\mbox{respectively} \quad \|\tThe\|_{\tilde{L}_t^\infty \dot{B}_{2,1}^s} +\nu' \|\tThe\|_{L_t^1 \dot{B}_{2,1}^{s+2}} \leq \|\tTheo\|_{\dot{B}_{2,1}^s}.)
 \end{equation}
 More generally for $s\in \R$ and $p,r\in[1,\infty]$, there exists a constant $C>0$ such that if $\tTheo\in \dot{B}_{p,r}^s(\R)$ then for all $q\in [1,\infty]$
 \begin{equation}
  \|\tThe\|_{\tilde{L}_t^q \dot{B}_{p,r}^{s+\frac2{q}}}\leq \frac{C}{(\nu')^\frac1{q}} \|\tTheo\|_{\dot{B}_{p,r}^s}.
   \label{estimThetaBspr}
 \end{equation}
 }
 \label{ThHeat}
\end{thm}
\begin{rem}
 \sl{Thanks to this result, the previously defined $\tThe$, $\tKe$ and $\tThee$ are global and satisfy similar estimates.
}
\end{rem}
On the other hand, thanks to \eqref{eqomtilde3} and Remark \ref{Remomtilde}, we observe that System \eqref{SNS3} (which is System (4.3) from \cite{Scro3}) is very close to the quasi-geostrophic system, and  we can easily adapt Theorem 1 from \cite{FCPAA} and state the following theorem which generalizes the results from Section 4 in \cite{Scro3} (as we need less initial regularity).

\begin{thm}
 \sl{Let $\delta>0$ and $\tvo^h\in H^{\frac12+\delta}$ a $\R^2$-valued vectorfield such that $\div_h \tvo^h=0$. Then System \eqref{SNS3} has a unique global solution $\tv^h \in E^{\frac12+\delta}=\dot{E}^0 \cap \dot{E}^{\frac12+\delta}$ and there exists a constant $C=C_{\delta, \nu}>0$ such that for all $t\geq 0$, we have:
 \begin{multline}
  \|\tv^h\|_{L^\infty H^{\frac12+\delta}}^2 +\nu \|\n \tv^h\|_{L^2 H^{\frac12+\delta}}^2 \leq C_{\delta, \nu} \|\tvo^h\|_{H^{\frac12+\delta}}^2 \max (1,\|\tvo^h\|_{H^{\frac12+\delta}}^\frac1{\delta})\\
  \leq C_{\delta, \nu} \max (1,\|\tvo^h\|_{H^{\frac12+\delta}})^{2+\frac1{\delta}},
 \end{multline}
 Moreover, we can also bound the term $\tG$ introduced in \eqref{defD}: for all $s\in[0,\frac12+\delta]$,
 \begin{equation}
  \int_0^\infty \|\tG(\tau)\|_{\dot{H}^s} d\tau \leq C_{\delta, \nu} \max (1,\|\tvo^h\|_{H^{\frac12+\delta}})^{2+\frac1{\delta}}.
 \end{equation}
 }
 \label{ThSNS}
\end{thm}
\textbf{Proof:} Thanks to Remark \ref{Remomtilde} and the Biot-Savart law (which is similar to the one featured in the quasi-geostrophic system), the proof is very close to what we did in Section 2.1 from \cite{FCPAA}. $\blacksquare$

\section{Existence of global weak solutions}

This section is devoted to prove the following result, which is the analoguous of the Leray theorem for \eqref{StratifD} and provides global weak solutions for any $\ee>0$.

\begin{thm} (Existence of Leray weak solutions)
 \sl{For any $\delta>0$ and $\Co\geq 1$, let $\tv_0^h\in H^{\frac12+\delta}(\R^3)$ (with $\div_h \tv_0^h=0$), $\tTheeo \in \dot{B}_{2,1}^{-\frac12}(\R)$ (for all $\ee>0$) such that:
 $$
 \|\tv_0^h\|_{H^{\frac12+\delta}(\R^3)}\leq \Co \quad \mbox{and} \quad \sup_{\ee>0} \|\tTheeo\|_{\dot{B}_{2,1}^{-\frac12}(\R)} \leq \Co.
 $$
 Thanks to Theorems \ref{ThHeat} and \ref{ThSNS}, $\tv^h$ and $\tThee$ globally exist (for all $\ee>0$) and respectively belong to $\dot{E}^0 \cap \dot{E}^{\frac12+\delta}$ and $\dot{B}^{-\frac12}$.

 Moreover there exists a constant $C_{\delta,\nu,\nu'}>0$ such that for any fixed $\ee>0$, if $\Uoe\in L^2(\R^3)$, then there exists a weak global solution of \eqref{StratifD} $(\De,\qe)$ with $\De \in \dot{E}^0$ and $\qe \in \dot{E}^1+L^\frac43 (\R_+,L^2)$, satisfying for all $t\geq 0$,
 \begin{multline}
  \|\De(t)\|_{L^2}^2 +\nu_0 \int_0^t \|\n \De (\tau)\|_{L^2}^2 d\tau\\
\leq \left(\|\Uoeosc\|_{L^2}^2 +\|\UoeS^h-\tvo^h\|_{L^2}^2 +C_{\delta,\nu,\nu'} \Co^{2+\frac1{\delta}}\right) e^{C_{\delta,\nu,\nu'} \Co^{2+\frac1{\delta}}}.
\label{estimaprioriLeray}
 \end{multline}
 }
 \label{ThLeray}
\end{thm}
\begin{rem}
 \sl{At this stage we do not ask that $\UoeS^h \rightarrow \tv_0^h$ and we could prove the result only substracting $(0,0,0,\tThee)$, but the present formulation will be helpful for the convergence (in the rest of the present article and in \cite{FCStratif2}).}
\end{rem}
In what follows, we will state a priori estimates, give the Friedrichs scheme and quickly outline what is new when adapting the proof of the classical Leray theorem from the Navier-Stokes case.

\subsection{A priori estimates}

\begin{prop}
 Let $\delta>0$, $\Co\geq 1$, $\tv_0^h\in H^{\frac12+\delta}(\R^3)$ (with $\div_h \tv_0^h=0$), $\tTheeo \in \dot{B}_{2,1}^{-\frac12}(\R)$ (for all $\ee>0$) with:
 $$
 \|\tv_0^h\|_{H^{\frac12+\delta}(\R^3)}\leq \Co \quad \mbox{and} \quad \sup_{\ee>0} \|\tTheeo\|_{\dot{B}_{2,1}^{-\frac12}(\R)} \leq \Co
 $$
 Then for any fixed $\ee>0$, if we set $\nu_0=\min(\nu,\nu')>0$, there exist $C_0>0$ and $C_{\delta,\nu,\nu'}>0$ such that for any solution $\De \in \dot{E}^0$ of System \eqref{StratifD}, and for all $t\geq 0$:
 \begin{multline}
  \|\De(t)\|_{L^2}^2 +\nu_0 \int_0^t \|\n \De (\tau)\|_{L^2}^2 d\tau \leq \left(\|\De(0)\|_{L^2}^2 +\frac12 \int_0^t \|\tG(t')\|_{L^2} dt'\right)\\
\times \exp\left(C_0 \int_0^t \Big[\|\tG(t')\|_{L^2} +\frac1{\nu_0} \|\n \tv^h(t')\|_{\dot{H}^\frac12}^2 +\frac1{\nu_0^{\frac13}} \|\tThee(t')\|_{\dot{H}^1}^\frac43\Big]dt' \right)\\
\leq \left(\|\De(0)\|_{L^2}^2 +C_{\delta,\nu,\nu'} \Co^{2+\frac1{\delta}}\right) e^{C_{\delta,\nu,\nu'} \Co^{2+\frac1{\delta}}}.
 \end{multline}
\label{estimapriori}
\end{prop}
\textbf{Proof:} taking the innerproduct in $L^2$ of \eqref{StratifD} with $\De=(\Ve,\He)$, and thanks to the fact that $\div \Ve=0=\div_h \tv^h$, we obtain that $(\De \cdot \n \De|\De)_{L^2}=0=(\tv^h \cdot \n_h \De|\De)_{L^2}$ and:
\begin{multline}
 \frac12 \frac{d}{dt}\|\De\|_{L^2}^2+\nu_0 \|\n\De\|_{L^2}^2 =-(\De \cdot \n \tv^h|\De^h)_{L^2} -(\De^3\cdot \d_3 \tThee|\De^4)_{L^2}+(\tG|\De)_{L^2}\\
 \overset{def}{=} A+B+C.
 \label{estimABC}
\end{multline}
The terms $A$ and $C$ are easy to bound and we immediately get that:
\begin{multline}
 |A| \leq \|\De \cdot \n \tv^h\|_{L^2} \|\De\|_{L^2} \leq \|\De\|_{L^6} \|\n \tv^h\|_{L^3} \|\De\|_{L^2} \leq \|\De\|_{\dot{H}^1} \|\n \tv^h\|_{\dot{H}^\frac12} \|\De\|_{L^2}\\
 \leq \frac{\nu_0}4 \|\De\|_{\dot{H}^1}^2 +\frac{C}{\nu_0} \|\n \tv^h\|_{\dot{H}^\frac12}^2 \|\De\|_{L^2}^2,
\end{multline}
and
\begin{equation}
|C|\leq \|\tG\|_{L^2} \|\De\|_{L^2} \leq \frac12 \|\tG\|_{L^2} +\frac12 \|\tG\|_{L^2} \|\De\|_{L^2}^2.
\end{equation}
The term $B$ is not classical as it features functions of one variable. Thanks to the Minkowski and Young estimates (using twice the coefficients $(\frac43,4)$ in the last lines), the Sobolev interpolation estimates and the 1d-Sobolev injection $\dot{H}^\frac14 (\R) \hookrightarrow L^4(\R)$:
\begin{multline}
|B|\leq \int_{\R^3} |\De^3(x_h,x_3)|\cdot |\d_3 \tThee (x_3)|\cdot |\De^4(x_h,x_3)| dx \leq \int_{\R^2} \left(\int_{\R} |\De(x_h,x_3)|^2 |\d_3 \tThee (x_3)| dx_3 \right) dx_h\\
 \leq \int_{\R^2} \|\De(x_h,\cdot)\|_{L^4(\R)}^2 \|\d_3 \tThee\|_{L^2(\R)} dx_h \leq C\|\d_3 \tThee\|_{L^2(\R)} \int_{\R^2} \|\De(x_h,\cdot)\|_{\dot{H}^\frac14(\R)}^2  dx_h\\
 \leq C\|\tThee\|_{\dot{H}^1(\R)} \int_{\R^2} \|\De(x_h,\cdot)\|_{L^2(\R)}^\frac32 \|\De(x_h,\cdot)\|_{\dot{H}^1(\R)}^\frac12 dx_h\\
 \leq C\|\tThee\|_{\dot{H}^1(\R)} \left(\int_{\R^2} \|\De(x_h,\cdot)\|_{L^2(\R)}^2 dx_h \right)^\frac34 \left(\int_{\R^2} \|\d_3 \De(x_h,\cdot)\|_{L^2(\R)}^2 dx_h\right)^\frac14\\
 \leq C\|\tThee\|_{\dot{H}^1(\R)} \|\De\|_{L^2(\R^3)}^\frac32 \|\d_3 \De\|_{L^2(\R^3)}^\frac12 \leq C\|\tThee\|_{\dot{H}^1(\R)} \|\De\|_{L^2}^\frac32 \|\n \De\|_{L^2}^\frac12\\
 \leq \frac{\nu_0}4 \|\De\|_{\dot{H}^1}^2 +\frac{C}{\nu_0^\frac13} \|\tThee\|_{\dot{H}^1(\R)}^\frac43 \|\De\|_{L^2}^2.
\end{multline}
Gathering these estimates into \eqref{estimABC} leads to:
$$
 \frac{d}{dt} \|\De\|_{L^2}^2+\nu_0 \|\n\De\|_{L^2}^2 \leq \|\tG\|_{L^2} +\|\De\|_{L^2}^2 \left(\|\tG\|_{L^2} +\frac{C}{\nu_0} \|\n \tv^h\|_{\dot{H}^\frac12}^2 +\frac{C}{\nu_0^\frac13} \|\tThee\|_{\dot{H}^1(\R)}^\frac43 \right),
$$
which entails the first part of the result after integration in time and thanks to the Gronwall Lemma. Thanks to Theorem \ref{ThSNS}, we bound the terms involving $\tG$ and $\tv_0^h$ and to bound the last term, we use the estimates from Theorem \ref{ThHeat} and more precisely \eqref{estimThetaBspr}.

We outline that it is not possible to directly use the estimate in $\dot{H}^s$ as the time integration exponent $p$ cannot reach $\frac43$. We cannot either use directly \eqref{estimThetaBspr} with $p=r=2$ and $q=\frac43$ as (thanks to Proposition \ref{Propermut}) we have:
$$
\|\tThee\|_{\tilde{L}_t^\frac43 \dot{H}^{s+\frac32}}\leq \|\tThee\|_{L_t^\frac43 \dot{H}^{s+\frac32}},
$$
and \eqref{estimThetaBspr} only bounds the left-hand side of the estimates (we need to bound the right-hand side). To simplify, we use the estimates in the case $(s,p,r,q)=(-\frac12,2,1,\frac43)$, and thanks to Proposition \ref{Propermut}, we have for all $t\geq 0$:
$$
\|\tThee\|_{L_t^\frac43 \dot{H}^1}\leq \|\tThee\|_{L_t^\frac43 \dot{B}_{2,1}^1}\leq \|\tThee\|_{\tilde{L}_t^\frac43 \dot{B}_{2,1}^1}\leq \frac{C}{(\nu')^\frac34} \|\tTheeo\|_{\dot{B}_{2,1}^{-\frac12}} \leq \frac{C \Co}{(\nu')^\frac34},
$$
which concludes the proof. $\blacksquare$

\subsection{The Friedrichs scheme}

In order to properly write the scheme corresponding to System \eqref{StratifD}, the first step is to express the pressure $\qe$ in terms of $\De=(\Ve,\He)$:
\begin{multline}
 \qe=-\frac1{\ee} \d_3 \D^{-1} \He -\D^{-1} \div \div \left(\De\otimes \De -(\tv^h,0) \otimes (\tv^h,0)\right)\\
 =-\frac1{\ee} \d_3 \D^{-1} \He -\D^{-1} \div \div \left(\Ve\otimes \Ve +\Ve \otimes (\tv^h,0) +(\tv^h,0) \otimes \Ve\right)\\
 = -\frac1{\ee} \d_3 \D^{-1} \He -\D^{-1} \left(\sum_{i,j=1}^3 \d_i \d_j (\Ve^i \Ve^j) +\sum_{i=1}^3 \sum_{j=1}^2 \d_i \d_j (\Ve^i \tv^j) +\sum_{i=1}^2 \sum_{j=1}^3 \d_i \d_j (\tv^i \Ve^j)\right),
 \label{defqe}
\end{multline}
which leads to the following scheme (where for $n\in\N$, $J_n$ is the Fourier truncation operator on the ball centered at zero and with radius $n$) for $\Den=(\Ven, \Hen)$:
\begin{multline}
 \d_t \Den-L J_n \Den +\frac1{\ee} \cB J_n\Den +\left(\begin{array}{c}J_n \div\Big[J_n \Ven \otimes J_n \Ven +J_n \Ven \otimes (\tv^h,0) +(\tv^h,0) \otimes J_n \Ven\Big] \\ J_n \div\Big[J_n \Ven \cdot (J_n \Hen +\tThee)\Big] +J_n \div_h\Big(\tv^h\cdot J_n \Hen \Big) \end{array}\right)\\
 = J_n \tG-\left(\begin{array}{c}\n J_n \qen \\0\end{array}\right),
\end{multline}
with initial data ${\Den}_{|t=0}=J_n \left(\Uoeosc +(\UoeS^h-\tvo^h,0,0)\right)$, the pressure $\qen$ being defined according to \eqref{defqe}.

\subsection{Sketch of the proof}

The rest of the proof of Theorem \ref{ThLeray} is classical and we only give a sketch of it, pointing out what is different.
\begin{itemize}
 \item The Friedrichs scheme is an ODE in $L_n^2=\{u\in L^2(\R^3), \mbox{ with }\mbox{supp }\hat{u}\subset B(0,n)\}$ and its solution $\Den$ belongs to $\mathcal{C}([0, T_{\ee,n}^*[, L_n^2(\R^3))$ for some lifespan $T_{\ee,n}^*$.
 \item We have $J_n \Den=\Den$ and, taking the divergence of the system, we show that $\div \Den=0$.
 \item Then the scheme can be rewritten in a form close to \eqref{StratifD} and Proposition \ref{estimapriori} implies that \eqref{estimaprioriLeray} is true for $\Den$, which provides a uniform bound in $n$ for $\|\De^n\|_{L^\infty L^2}$ which implies $T_{\ee,n}^*=+\infty$ thanks to the blow-up criterion for the ODE.
 \item As $\Den$ is bounded in $\dot{E}^0$, we can extract a subsequence that weakly converges to  some $\De\in \dot{E}^0$ satisfying \eqref{estimaprioriLeray} with the same bound.
 \item All that remains is to show the weak limit solves \eqref{StratifD}, the only difficulty being to deal with the nonlinear terms limits, for which the classical argument is to use the Ascoli theorem. This is here that some adaptation is needed. In the Navier-Stokes case, it is shown that every term in the right-hand side is in $L_{loc}^\frac43 \dot{H}^{-1}$. It is easily proved for the terms involving $\tv^h$ and all that remains is to check the term involving $\tThee$. With the same ideas as in  the proof of Proposition \ref{estimapriori}:
 \begin{multline}
  \|J_n \div (J_n \Ven \cdot \tThee)\|_{\dot{H}^{-1}}^2 \leq \|\Ven \cdot \tThee\|_{L^2}^2\\
  = \int_{\R^2} \left(\int_{\R} |\Ven(x_h,x_3|^2 |\tThee(x_3)|^2 dx_3 \right)dx_h \leq \int_{\R^2} \|\Ven(x_h,\cdot)\|_{L^4(\R)}^2 \|\tThee\|_{L^4(\R)}^2 dx_h\\
  \leq C \|\tThee\|_{\dot{H}^\frac14}^2 \int_{\R^2} \|\Ven (x_h,\cdot)\|_{L^2(\R)}^\frac32 \|\d_3 \Ven (x_h,\cdot)\|_{L^2(\R)}^\frac12 dx_h\\
  \leq C \|\tThee\|_{\dot{H}^\frac14}^2 \|\Ven\|_{L^2(\R^3)}^\frac32 \|\n \Ven\|_{L^2(\R^3)}^\frac12.
 \end{multline}
Thanks to the Young estimates (with $(\frac43,4)$), we obtain that:
$$
 \|J_n \div (J_n \Ven \cdot \tThee)\|_{L_t^2 \dot{H}^{-1}} \leq C \|\tThee\|_{L_t^\frac83 \dot{H}^\frac14} \|\Den\|_{L_t^\infty L^2(\R^3)}^\frac34 \|\n \Den\|_{L_t^2 L^2(\R^3)}^\frac14.
$$
The first term is bounded thanks to \eqref{estimThetaBspr} and the fact that $\tTheo\in \dot{B}_{2,1}^{-\frac12}$:
$$
\|\tThee\|_{L_t^\frac83 \dot{H}^\frac14}\leq \|\tThee\|_{L_t^\frac83 \dot{B}_{2,1}^\frac14}\leq \|\tThee\|_{\tilde{L}_t^\frac83 \dot{B}_{2,1}^\frac14}\leq \frac{C}{(\nu')^\frac38} \|\tTheeo\|_{\dot{B}_{2,1}^{-\frac12}} \leq \frac{C \Co}{(\nu')^\frac38},
$$
so that we have $J_n \div (J_n \Ven \cdot \tThee) \in L^2(\R_+, \dot{H}^{-1}) \subset L_{loc}^\frac43 (\R, \dot{H}^{-1})$.
\item Concerning the regularity of the pressure term $\qe$, only one term is different, namely $-\frac1{\ee} \d_3 \D^{-1} \De^4$, which obviously belongs to $\dot{E}^1$.
\end{itemize}

\section{Convergence of $\De$}

Let us first state more precisely the announced convergence result.

\begin{thm} (Convergence)
 \sl{For any $\delta>0$, $\Co\geq 1$, $\tTheo\in \dot{B}_{2,1}^{-\frac12}(\R)$ , $\tv_0^h\in H^{\frac12+\delta}(\R^3)^2$ (with $\div_h \tv_0^h=0$ or, equivalently, $ \tv_0^h=\bP_2\tv_0^h$) and any $\Uoe\in L^2$ (divergence-free) and $\tTheeo \in \dot{B}_{2,1}^{-\frac12}(\R)$ (for all $\ee>0$) with:
\begin{equation}
  \begin{cases}
  \vspace{0.1cm}
  \|\tv_0^h\|_{H^{\frac12+\delta}(\R^3)}\leq \Co,\\
  \vspace{0.1cm}
  \sup_{\ee>0} \|\Uoe\|_{L^2} \leq \Co,\\
  \|\UoeS^h-\tv_0^h\|_{L^2} \underset{\ee\rightarrow 0}{\longrightarrow} 0,
 \end{cases}
\quad \mbox{and} \quad
\begin{cases}
 \vspace{0.1cm}
 \|\tTheo\|_{\dot{B}_{2,1}^{-\frac12}}\leq \Co,\\
 \sup_{\ee>0} \|\tTheeo\|_{\dot{B}_{2,1}^{-\frac12}(\R)} \leq \Co,\\
 \|\tTheeo-\tTheo\|_{\dot{B}_{2,1}^{-\frac12}(\R)} \underset{\ee\rightarrow 0}{\longrightarrow} 0,
\end{cases}
\end{equation}
the global weak solution $\Ue$ (constructed in Theorem \ref{ThLeray}) converges to $(\tv^h,0,\tThe)$ (where $\tv^h$ and $\tThe$ are the global solutions of Systems \eqref{SNS3} and \eqref{SNS4}) in the following sense: if $\De=\Ue-(\tv^h,0, \tThee)$ (where $\tThee$ is the global solution of \eqref{SysttThee}), then
\begin{itemize}
 \item the stratified part $D_{\ee,S}=\bP_2 \De$ converges to zero: for all $q\in]2,6[$,
$$
\|D_{\ee,S}\|_{L_{loc}^2(\R_+, L_{loc}^q(\R^3)} \underset{\ee\rightarrow 0}{\longrightarrow} 0,
$$
\item the oscillating part $\Dosc =(I_d-\bP_2) \De$ converges to zero: for all $q\in]2,6[$, there exists $\ee_1=\ee_1(\nu,\nu',q)>0$ and, for all $t\geq 0$, a constant $\mathbb{D}_t=\mathbb{D}_{t,\delta,\nu,\nu',q,\Co}$ such that for all $\ee\in]0, \ee_1]$,

\begin{equation}
 \|\Dosc\|_{L_t^2 L^q} =\|\Dosc\|_{L^2([0,t], L^q(\R^3)} \leq \mathbb{D}_t \ee^{\frac{K(q)}{640}},\quad \mbox{with} \quad K(q)\overset{def}{=} \frac{\min(\frac6{q}-1,1-\frac2{q})^2}{ (\frac6{q}-1)}.
 \end{equation}
\end{itemize}
Moreover, when $\nu=\nu'$, the previous estimates can be upgraded into $\|\Dosc\|_{L_t^2 L^q} \leq \mathbb{D}_t \ee^{\frac{K(q)}{544}}$ (now valid for all $\ee>0$) and we can obtain global-in-time estimates with better convergence rate: there exists a constant $C=C_{\nu,\delta,\Co}>0$ such that, for any $\ee>0$,
$$
\|\Dosc\|_{\tilde{L}^\frac43 \dot{B}_{8, 2}^0 +\tilde{L}^1 \dot{B}_{8, 2}^0} \leq C \ee^\frac3{16}.
$$
 }
 \label{ThCV}
\end{thm}
\begin{rem}
 \sl{Choosing $\tTheeo=\tTheo$ for each $\ee>0$ leads to the particular case described by the Theorem from the introduction.}
\end{rem}

In the rest of this section, we begin to prove the convergence of the oscillating part $\Dosc$, and then we use it to prove that the stratified part $\DeS$ also goes to zero.

\subsection{Convergence of $\De$ in the case $\nu\neq\nu'$}

\subsubsection{Convergence of the oscillating part}

We refer to \eqref{defS} and \eqref{defosc} for the definition of $\DeS$ and $\Dosc$ and to Proposition \ref{PropSosc} for properties of this orthogonal decomposition. Thanks to Propositions \ref{PropSosc} (Points 6,7) and \ref{PropG}, we have $\cB \De=\cB \Dosc$ and $\bP (\n \qe,0)=0=\bP \tG$, which entails that, applying $\bP (I_d-\bP_2)$ to System \eqref{StratifD}, $\Dosc$ solves (recall that $(I_d-\bP_2)\Dosc=\Dosc$ and $\bP_2 \Dosc=0$):

\begin{equation}
\begin{cases}
 \d_t \Dosc -(L- \frac1{\ee} \bP\cB) \Dosc = -\bP (I_d-\bP_2) \left[\De\cdot \n \De
 +\left(\begin{array}{c}
  \De \cdot \n \tv^h \\0\\ \De^3\cdot \d_3 \tThee \end{array}\right) +\tv^h\cdot \n_h \De \right] +\tG,\\
 D_{\ee,osc|t=0}=\Uoeosc.
\end{cases}
\label{SystDosc}
\end{equation}
Then, referring to Section \ref{Troncatures} for the definition of $\cPrR$, as in \cite{CDGG, FC1} we truncate in frequency and split $\Dosc$ into three parts:
\begin{equation}
 \Dosc= \Big(I_d-\chi(\frac{|D|}{\Re})\Big) \Dosc +\chi(\frac{|D|}{\Re})\chi(\frac{|D_h|}{2\re}) \Dosc +\cPrR \Dosc.
\end{equation}
The first two terms are estimated as in \cite{CDGG, FC1} using Lemma \ref{lemaniso} and Propositions \ref{PropSosc} and \ref{estimaprioriLeray}, for all $t\geq 0$ and $q\in]2,6[$:
\begin{equation}
 \begin{cases}
 \vspace{0.2cm}
  \|(I_d-\chi(\frac{|D|}{\Re})) \Dosc\|_{L_t^2 L^q} \leq \Re^{-\frac12(\frac6{q}-1)} \|\Dosc\|_{L_t^2 \dot{H}^1} \leq \Re^{-\frac12(\frac6{q}-1)} \|\De\|_{L_t^2 \dot{H}^1},\\
  \|\chi(\frac{|D|}{\Re}) \chi(\frac{|D_h|}{2\re})\Dosc\|_{L_t^2 L^q} \leq
C t^\frac12 (\Re \re^2)^{\frac12-\frac1{q}} \|\Dosc\|_{L_t^\infty L^2} \leq
C t^\frac12 (\Re \re^2)^{\frac12-\frac1{q}} \|\De\|_{L_t^\infty L^2}. \end{cases}
\label{Tronc12}
\end{equation}
The third part is bounded thanks to the Strichartz estimates proved in Proposition \ref{Estimdispnu2}: when $\ee\leq \ee_1(\nu,\nu',m,M)$, choosing $(p,r)=(2,q)$ in \eqref{Estimdispnu3} and defining $\Fe$ as the right-hand side of System \eqref{SystDosc}, we obtain that:
\begin{equation}
  \|\cPrR \Dosc\|_{L_t^2 L^q} \leq C_{\nu, \nu',q} \frac{\Re^{7-\frac9{q}}}{\re^{\frac{15}2-\frac7{q}}} \ee^{\frac18 (1-\frac2{q})}\Big( \|\cPrR \Uoeosc\|_{L^2} +\|\cPrR \Fe\|_{L_t^1 L^2}\Big),
  \label{Tronc3}
 \end{equation}
All that remains is then to bound the external force terms. Three of them can be bounded using Lemma \ref{lemaniso}:
$$
\|\cPrR (fg)\|_{L_t^1 L^2}\leq C\Re^\frac32 \|fg\|_{L_t^1 L^1} \leq C \Re^\frac32 t^\frac12 \|f\|_{L_t^\infty L^2} \|g\|_{L_t^2 L^2},
$$
so that we immediately obtain:
\begin{multline}
\|\cPrR \left( \De\cdot \n \De +\De \cdot \n \tv^h +\tv^h\cdot \n_h \De \right)\|_{L_t^1 L^2}\\
\leq C \Re^\frac32 t^\frac12 \left((\|\De\|_{L_t^\infty L^2}  +\|\tv^h\|_{L_t^\infty L^2}) \|\n \De\|_{L_t^2 L^2} +\|\De\|_{L_t^\infty L^2} \|\n \tv^h\|_{L_t^2 L^2}\right).
\end{multline}
The last term is bounded as we did for the Friedrichs scheme:
\begin{multline}
\|\cPrR (\De^3\cdot \d_3 \tThee)\|_{L^2}^2 \leq \int_{\R^2} \|\De^3(x_h,\cdot)\|_{L^4(\R)}^2\|\d_3 \tThee\|_{L^4(\R)}^2 dx_h\\
\leq C\|\tThee\|_{\dot{H}^\frac54(\R)}^2 \int_{\R^2} \|\De^3(x_h,\cdot)\|_{\dot{H}^\frac14(\R)}^2 dx_h \leq C\|\tThee\|_{\dot{H}^\frac54(\R)}^2 \|\De\|_{L^2(\R^3)}^\frac32 \|\n \De\|_{L^2(\R^3)}^\frac12,
\end{multline}
and
\begin{equation}
 \|\cPrR (\De^3\cdot \d_3 \tThee)\|_{L_t^1 L^2} \leq \|\tThee\|_{L_t^\frac87 \dot{H}^\frac54} \|\De\|_{L_t^\infty L^2}^\frac34 \|\n \De\|_{L_t^2 L^2}^\frac14.
 \label{estimL1L2}
\end{equation}
Thanks to \eqref{estimThetaBspr}, as $\tTheo\in \dot{B}_{2,1}^{-\frac12}$:
\begin{equation}
\|\tThee\|_{L_t^\frac87 \dot{H}^\frac54}\leq \|\tThee\|_{L_t^\frac87 \dot{B}_{2,1}^\frac54}\leq \|\tThee\|_{\tilde{L}_t^\frac87 \dot{B}_{2,1}^\frac54}\leq \frac{C}{(\nu')^\frac78} \|\tTheeo\|_{\dot{B}_{2,1}^{-\frac12}} \leq \frac{C \Co}{(\nu')^\frac78},
\label{CondTheta}
\end{equation}
so that, thanks to Theorem \ref{ThSNS} and \eqref{estimaprioriLeray} we can write that:
\begin{equation}
 \|\cPrR \Fe\|_{L_t^1 L^2} \leq C_{\delta,\nu,\nu',q} (1+t^\frac12) \Re^\frac32 \left(\|\Uoeosc\|_{L^2}^2 +\|\UoeS^h-\tvo^h\|_{L^2}^2 + \Co^{2+\frac1{\delta}}\right) e^{C_{\delta,\nu,\nu'} \Co^{2+\frac1{\delta}}}.
 \label{Tronc4}
 \end{equation}
Gathering \eqref{Tronc12}, \eqref{Tronc3} and \eqref{Tronc4}, and thanks to the assumptions on the initial data, we finally obtain that for all $t\geq 0$ and $q\in]2,6[$ (replacing $\re,\Re$ by their expressions in $\ee$):
\begin{multline}
 \|\Dosc\|_{L_t^2 L^q} \leq C_{\delta,\nu,\nu',q} e^{C_{\delta,\nu,\nu'} \Co^{2+\frac1{\delta}}} \left(\max(1,\|\Uoeosc\|_{L^2})^2 +\max(1,\|\UoeS^h-\tvo^h\|_{L^2})^2 +\Co^{2+\frac1{\delta}}\right)\\
 \times \left(\Re^{-\frac12(\frac6{q}-1)} +t^\frac12 (\Re \re^2)^{\frac12-\frac1{q}} +\frac{\Re^{\frac{17}2-\frac9{q}}}{\re^{\frac{15}2-\frac7{q}}} \ee^{\frac18 (1-\frac2{q})}(1+t^\frac12) \right)\\
 \leq \mathbb{D}_t \Big(\max(1,\|\Uoeosc\|_{L^2})^2 +\max(1,\|\UoeS^h-\tvo^h\|_{L^2})^2\Big)\\
 \times \left(\ee^{\frac{M}2(\frac6{q}-1)}+\ee^{(2m-M)(\frac12-\frac1{q})} +\ee^{\frac18 (1-\frac2{q})-M(\frac{17}2-\frac9{q})-m(\frac{15}2-\frac7{q})}\right)\\
 \leq \mathbb{D}_t \ee^{J(m,M,q)},
 \label{estimDosc1}
\end{multline}
where $\mathbb{D}_t=\mathbb{D}_{t,\delta,\nu,\nu',q,\Co}$ and we define
\begin{multline}
 J(m,M,q)\\
 \overset{def}{=}\min\left(\frac{M}2(\frac6{q}-1), (2m-M)(\frac12-\frac1{q}), \frac18 (1-\frac2{q})-M(\frac{17}2-\frac9{q})-m(\frac{15}2-\frac7{q})\right).
\end{multline}
Let us recall that $m,M$ also have to satisfy:
\begin{equation}
 \begin{cases}
  3M+m<1,\\
  M<2m.
 \end{cases}
 \label{CondmM3}
\end{equation}
In order to provide a more explicit convergence rate, we can observe that the third term in the minimum becomes positive when $m,M$ are small enough, but we have to care that, due to the first two terms, they are not too small. Introducing $A,B$ as below (both of them are positive as $q\in]2,6[$),
$$
A= \frac6{q}-1 \quad \mbox{and} \quad B=1-\frac2{q},
$$
we rewrite the third term in the minimum as a function of $A$ and $B$, which leads to:
\begin{equation}
 J(m,M,q)= A\min\left(\frac{M}2, (m-\frac{M}2)\frac{B}{A}, (\frac18-\frac{21}2 M-\frac{19}2m)\frac{B}{A}-2(m+M)\right).
\end{equation}
Choosing $m=M$ (and taking into account \eqref{CondmM3}), we are reduced, given $A,B>0$, to choose $m$ small but with the view that the following quantity has to be as large as possible:
\begin{equation}
 J(m,m,q)= A\min\left(\frac{m}2, \frac{m}2\frac{B}{A}, (\frac18-20m)\frac{B}{A}-4m\right).
\end{equation}
Observing that:
$$
m\leq \frac1{320} \Longrightarrow
\begin{cases}
\vspace{0.2cm}
 \frac18-20 m \geq \frac1{16},\\
 J(m,m,q)\geq A\min\left(\frac{m}2, \frac{m}2\frac{B}{A}, \frac1{16}\frac{B}{A}-4m\right),
\end{cases}
$$
we obtain that if in addition $m=M\leq \frac1{128}\frac{B}{A}$, then $\frac1{16}\frac{B}{A}-4m \geq \frac1{32}\frac{B}{A}$ and
\begin{equation}
 J(m,m,q) \geq A \min\left(\frac{m}2, \frac{m}2\frac{B}{A}, \frac1{32}\frac{B}{A}\right) \geq A \min(1,\frac{B}{A}) \min\left(\frac{m}2, \frac1{32}\right)=\min (A,B) \frac{m}2.
\end{equation}
Finally, if we choose $m_q=M_q=\min(\frac1{320}, \frac1{128}\frac{B}{A})$, we obtain:
\begin{equation}
 J(m_q,M_q,q) \geq \frac1{640} \min(1, \frac{B}{A})\min (A,B) =\frac{\min(A,B)^2}{640 A},
\end{equation}
so that, plugging into \eqref{estimDosc1}, we end-up with the fact that for all $t>0$, all $q\in]2,6[$ and all $\ee\leq \ee_1(\nu,\nu',q)$:
$$
 \|\Dosc\|_{L_t^2 L^q} \leq \mathbb{D}_{t,\delta,\nu,\nu',q,\Co} \ee^{\frac{\min(\frac6{q}-1,1-\frac2{q})^2}{640 (\frac6{q}-1)}}. \blacksquare
$$
\begin{rem}
 \sl{We emphasize that $\ee_1=\ee_1(\nu,\nu',m,M) \overset{def}{=}\left(\frac{\sqrt{2}}{|\nu-\nu'|}\right)^\frac1{1-(3M+m)}$ turns into $\ee_1(\nu, \nu',q)$ as $m,M$ are fixed only depending on $q$.}
\end{rem}

\subsubsection{Convergence of $\DeS$}

Thanks to the assumptions on the initial data, the estimates \eqref{estimaprioriLeray} can be bounded uniformly in $\ee$ as follows: for all $t\geq 0$, and all $\ee>0$
\begin{equation}
  \|\De(t)\|_{L^2}^2 +\nu_0 \int_0^t \|\n \De (\tau)\|_{L^2}^2 d\tau\\
\leq C_{\delta,\nu,\nu'} \Co^{2+\frac1{\delta}} e^{C_{\delta,\nu,\nu'} \Co^{2+\frac1{\delta}}}.
\label{estimaprioriLeray2}
 \end{equation}
Therefore we can extract from $(\De)_{\ee>0}$ a subsequence that weakly converges to some $\tD\in \dot{E}^0$ also satisfying \eqref{estimaprioriLeray2}. To simplify, let us also denote $(\De)_{\ee>0}$ this subsequence. We showed in the previous section that for all $q\in]2,6[$ and all $t>0$, $\|\Dosc\|_{L_t^2 L^q} \underset{\ee \rightarrow 0}{\longrightarrow}0$, which in particular implies that $\Dosc \underset{\ee \rightarrow 0}{\rightharpoonup}0$ and $\DeS \underset{\ee \rightarrow 0}{\rightharpoonup} \tD$.

Applying the operator $\bP_2$ to System \eqref{StratifD}, we obtain that $\DeS=(\DeS^h,0,0)$ satisfies (using once more Proposition \ref{PropSosc}):
\begin{multline}
 \d_t \DeS -L \DeS = -\bP_2\Big[\DeS\cdot \n \DeS +\DeS \cdot \n \left(\begin{array}{c}\tv^h\\ 0\\0\end{array}\right) +\tv^h\cdot \n_h \DeS\\
 +\DeS\cdot \n \Dosc +\Dosc\cdot \n \DeS +\Dosc\cdot \n \Dosc +\Dosc \cdot \n \left(\begin{array}{c}\tv^h\\ 0\\0\end{array}\right) +\tv^h\cdot \n_h \Dosc\Big],
\end{multline}
with initial data: $D_{\ee,S|t=0}=(\UoeS^h-\tv_0^h,0,0)$.
\\

Now, as $\DeS \underset{\ee \rightarrow 0}{\rightharpoonup} \tD$ (weakly) and $\Dosc \underset{\ee \rightarrow 0}{\longrightarrow}0$ (strongly in $L_{loc}^2 L^q$), every term involving $\Dosc$ weakly converges to zero. All that remains is to show $\DeS\cdot \n \DeS$ weakly converges to $\tD\cdot \n \tD$ which is similar to what is done in the classical Navier-Stokes case. We obtain that $\tD$ satisfies the following system:
\begin{equation}
\begin{cases}
  \d_t \tD -L \tD = -\bP_2\Big[\tD\cdot \n \tD +\tD \cdot \n \left(\begin{array}{c}\tv^h\\ 0\\0\end{array}\right) +\tv^h\cdot \n_h \tD\Big],\\
  \tD_{|t=0}=0.
\end{cases}
\end{equation}
Thanks to Proposition \ref{estimapriori} (with zero initial data and external force $\tG$ replaced by zero), we obtain that $\tD=0$. More precisely, we proved that the only possible weak limit for any extracted sequence is zero, which implies the whole sequence itself weakly converges to zero. To finish the proof, we have to upgrade the convergence of $\DeS$ towards zero: we repeat the arguments from \cite{FC1} and manage to finally obtain that for all $q\in]2,6[$:
$$
\DeS \underset{\ee \rightarrow 0}{\longrightarrow}0 \quad \mbox{in }L_{loc}^2 L_{loc}^q.
$$

\subsection{Convergence of $\Dosc$ in the case $\nu=\nu'$}

In the case $\nu=\nu'$, we can use the previous arguments but we do not need anymore any smallness condition on $\ee$ and (taking advantage of the fact that the $\bP_k$ are orthogonal) we can save $\Re/\re$ in the Strichartz estimates \eqref{Estimdispnu3} and in the final estimates \eqref{estimDosc1}, which modifies $J(m,M,q)$ into:
\begin{multline}
 K(m,M,q)\\
 \overset{def}{=}\min\left(\frac{M}2(\frac6{q}-1), (2m-M)(\frac12-\frac1{q}), \frac18 (1-\frac2{q})-M(\frac{15}2-\frac9{q})-m(\frac{13}2-\frac7{q})\right),
\end{multline}
which allows a slight improvement in the final convergence rate (following the same steps as before):
 for all $t>0$, all $q\in]2,6[$ and all $\ee>0$:
$$
 \|\Dosc\|_{L_t^2 L^q} \leq \mathbb{D}_{t,\delta,\nu,\nu',q,\Co} \ee^{\frac{\min(\frac6{q}-1,1-\frac2{q})^2}{544 (\frac6{q}-1)}}.
$$
Nevertheless, when $\nu=\nu'$ we can obtain much better Strichartz estimates (see Proposition \ref{Estimdispnu} in the appendix) which do not require frequency truncations. So the convergence for $\Dosc$ is directly given by this proposition, the only restrictions in choosing the coefficients $(d,p,r,\theta)$ being guided by the possible bounds for the external force term. These new Strichartz estimates will allow global in time estimates and a better convergence rate.
\\

To do this we remark that, as in \cite{FCRF}, among the various external force terms, we can observe two distinct regularities, which suggests us to decompose $\Dosc=\Dosc^1+\Dosc^2$ which respectively solve:
\begin{equation}
\begin{cases}
 \d_t \Dosc^1 -(L- \frac1{\ee} \bP\cB) \Dosc^1 = -\bP (I_d-\bP_2)
 \left(\begin{array}{c} 0\\0\\0\\ \De^3\cdot \d_3 \tThee
 \end{array}\right) +\tG,\\
 D_{\ee,osc|t=0}^1=\Uoeosc,
\end{cases}
\label{SystDosc2}
\end{equation}
and
\begin{equation}
\begin{cases}
 \d_t \Dosc^2 -(L- \frac1{\ee} \bP\cB) \Dosc^2 = -\bP (I_d-\bP_2) \left[\De\cdot \n \De +\left(\begin{array}{c} \De \cdot \n \tv^h \\0\\0
 \end{array}\right) +\tv^h\cdot \n_h \De \right],\\
 D_{\ee,osc|t=0}^2=0.
\end{cases}
\label{SystDosc3}
\end{equation}
The initial data and external force terms in the first system can be respectively bounded in $L^2$ and $L^1(\R_+,L^2)$ (globally in time, we refer to \eqref{estimL1L2} and Theorem \ref{ThSNS}), whereas the external force terms in the second system can be bounded in $L^1(\R_+,\dot{H}^{-\frac12})$ according to the following product law: if $\div f=0$,
$$
\|f\cdot \n g\|_{\dot{H}^{-\frac12}} \leq \|fg\|_{\dot{H}^\frac12} \leq \|f\|_{\dot{H}^1} \|g\|_{\dot{H}^1}.
$$
Thanks to Proposition \ref{Estimdispnu}, for any $(p_1,p_2,r_1,r_2)$ satisfying:
\begin{equation}
\frac2{r_1}+\frac1{p_1}=1 \quad \mbox{and} \quad \frac2{r_2}+\frac1{p_2}=\frac54,
\label{condpr}
\end{equation}
we have (thanks to Theorem \ref{ThSNS}, \eqref{estimL1L2}, \eqref{estimaprioriLeray} and \eqref{CondTheta}) that for all $t\geq 0$:
\begin{multline}
 \|\Dosc^1\|_{\tilde{L}_t^{p_1}\dot{B}_{r_1, 2}^0} \leq C_{p_1,r_1,\nu} \ee^{\frac1{4}(1-\frac{2}{r_1})} \left(\|\Uoeosc\|_{L^2}+\int_0^t (\|\De^3\cdot \d_3 \tThee\|_{L^2}+\|\tG\|_{L^2})dt'\right)\\
 \leq C_{p_1,r_1,\nu} \ee^{\frac1{4}(1-\frac2{r_1})} \left(\|\Uoeosc\|_{L^2}+ \|\tThee\|_{L_t^\frac87 \dot{H}^\frac54} \|\De\|_{L_t^\infty L^2}^\frac34 \|\n \De\|_{L_t^2 L^2}^\frac14 +\|\tG\|_{L_t^1 L^2}\right)\\
 \leq C_{p_1,r_1,\nu} \ee^{\frac1{4}(1-\frac{2}{r_1})} \left(\|\Uoeosc\|_{L^2} +C_{\delta,\nu} \Co^{2+\frac1{\delta}} \right) e^{C_{\delta,\nu} \Co^{2+\frac1{\delta}}}\\
 \leq C_{p_1,r_1,\nu,\delta, \Co} \ee^{\frac1{4}(1-\frac{2}{r_1})} (\|\Uoeosc\|_{L^2} +1),
\end{multline}
and
\begin{multline}
 \|\Dosc^2\|_{\tilde{L}_t^{p_2}\dot{B}_{r_2, 2}^0} \leq C_{p_2,r_2,\nu} \ee^{\frac1{4}(1-\frac{2}{r_2})} \int_0^t \left(\|\De\|_{\dot{H}^1}^2+2 \|\De\|_{\dot{H}^1} \|\tv^h\|_{\dot{H}^1}\right)dt'\\
 \leq C_{p_2,r_2,\nu,\delta, \Co} \ee^{\frac1{4}(1-\frac{2}{r_2})} (\|\Uoeosc\|_{L^2}^2 +1).
\end{multline}
If we wish the largest possible exponent for $\ee$, we need $r_1,r_2$ to be as large as possible. Let us recall that $p_2, r_2$ also satisfy (see Proposition \ref{Estimdispnu}) $p_2\leq 4/(1-2/r_2)$ which is satisfied as \eqref{condpr} (and the fact that $p_2\geq 1$) leads to:
$$
\frac85 \leq r_2 \leq 8.
$$
Choosing $r_1=r_2=8$, we finally obtain $(p_1,p_2)=(\frac43,1)$ (the conditions from Proposition \ref{Estimdispnu} are satisfied) and:
$$
\|\Dosc\|_{\tilde{L}^\frac43 \dot{B}_{8, 2}^0 +\tilde{L}^1 \dot{B}_{8, 2}^0} \leq C_{\nu,\delta,\Co} \ee^\frac3{16} (\|\Uoeosc\|_{L^2} +1) \leq C_{\nu,\delta,\Co} \ee^\frac3{16},
$$
which concludes the proof. $\blacksquare$

\section{Appendix 1}

\subsection{Notations, Sobolev spaces and Littlewood-Paley decomposition}

For a complete presentation of the Sobolev spaces and the Littlewood-Paley decomposition, we refer to \cite{Dbook}. We will use the same notations as in the appendix of \cite{FCPAA}. Let us first mention the following lemma:
\begin{prop}
 \sl{(\cite{Dbook} Chapter 2) We have the following continuous injections:
$$
 \begin{cases}
\mbox{For any } p\geq 1, & \dot{B}_{p,1}^0 \hookrightarrow L^p,\\
\mbox{For any } p\in[2,\infty[, & \dot{B}_{p,2}^0 \hookrightarrow L^p,\\
\mbox{For any } p\in[1,2], & \dot{B}_{p,p}^0 \hookrightarrow L^p.
\end{cases}
$$
}
 \label{injectionLr}
\end{prop}
Sometimes it is more convenient to work in a slight modification of the classical $L_t^p \dot{B}_{q,r}^s$ Spaces: the Chemin-Lerner time-space Besov spaces. As explained in the following definition, the integration in time is performed before the summation with respect to the frequency decomposition index:
\begin{defi} \cite{Dbook}
 \sl{For $s,t\in \R$ and $a,b,c\in[1,\infty]$, we define the following norm
 $$
 \|u\|_{\tilde{L}_t^a \dot{B}_{b,c}^s}= \Big\| \left(2^{js}\|\ddj u\|_{L_t^a L^b}\right)_{j\in \Z}\Big\|_{l^c(\Z)}.
 $$
 The space $\tilde{L}_t^a \dot{B}_{b,c}^s$ is defined as the set of tempered distributions $u$ such that $\lim_{j \rightarrow -\infty} S_j u=0$ in $L^a([0,t],L^\infty(\R^d))$ and $\|u\|_{\tilde{L}_t^a \dot{B}_{b,c}^s} <\infty$.
 }
 \label{deftilde}
\end{defi}
We refer once more to \cite{Dbook} (Section 2.6.3) for more details and will only recall the following proposition:
\begin{prop}
\sl{
For all $a,b,c\in [1,\infty]$ and $s\in \R$:
     $$
     \begin{cases}
    \mbox{if } a\leq c,& \forall u\in L_t^a \dot{B}_{b,c}^s, \quad \|u\|_{\tilde{L}_t^a \dot{B}_{b,c}^s} \leq \|u\|_{L_t^a \dot{B}_{b,c}^s}\\
    \mbox{if } a\geq c,& \forall u\in\tilde{L}_t^a \dot{B}_{b,c}^s, \quad \|u\|_{\tilde{L}_t^a \dot{B}_{b,c}^s} \geq \|u\|_{L_t^a \dot{B}_{b,c}^s}.
     \end{cases}
     $$
     \label{Propermut}
     }
\end{prop}

\subsection{Truncations}
\label{Troncatures}
In this section we precise the truncation operator that we will abundantly use: let us choose a function $\chi\in \mathcal{C}_0^\infty (\R, \R)$ taking values into $[0,1]$ and such that:
$$
\begin{cases}
 \mbox{supp }\chi \subset [-1,1],\\
 \chi \equiv 1 \mbox{ near } [-\frac12,\frac12].
\end{cases}
$$
Given $0<r<R$ we will denote by $\cC_{r,R}$ the following set (where $\xi=(\xi_h,\xi_3)$):
\begin{equation}
\cC_{r,R} =\{\xi \in \R^3, \quad |\xi|\leq R \mbox{ and } |\xi_h|\geq r\}.
 \label{CrR}
\end{equation}
If we define $f_{r,R}(\xi)=\chi (\frac{|\xi|}{R})\big(1-\chi (\frac{|\xi_h|}{2r})\big)$, then:
\begin{equation}
 \begin{cases}
 \mbox{supp }f_{r,R} \subset \cC_{r,R},\\
 f_{r,R}\equiv 1 \mbox{ on } \cC_{2r,\frac{R}2}.
\end{cases}
\label{frR}
\end{equation}
Now we introduce the following frequency truncation operator on $\cC_{r,R}$ ($\mathcal{F}^{-1}$ is the inverse Fourier transform and $|D|^s$ the classical derivation operator: $|D|^s f =\mathcal{F}^{-1} (|\xi|^s \hat{f}(\xi))$.):
\begin{multline}
 \cP_{r,R} u= f_{r,R}(D) u =\chi (\frac{|D|}{R})\big(1-\chi (\frac{|D_h|}{2r})\big) u \\
 =\mathcal{F}^{-1} \Big(f_{r,R}(\xi) \hat{u}(\xi)\Big) = \mathcal{F}^{-1} \Big(\chi (\frac{|\xi|}{R})\big(1-\chi (\frac{|\xi_h|}{2r})\big) \hat{u}(\xi)\Big),
\label{PrR}
\end{multline}
Thanks to \eqref{frR}, we have:
\begin{equation}
 f_{\frac{r}2,2R}(D) f_{r,R}(D) u =f_{r,R}(D) u.
\end{equation}
In what follows (and as in \cite{FC5, FCPAA}) we will use it for particular radii $r_\ee=\ee^m$ and $R_\ee =\ee^{-M}$, where $m$ and $M$ will be precised later. Let us end with the following anisotropic Bernstein-type result (easily adapted from \cite{FC1}, see \cite{Dragos1} for more general anisotropic estimates):
\begin{lem}
\sl{There exists a constant $C>0$ such that for all function $f$, $\aa>0$, $1\leq q \leq p \leq \infty$ and all $0<r<R$, we have
\begin{equation}
\|\chi (\frac{|D|}{R}) \chi (\frac{|D_h|}{r}) f\|_{L^p} \leq C(R r^2)^{\frac{1}{q}-\frac{1}{p}} \|\chi (\frac{|D|}{R}) \chi (\frac{|D_h|}{r}) f\|_{L^q}.
\end{equation}
Moreover if $f$ has its frequencies located in $\cC_{r,R}$, then
$$
\||D|^\aa f\|_{L^p} \leq C R^\aa \|f\|_{L^p}. \blacksquare
$$
}
\label{lemaniso}
\end{lem}

\subsection{About the linearized system}

Consider the linearized system ($f_0, \Fe$ being divergence-free, the second form is obtained using the Leray projector $\bP$):
\begin{equation}
\begin{cases}
\d_t f-(L-\frac{1}{\ee} \cB) f=\Fe,\\
\div f=0,\\
f_{|t=0}=f_0.
\end{cases}
\Longleftrightarrow
\begin{cases}
\d_t f-(L-\frac{1}{\ee} \bP \cB) f=\Fe,\\
f_{|t=0}=f_0.
\end{cases}
\label{systdisp}
\end{equation}
Applying the Fourier transform turns the equation into (as in \cite{FC1, Scro3}):
$$
\d_t \hat{f}- \mathbb{B}(\xi, \ee)\hat{f}=\hat{\Fe},
$$
where
$$\mathbb{B}(\xi, \ee)= \hat{L-\frac{1}{\ee} \bP \cB} =\left(
\begin{array}{cccc}
-\nu(\xi_2^2+\xi_3^2) & \nu \xi_1 \xi_2 & \nu \xi_1 \xi_3 & \displaystyle{\frac{\xi_1\xi_3}{\ee |\xi|^2}}\\
\nu \xi_1 \xi_2 & -\nu(\xi_1^1+\xi_3^2) & \nu \xi_2 \xi_3 & \displaystyle{\frac{\xi_2\xi_3}{\ee |\xi|^2}}\\
\nu \xi_1 \xi_3 & \nu \xi_2 \xi_3 & -\nu(\xi_1^1+\xi_2^2) & \displaystyle{-\frac{\xi_1^2+\xi_2^2}{\ee |\xi|^2}}\\
0 & 0 & \displaystyle{\frac{1}{\ee}} & \displaystyle{-\nu'|\xi|^2}
\end{array}
\right).
$$
\begin{rem}
 \sl{
 Note that in \cite{Scro3, LT, MuScho} the authors consider the matrix $\hat{L-\frac{1}{\ee} \bP \cB} \bP$. Doing this will gather our first two eigenvalues (see below) into the same double eigenvalue}.
\end{rem}

\subsubsection{Eigenvalues, projectors}

We begin with the eigenvalues and eigenvectors of matrix $\mathbb{B}(\xi, \ee)$. The methods are similar to those in \cite{FC1, FC2, FC3, FC5, FCpochesLp}. As there are some differences we will give precise results and skip details. First, the characteristic polynomial is:
$$
\mbox{det} \big(X I_4-\mathbb{B}(\xi, \ee)\big)= X(X+\nu |\xi|^2)\Big(X^2+(\nu+\nu')|\xi|^2 X+ \nu \nu' |\xi|^4 +\frac{\xi_1^2+\xi_2^2}{\ee^2 |\xi|^2} \Big),
$$
whose roots are much simpler to obtain compared to System $(PE_\ee)$. The discriminent of the degree 2 factor is:
$$
D=(\nu-\nu')^2|\xi|^4-4\frac{\xi_1^2+\xi_2^2}{\ee^2 |\xi|^2}.
$$
Which is nonpositive if, and only if, $|\nu-\nu'| \ee |\xi|^3 \leq 2|\xi_h|$. If $\xi \in \cC_{r,R}$ it is sufficient to ask that $|\nu-\nu'| \ee R^3 \leq 2 r$. This asks to split the discussion into two cases.

\begin{itemize}
\item When $\nu\neq\nu'$, if we choose $(\re,\Re)=(\ee^{m}, \ee^{-M})$ ($m,M$ being precised later) then $D<0$ on $\cC_{\re,\Re}$ as soon as $|\nu-\nu'| \ee \Re^3 \leq 2 \re$, which is equivalent to:
\begin{equation}
 \ee^{1-(3M+m)} \leq \frac{2}{|\nu-\nu'|}.
\end{equation}
So, when
\begin{equation}
 3M+m<1 \quad \mbox{and} \quad 0<\ee \leq \ee_0=\ee_0(\nu,\nu',m,M) \overset{def}{=}\left(\frac{2}{|\nu-\nu'|}\right)^\frac1{1-(3M+m)},
\label{CondmM1}
 \end{equation}
we obtain the following eigenvalues:
$$
\lambda_1\exi=0, \quad \lambda_2\exi=-\nu |\xi|^2, \quad \lambda_{3,4}\exi=\lambda_{\pm}\exi,
$$
with, for $\eta\in\{-1,1\}$:
\begin{equation}
\lambda_{\eta} \exi=-\frac{\nu+\nu'}2 |\xi|^2+i\eta \frac{|\xi_h|}{\ee|\xi|} \sqrt{1-\frac{(\nu-\nu')^2 \ee^2 |\xi|^6}{4|\xi_h|^2}}.
\label{vp34}
 \end{equation}
The eigenvectors are as follows:
$$
V_1\exi=\left(\begin{array}{c} (\ee \nu \nu' |\xi|^2 +\frac1{\ee|\xi|^2})\xi_1\\ (\ee \nu \nu' |\xi|^2 +\frac1{\ee|\xi|^2})\xi_2\\ \ee \nu \nu' |\xi|^2 \xi_3,\\ \nu \xi_3 \end{array}\right), \quad V_2\exi=\frac1{|\xi_h|}\left(\begin{array}{c} -\xi_2 \\ \xi_1\\ 0\\ 0 \end{array}\right),
$$
and for $\eta\in\{-1,1\}$, $V_{3,4}\exi=V_{\pm}\exi$:
\begin{equation}
V_\eta\exi= \left(\begin{array}{c}\frac{\xi_1 \xi_3}{\sqrt{2} |\xi| |\xi_h|}\\ \frac{\xi_2 \xi_3}{\sqrt{2} |\xi| |\xi_h|}\\ -\frac{|\xi_h|}{\sqrt{2} |\xi|}\\  \frac{(\nu-\nu')\ee|\xi|^3}{2\sqrt{2} |\xi_h|} + i\frac{\eta}{\sqrt{2}} \sqrt{1-\frac{(\nu-\nu')^2 \ee^2 |\xi|^6}{4|\xi_h|^2}}\big)\end{array}\right).
\end{equation}
 \item When $\nu=\nu'$, $D<0$ on $\R^3$ and the eigenvalues become:
$$
\lambda_1\exi=0, \quad \lambda_2\exi=-\nu |\xi|^2, \quad \lambda_{3,4}\exi=\lambda_{\pm}\exi,
$$
with, for $\eta\in\{-1,1\}$:
 \begin{equation}
  \lambda_{\eta}\exi=-\nu |\xi|^2+i\eta \frac{|\xi_h|}{\ee|\xi|},
 \end{equation}
with the corresponding eigenvectors:
 \begin{equation}
 V_1\exi=\left(\begin{array}{c} (\ee \nu^2 |\xi|^2 +\frac1{\ee|\xi|^2})\xi_1\\ (\ee \nu^2 |\xi|^2 +\frac1{\ee|\xi|^2})\xi_2\\ \ee \nu^2 |\xi|^2 \xi_3,\\ \nu \xi_3 \end{array}\right),
 \quad V_2\exi=\frac1{|\xi_h|}\left(\begin{array}{c} -\xi_2 \\ \xi_1\\ 0\\ 0 \end{array}\right),
 \quad V_\eta\exi= \left(\begin{array}{c}\frac{\xi_1 \xi_3}{\sqrt{2} |\xi| |\xi_h|}\\ \frac{\xi_2 \xi_3}{\sqrt{2} |\xi| |\xi_h|}\\ -\frac{|\xi_h|}{\sqrt{2} |\xi|}\\  i\frac{\eta}{\sqrt{2}} \end{array}\right)
\end{equation}
\end{itemize}
In both cases, the first eigenvector is not orthogonal to $(\xi,0)$ so does not correspond to a divergence-free quantity and will then play no role in the study of the linearized system. We can easily see that in both cases $V_2$ is orthogonal to $V_3$ and $V_4$, but $V_3$ and $V_4$ are orthogonal only when $\nu=\nu'$.

As in \cite{FC1}, any $\R^4$-valued divergence-free function f can have its Fourier transform decomposed into the family $(V_2\exi, V_3\exi, V_4\exi)$ according to:
$$
\hat{f}(\xi)=\sum_{k=2}^4 a_k\exi V_k\exi.
$$
Denoting $\mathcal{P}_k\exi= a_k\exi V_k\exi$ we introduce the corresponding projectors: for $k\in \{2,3,4\}$ we define the Fourier projector $\bP_k=\bP_k(\ee,D)$ as:
\begin{equation}
\bP_k f =\cF^{-1} \Big(\mathcal{P}_k\exi (\hat{f}(\xi))\Big).
 \label{defPk}
\end{equation}
Due to the fact that $V_2$ is orthogonal to the other two eigenvectors, we have:
\begin{equation}
 \bP_2 f =\cF^{-1} \Big((\hat{f}(\xi)\cdot V_2\exi) V_2\exi\Big),
\label{defP2}
 \end{equation}
and it is only in the case $\nu=\nu'$ that similar properties also stand for $k=3,4$:
$$
\bP_k f =\cF^{-1} \Big((\hat{f}(\xi)\cdot V_k\exi) V_k\exi\Big).
$$
We gather in the following proposition the properties we will use to obtain the Strichartz estimates.

\begin{prop}
\label{estimvp}
\sl{
If $\nu\neq \nu'$, for all $m,M>0$ with $3M+m<1$, for all $\ee< \ee_1= \left(\frac{\sqrt{2}}{|\nu-\nu'|}\right)^\frac1{1-(3M+m)}$, if $\re=\ee^m$ and $\Re =\ee^{-M}$, then for all $\xi \in \mathcal{C}_{r_{\ee}, R_{\ee}}$, the matrix $\mathbb{B}(\xi, \ee) = \widehat{L-\frac{1}{\ee} \bP \cB}$ is diagonalizable and its eigenvalues satisfy:
\begin{equation}
\label{vp}
\begin{cases}
\vspace{0.2cm} \lambda_1\exi=0,\\
\vspace{0.2cm} \lambda_2\exi=-\nu |\xi|^2,\\
\vspace{0.2cm} \lambda_3\exi = -\frac{\nu+\nu'}2 |\xi|^2+i\frac{|\xi_h|}{\ee|\xi|} -i \ee D\exi,\\
\vspace{0.2cm} \lambda_4\exi =\overline{\lambda_3\exi},
\end{cases}
\end{equation}
where $D\exi$ satisfies for all $\xi \in\cC_{\re, \Re}$ (with $k\in\{1,2\}$):
$$
\begin{cases}
\vspace{0.2cm}
|D\exi| \leq (\nu-\nu')^2\frac1{4\sqrt{2}} \frac{|\xi|^5}{|\xi_h|} \leq C_0 (\nu-\nu')^2 \frac{\Re^5}{\re}=C_0 (\nu-\nu')^2 \ee^{-(5M+m)},\\
\vspace{0.2cm}
|\d_{\xi_k}D\exi| \leq (\nu-\nu')^2\frac{9}{2\sqrt{2}} \frac{|\xi|^5}{|\xi_h|^2} \leq C_0 (\nu-\nu')^2 \frac{\Re^5}{\re^2}=C_0 (\nu-\nu')^2 \ee^{-(5M+2m)},\\
|\d_{\xi_3}D\exi| \leq (\nu-\nu')^2\frac{15}{4\sqrt{2}} \frac{|\xi|^4}{|\xi_h|} \leq C_0 (\nu-\nu')^2 \frac{\Re^4}{\re}=C_0 (\nu-\nu')^2 \ee^{-(4M+m)},
\end{cases}
$$
Moreover, the projectors $\bP_k=\bP_k(\ee,D)$ satisfy that, for any divergence-free $\R^4$-valued vectorfield $f$, we have:
\begin{equation}
\begin{cases}
\vspace{0.2cm}
 \bP_2 f =(\n_h^\perp \D_h^{-1} \omega(f),0,0), \quad \mbox{with }\omega(f)=\d_1 f^2-\d_2 f^1,\\
 \|\bP_2 f\|_{\Hs} \leq \|(f^1,f^2)\|_{\Hs} \leq \|f\|_{\Hs}, \quad \mbox{for any }s\in \R.
\end{cases}
\label{estimvp2}
\end{equation}
and
\begin{equation}
\begin{cases}
\vspace{0.2cm}
 (I_d-\bP_2) f =(\n_h \D_h^{-1} \div_h f^h,f^3,f^3), \quad \mbox{with }\div_h f^h=\d_1 f^1+\d_2 f^2,\\
 \|(I_d-\bP_2) f\|_{\Hs} \leq \|f\|_{\Hs}, \quad \mbox{for any }s\in \R.
\end{cases}
\end{equation}
Finally for $k=3,4$,
\begin{equation}
\label{estimvp34}
\|\bP_k \cP_{\re, \Re} f\|_{\Hs }\leq \sqrt{2} \frac{\Re}{\re} \|\cP_{\re, \Re} f\|_{\Hs } =\sqrt{2} \ee^{-(m+M)} \|\cP_{\re, \Re} f\|_{\Hs }.
\end{equation}
If $\nu=\nu'$, there is no need anymore of a frequency truncation or an expansion for the last two eigenvalues (no $\ee_1$ either is necessary), and the $\bP_k$ ($k\in\{2,3,4\}$) are orthogonal so for any divergence-free $\R^4$-valued vectorfield $f$, we have:
$$
\|\bP_k f\|_{\Hs} \leq \|f\|_{\Hs}, \quad \mbox{for any }s\in \R.
$$}
\end{prop}
\textbf{Proof: }to obtain the developpment of $\lambda_j\exi$ ($j\in \{3,4\}$) we use \eqref{vp34} and the following Taylor expansion of order one:
$$
\sqrt{1+x}=1+\frac{x}2 \int_0^1 \frac{du}{\sqrt{1+xu}},
$$
so that we immediately get:
\begin{equation}
 D\exi= (\nu-\nu')^2 \frac{|\xi|^5}{8|\xi_h|} \int_0^1 \frac{du}{\sqrt{1-u\frac{(\nu-\nu')^2 \ee^2 |\xi|^6}{4|\xi_h|^2}}}.
\end{equation}
To bound $D$ and its derivatives, we impose a little stronger condition than \eqref{CondmM1} on $\ee$, namely we ask that:
\begin{equation}
 \Big( \ee\frac{|\nu-\nu'| |\xi|^3}{2|\xi_h|}\Big)^2 \leq \Big( \ee\frac{|\nu-\nu'| \Re^3}{2 \re}\Big)^2\leq \frac12, \quad \mbox{which leads to} \quad \ee \leq \left(\frac{\sqrt{2}}{|\nu-\nu'|}\right)^\frac1{1-(3M+m)}.
\label{CondmM2}
\end{equation}
To bound the projectors norms, the methods are similar to what we did in \cite{FC1} and we also use \eqref{CondmM2} in the cases $i=3,4$, as we obtain:
$$
|\cF(\bP_i \cP_{\re, \Re} f)(\xi)|\leq \frac{|\xi|}{\sqrt{2}\|\xi_h|}\frac1{\sqrt{1-\frac{(\nu-\nu')^2 \ee^2 |\xi|^6}{4|\xi_h|^2}}} |\cF\left(\cP_{\re, \Re} (f^3,f^4)\right)|,
$$
which gives \eqref{estimvp34}.

\subsubsection{Strichartz estimates when $\nu\neq\nu'$}

The aim of this section is to prove the following Strichartz estimates:
\begin{prop}
 \sl{For any $d\in \R$, $r\geq2$, $q\geq 1$ and $p\in[1,\frac{8}{1-\frac2{r}}]$, there exists a constant $C_{p,r}>0$ such that for any $\ee\in]0,\ee_1]$ (where $\ee_1= \left(\sqrt{2}/|\nu-\nu'|\right)^\frac1{1-(3M+m)}$) and any $f$ solving \eqref{systdisp} with initial data $f_0$ and external force $\Fe$ such that $\div f_0=\div \Fe=0$, then for $k=3,4$,
 \begin{multline}
  \||D|^d \bP_k \cPrR f\|_{\tilde{L}_t^p\dot{B}_{r, q}^0}\\
  \leq \frac{C_{p,r}}{(\nu+\nu')^{\frac1{p}-\frac18(1-\frac2{r})}}  \frac{\Re^{7 -\frac9{r}}}{\re^{\frac{13}2+\frac2{p}-\frac7{r}}} \ee^{\frac18(1-\frac{2}{r})}\left( \|\cPrR f_0\|_{\dot{B}_{2, q}^d}+ \|\cPrR \Fe\|_{L^1 \dot{B}_{2, q}^d}\right).
 \end{multline}
\label{Estimdispnu2}
}
\end{prop}
\begin{rem}
 \sl{The condition $\ee\leq \ee_1$ is needed only for us to use the estimates from Proposition \ref{estimvp}.}
\end{rem}
\textbf{Proof:} As in \cite{FC1,FCPAA,FCRF}, we first assume $\Fe=0$ (the inhomogeneous case is dealt reproducing the same steps on the Duhamel term), and we will prove that under the previous assumptions,
\begin{equation}
  \|\bP_k \cPrR f\|_{L_t^p L^r} \leq C_{\nu, \nu',p,r} \frac{\Re^{7 -\frac9{r}}}{\re^{\frac{13}2+\frac2{p}-\frac7{r}}} \ee^{\frac18(1-\frac{2}{r})}\Big( \|\cPrR f_0\|_{L^2} +\|\cPrR \Fe\|_{L^1 L^2}\Big),
  \label{Estimdispnu3}
 \end{equation}
which will give the result if applied to $\ddj |D|^d f$. The proof is close to the one of Proposition 51 from \cite{FCPAA}, but also features improvements coming from \cite{FCcompl} (as explained in this article, using the Riesz-Thorin theorem allows to upgrade the condition $r>4$ into $r\geq 2$). We will skip details and point out what is new. Let $\cA$ be the set:
$$
\cA\overset{def}{=}\{\psi \in \cC_0^\infty (\R_+\times \R^3, \R), \quad \|\psi\|_{L^{\bar{p}}(\R_+, L^{\bar{r}}(\R^3))}\leq 1\}.
$$
We follow the steps from \cite{FCcompl}: taking $k=3$, thanks to Plancherel and \eqref{PrR},
\begin{multline}
 \|\bP_k \cPrR f\|_{L^p L^r}= \sup_{\psi \in \cA} \int_0^\infty \int_{\R^3} \bP_k \cPrR f (t,x) \psi(t,x) dx dt\\
 =C \sup_{\psi \in \cA} \int_0^\infty \int_{\R^3} e^{-\frac{\nu+\nu'}2 t|\xi|^2+i\frac{t}{\ee}\frac{|\xi_h|}{|\xi|}-i t \ee D\exi} \cF\big(\bP_k \cPrR f_0\big)(\xi) f_{\frac{\re}2, 2\Re} (\xi)\hat{\psi}(t,\xi) d\xi dt\\
 \leq C \|\bP_k \cPrR f_0\|_{L^2} \sup_{\psi \in \cA} \Bigg[\int_0^\infty \int_0^\infty \int_{\R^3} f_{\frac{\re}2, 2\Re} (\xi)\hat{\psi}(t,\xi)\\
 \times  e^{-\frac{\nu+\nu'}2 (t+t')|\xi|^2+i\frac{t-t'}{\ee}\frac{|\xi_h|}{|\xi|}-i (t-t') \ee D\exi}  \overline{f_{\frac{\re}2, 2\Re} (\xi)\hat{\psi}(t',\xi)}dxdtdt'\Bigg]^{\frac12}.
 \end{multline}
Dispatching in both terms the heat semi-group, and using once more Plancherel and the H\"older inequality, we obtain:
\begin{multline}
 \|\bP_k \cPrR f\|_{L^p L^r} \leq C \|\bP_k \cPrR f_0\|_{L^2} \\
 \times \sup_{\psi \in \cA} \left[\int_0^\infty \int_0^\infty \|L_{\ee,t,t'} \psi(t,.)\|_{L^r} \|e^{\frac{\nu+\nu'}4 (t+t')\D} \overline{\cPrRb \psi (t',.)}\|_{L^{\bar{r}}} dtdt'\right]^{\frac12},
 \label{estimTT2}
\end{multline}
where we have defined for some $g$:
\begin{multline}
 \left(L_{\ee,t,t'} g\right)(x)= \int_{\R^3} e^{ix\cdot \xi} e^{-\frac{\nu+\nu'}4 (t+t')|\xi|^2+i\frac{t-t'}{\ee}\frac{|\xi_h|}{|\xi|}-i(t-t') \ee D\exi} f_{\frac{\re}2, 2\Re} (\xi)\hat{g}(\xi) d\xi\\
 = K_{\ee,t,t'}\star g,
\end{multline}
with
$$
K_{\ee,t,t'} (x)= (2\pi)^{-3}\int_{\R^3} e^{ix\cdot \xi} e^{-\frac{\nu+\nu'}4 (t+t')|\xi|^2+i\frac{t-t'}{\ee}\frac{|\xi_h|}{|\xi|}-i(t-t') \ee D\exi} f_{\frac{\re}2, 2\Re} (\xi) d\xi.
$$
To estimate the second norm in \eqref{estimTT2}, we cannot directly use the well-known smoothing effect of the heat flow described in Lemma 2.3 in \cite{Dbook} (see Section 2.1.2) as the frequencies are not supported in an annulus but in the set $\cC_{\frac{\re}2,2\Re}$, which needs us to adapt the result as stated in the following proposition (see \eqref{CrR} for the definition of the set $\cC_{r,R}$):
\begin{prop}
 \sl{Let $0<r<R$ be fixed. There exists a constant $C$ such that for any $p\in [1,\infty]$, $t\geq 0$ and any function $u$ we have:
 $$
 \mbox{Supp }\hat{u} \subset \cC_{r,R} \Rightarrow \|e^{t\D} u\|_{L^p} \leq C\frac{R^3}{r^4} e^{-\frac{t}2 r^2} \|u\|_{L^p}.
 $$
 }
 \label{Heatflow}
\end{prop}
We refer to the end of Section \ref{FinApp1} for the proof of this proposition, which allows us to write that:
\begin{equation}
 \|e^{\frac{\nu+\nu'}4 (t+t')\D} \overline{\cPrRb \psi (t',.)}\|_{L^{\bar{r}}} \leq C\frac{\Re^3}{\re^4} e^{-\frac{\nu+\nu'}{32} (t+t')\re^2} \|\psi(t',.)\|_{L^{\bar{r}}}.
 \label{estimheat}
\end{equation}
We start to deal with the other term through the same steps: with a view to use the Riesz-Thorin theorem, we bound the following norms:
\begin{equation}
 \begin{cases}
 \|L_{\ee,t,t'}\|_{L^2 \rightarrow L^2}\leq C_0 e^{-\frac{\nu+\nu'}{16} (t+t')\re^2},\\
 \|L_{\ee,t,t'}\|_{L^1 \rightarrow L^\infty}\leq  \|K_{\ee,t,t'}\|_{L^\infty}.
\end{cases}
\label{RieszTho}
\end{equation}
Thanks to the definition of the kernel, we easily obtain that
\begin{equation}
 \|K_{\ee,t,t'}\|_{L^\infty} \leq C_0 \Re^3 e^{-\frac{\nu+\nu'}{16} (t+t')\re^2}.
 \label{normeLinfiniK}
\end{equation}
In order to bound this norm with a negative power of $\frac{|t-t'|}{\ee}$, we cannot use the classical argument developped by Chemin, Desjardins, Gallagher and Grenier in \cite{CDGG, CDGG2, CDGGbook}. Reproducing the method would lead to:
$$
 |K_{\ee,t,t'}(x)|\leq \frac{C_0}{\re^2} e^{-\frac{\nu+\nu'}{32} (t+t')\re^2} \int_{\cC_{\frac{\re}2,2\Re}} \frac{d\xi}{1+\frac{|t-t'|}{\ee} \frac{\xi_2^2 \xi_3^4}{|\xi_h|^2|\xi|^6}}.
$$
Unfortunately, due to the definition of $\cC_{r,R}$ we cannot follow anymore the steps from the classical method, as we cannot bound from below on $\cC_{\frac{\re}2,2\Re}$ the term $\frac{\xi_2^2 \xi_3^4}{|\xi_h|^2|\xi|^6}$ (neither $\xi_2$ nor $\xi_3$ can be). The first idea would be to make $|\xi_2|$ easily bounded from below and change the set $\cC_{r,R}$ into:
$$
\cC'_{r,R}=\{\xi\in \R^3, |\xi|\leq R, \mbox{ and } |\xi_2|\geq r\},
$$
but in this case we would not be able to perform the key argument (invariance by rotation around the $x_3$ axis) to assume $x_2=0$ anymore. As in \cite{Scro3}, the best possibility is to use the third variable: thanks to the change of variable $(\xi_1,\xi_2,\xi_3) \mapsto (\xi_1,\xi_2,-\xi_3)$, we easily obtain that $K_{\ee,t,t'}(x_1,x_2,-x_3)=K_{\ee,t,t'}(x_1,x_2,x_3)$ (we refer to the previous proof for the expression of $D$), so that we can assume that $x_3\geq 0$. Moreover for any $t,t',\ee$, we have:
$$
\|K_{\ee,t,t'}\|_{L^\infty(\R^3)}= \sup_{x\in \R^3} \|K_{\ee,t,t'}(x_1,x_2,\frac{t-t'}{\ee} x_3)\|,
$$
so that we are reduced to bound ($C=(2\pi)^{-3}$):
$$
K_{\ee,t,t'}(x_1,x_2,\frac{t-t'}{\ee} x_3) =C\int_{\R^3} e^{ix_h\cdot \xi_h} e^{-\frac{\nu+\nu'}4 (t+t')|\xi|^2+i\frac{t-t'}{\ee} a(\xi)-i(t-t') \ee D\exi} \chi(\frac{|\xi|}{2\Re}) (1-\chi(\frac{|\xi_h|}{\re})) d\xi,
$$
where we have denoted:
\begin{equation}
a(\xi) \overset{def}{=} x_3\cdot \xi_3+\frac{|\xi_h|}{|\xi|}.
\end{equation}
Then we can introduce the operator $\mathcal{L}$, defined as follows:
\begin{equation}
 \mathcal{L}f=
 \begin{cases}
 \vspace{0.2cm}
  \displaystyle{\frac1{1+\frac{t-t'}{\ee}\aa(\xi)^2} \left(f(\xi)+i\aa(\xi) \d_{\xi_3} f(\xi)\right)} & \mbox{if } t>t',\\
  \displaystyle{\frac1{1+\frac{t'-t}{\ee}\aa(\xi)^2} \left(f(\xi)-i\aa(\xi) \d_{\xi_3} f(\xi)\right)} & \mbox{else },
 \end{cases}
\end{equation}
where we defined $\aa$ as follows:
$$
\aa(\xi)=-\d_{\xi_3} a(\xi)= -(x_3-\frac{\xi_3 |\xi_h|}{|\xi|^3}).
$$
Assume $t>t'$ for instance ($t-t'=|t-t'|$ in what follows), and performing an integration by parts, we have (as $\mathcal{L}$ is taylored to leave invariant $e^{i x_h\cdot \xi_h + i\frac{t-t'}{\ee}a(\xi)} (1-\chi(\frac{|\xi_h|}{\re}))$):
\begin{multline}
 K_{\ee,t,t'}(x_1,x_2,\frac{t-t'}{\ee} x_3)\\
 =C\int_{\R^3} e^{ix_h\cdot \xi_h+i\frac{t-t'}{\ee}a(\xi)}\big(1-\chi(\frac{|\xi_h|}{\re})\big) ^t\mathcal{L} \left(e^{-\frac{\nu+\nu'}4 (t+t')|\xi|^2-i(t-t') \ee D\exi} \chi(\frac{|\xi|}{2\Re})\right) d\xi,
\end{multline}
where the transposed operator is:
$$
 ^t\mathcal{L} f= \left(\frac1{1+\frac{t-t'}{\ee}\aa(\xi)^2}-i\d_{\xi_3} \aa(\xi) \frac{1-\frac{t-t'}{\ee}\aa(\xi)^2}{(1+\frac{t-t'}{\ee}\aa(\xi)^2)^2} \right)f(\xi) -\frac{i\aa(\xi)}{1+\frac{t-t'}{\ee}\aa(\xi)^2} \d_{\xi_3}f(\xi).
$$
Using the estimates from Proposition \ref{estimvp} for $D$, $\d_{\xi_3} D$, getting rid of $t,t'$ with the fact that $xe^{-x}\leq \frac2{e}e^{-\frac{x}2}$ for $x\geq 0$, writing $|\nu-\nu'|\leq \nu+\nu'$ and using \eqref{CondmM2} (the details are omitted but close to what is done for instance in \cite{FC1,FC3}), there exists some constant $C_0>0$ (depending on $\|\chi'\|_{L^\infty}$) such that:
\begin{multline}
 | ^t\mathcal{L} \left(e^{-\frac{\nu+\nu'}4 (t+t')|\xi|^2-i(t-t') \ee D\exi} \chi(\frac{|\xi|}{2\Re})\right)| \leq C_0 \frac{e^{-\frac{\nu+\nu'}4 (t+t')|\xi|^2}}{1+\frac{t-t'}{\ee} \aa^2}\\
 \times \left[(1+|\d_{\xi_3} \aa|) +|\aa|\big(\frac{\nu+\nu'}2(t+t')|\xi_3| +|t-t'| \ee |\d_{\xi_3} D(\ee,\xi)|+\frac{|\xi_3|}{|\xi|}\frac1{\Re}\big) \right]\\
 \leq C_0 \frac{e^{-\frac{\nu+\nu'}8 (t+t')|\xi|^2}}{1+\frac{t-t'}{\ee} \aa^2} \times \left[\big(1+4\frac{|\xi_h|}{|\xi|^3}\big) +|\aa|\Big(\frac4{e}\big(\frac{|\xi_3|}{|\xi|^2}+ 2C_0 \ee\frac{|\nu-\nu'|^2}{\nu+\nu'}\frac{|\xi|^2}{|\xi_h|}\big)+\frac{|\xi_3|}{|\xi|}\frac1{\Re}\Big) \right]\\
 \leq C_0 \frac{e^{-\frac{\nu+\nu'}{32} (t+t')\re^2}}{1+\frac{t-t'}{\ee} \aa^2} \left(\frac1{\re^2}+ \frac{|\aa|}{\re}\right),
\end{multline}
which leads to
\begin{multline}
 |K_{\ee,t,t'}(x_1,x_2,\frac{t-t'}{\ee} x_3)|\\
 \leq C_0\int_{\frac{\re}2\leq |\xi_h| \leq 2\Re} e^{-\frac{\nu+\nu'}{32} (t+t')\re^2} \int_{-\sqrt{(2\Re)^2-|\xi_h|^2}}^{\sqrt{(2\Re)^2-|\xi_h|^2}} \frac1{1+\frac{t-t'}{\ee} \aa^2} \left(\frac1{\re^2}+ \frac{|\aa|}{\re}\right) d\xi_3 d\xi_h.
 \label{normeLinfiniK2}
\end{multline}
Thanks to the fact that $|\aa| =\left(\frac{t-t'}{\ee}\right)^{-\frac12} \left(\frac{t-t'}{\ee}\right)^{\frac12} |\aa|\leq  \frac12\left(\frac{t-t'}{\ee}\right)^{-\frac12}(1+\frac{t-t'}{\ee} \aa^2)$, we have:
\begin{multline}
 \int_{\frac{\re}2\leq |\xi_h| \leq 2\Re} e^{-\frac{\nu+\nu'}{32} (t+t')\re^2} \int_{-\sqrt{(2\Re)^2-|\xi_h|^2}}^{\sqrt{(2\Re)^2-|\xi_h|^2}} \frac1{1+\frac{t-t'}{\ee} \aa^2} \frac{|\aa|}{\re} d\xi_3 d\xi_h\\
 \leq \frac1{\re} \left(\frac{t-t'}{\ee}\right)^{-\frac12} \int_{\frac{\re}2\leq |\xi_h| \leq 2\Re} e^{-\frac{\nu+\nu'}{32} (t+t')\re^2} \int_{-\sqrt{(2\Re)^2-|\xi_h|^2}}^{\sqrt{(2\Re)^2-|\xi_h|^2}} 1 d\xi_3 d\xi_h\\
 \leq C_0 \frac{\Re^3}{\re} \left(\frac{t-t'}{\ee}\right)^{-\frac12} e^{-\frac{\nu+\nu'}{32} (t+t')\re^2}.
\end{multline}
To bound the other part of the right-hand side of \eqref{normeLinfiniK2} we decompose the integral in $\xi_3$ as follows:
$$
\int_{-\sqrt{(2\Re)^2-|\xi_h|^2}}^0 \frac1{1+\frac{t-t'}{\ee} \aa^2} d\xi_3 +\int_0^{\sqrt{(2\Re)^2-|\xi_h|^2}} \frac1{1+\frac{t-t'}{\ee} \aa^2} d\xi_3.
$$
The first part is easily bounded: as $-\xi_3$ and $x_3$ are nonnegative,
$$
|\aa(\xi)|= |x_3-\frac{\xi_3|\xi_h|}{|\xi|^3}|= x_3-\frac{\xi_3|\xi_h|}{|\xi|^3} \geq |\xi_3|\frac{|\xi_h|}{|\xi|^3} \geq |\xi_3| \frac{\re}{16 \Re^3},
$$
and (thanks to the change of variable $z=\left(\frac{t-t'}{\ee}\right)^\frac12 \frac{\re}{16 \Re^3} \xi_3$):
\begin{equation}
 \int_{-\sqrt{(2\Re)^2-|\xi_h|^2}}^0 \frac1{1+\frac{t-t'}{\ee} \aa^2} d\xi_3 \leq \int_{-\sqrt{(2\Re)^2-|\xi_h|^2}}^0 \frac1{1+\frac{t-t'}{\ee} \frac{\xi_3^2 \re^2}{16^2 \Re^6}} d\xi_3
 \leq C_0 \left(\frac{t-t'}{\ee}\right)^{-\frac12} \frac{\Re^3}{\re}.
\end{equation}
The second part is bounded thanks to Proposition \ref{PropIab}: there exists a constant $C_0>0$ such that
\begin{multline}
 \int_0^{\sqrt{(2\Re)^2-|\xi_h|^2}} \frac1{1+\frac{t-t'}{\ee} \aa^2} d\xi_3 =I_{|\xi_h|,x_3}^{2\Re}(\frac{t-t'}{\ee}) \leq C_0 \frac{(2\Re)^7}{|\xi_h|^\frac{11}2} \min\left(1, \big(\frac{t-t'}{\ee}\big)^{-\frac14}\right)\\
 \leq C_0 \frac{\Re^7}{\re^\frac{11}2} \min\left(1, \big(\frac{t-t'}{\ee}\big)^{-\frac14}\right).
\end{multline}
Coming back to \eqref{normeLinfiniK2} when $t-t'=|t-t'|\geq \ee$, we end-up with:
\begin{multline}
  |K_{\ee,t,t'}(x_1,x_2,\frac{t-t'}{\ee} x_3)| \leq C_0\left(\frac{\Re^9}{\re^\frac{15}2} \min\Big(1, \big(\frac{t-t'}{\ee}\big)^{-\frac14}\Big) +\frac{\Re^5}{\re^3} \big(\frac{t-t'}{\ee}\big)^{-\frac12}\right) e^{-\frac{\nu+\nu'}{32} (t+t')\re^2}\\
  \leq C_0 \frac{\Re^9}{\re^\frac{15}2} \big(\frac{\ee}{t-t'}\big)^\frac14 e^{-\frac{\nu+\nu'}{32} (t+t')\re^2},
\end{multline}
which, coupled with \eqref{normeLinfiniK} when $|t-t'|\leq \ee$, finally implies that:
\begin{multline}
\|K_{\ee,t,t'}\|_{L^\infty} \leq C_0 \frac{\Re^9}{\re^\frac{15}2} \min\left(1,\big(\frac{\ee}{|t-t'|}\big)^\frac14\right) e^{-\frac{\nu+\nu'}{32} (t+t')\re^2}\\
\leq C_0 \frac{\Re^9}{\re^\frac{15}2} \frac{\ee^\frac14}{|t-t'|^\frac14} e^{-\frac{\nu+\nu'}{32} (t+t')\re^2}.
\label{taux}
\end{multline}
Using this together with \eqref{RieszTho}, we obtain thanks to the Riesz-Thorin theorem that for any $r\in[2, \infty]$:
$$
 \|L_{\ee,t,t'} g\|_{L^r} \leq C_0 \left(\frac{\Re^9}{\re^\frac{15}2} \frac{\ee^\frac14}{|t-t'|^\frac14}\right) ^{1-\frac2{r}} e^{-\frac{\nu+\nu'}{32} (t+t')\re^2} \|g\|_{L^{\bar{r}}}.
$$
Gathering this estimates together with \eqref{estimheat}, and thanks to \eqref{estimvp34}, we can properly bound \eqref{estimTT2} and obtain that:
\begin{multline}
  \|\bP_k \cPrR f\|_{L^p L^r}\\
  \leq C_0 \|\cPrR f_0\|_{L^2} \sup_{\psi \in \cA} \frac{\Re^{1+\frac32 +\frac92 (1-\frac2{r})}}{\re^{1+2+\frac{15}4 (1-\frac2{r})}} \ee^{\frac18 (1-\frac2{r})} \left[\int_0^\infty \int_0^\infty \frac{h(t)h(t')}{|t-t'|^{\frac14(1-\frac2{r})}} dtdt'\right]^{\frac12},
\end{multline}
with $h(t)=e^{-\frac{\nu+\nu'}{16} t\re^2} \|\psi(t,.)\|_{L^{\bar{r}}}$. As in \cite{FCPAA, FCcompl}, the last term is bounded using the Hardy-Littlewood-Sobolev estimates like in \cite{FCPAA, FCcompl}: if $k_1\geq 1$ is such that:
$$
\frac2{k_1}+\frac14(1-\frac2{r})=2,
$$
then
$$
\int_0^\infty \int_0^\infty \frac{h(t)h(t')}{|t-t'|^{\frac14(1-\frac2{r})}} dtdt' \leq C \|h\|_{L^{k_1}}^2,
$$
so that if $\bb\geq 1$ is defined as $\frac1{\bb}+\frac1{\bar{p}}=\frac1{k_1}$, that is $\frac1{\bb}=\frac1{p}-\frac18(1-\frac2{r})$ (this is here that we need the condition $p\leq \frac8{1-\frac2{r}}$), we have:
$$
\|h\|_{L^{k_1}} \leq \|e^{-\frac{\nu+\nu'}{16}\re^2 \cdot}\|_{L^\bb} \|\psi\|_{L_t^{\bar{p}} L_x^{\bar{r}}} \leq C \left(\frac{16}{(\nu+\nu')\re^2} \frac1{\bb}\right)^\frac1{\bb}\|\psi\|_{L_t^{\bar{p}} L_x^{\bar{r}}}.
$$
This finally entails that:
$$
 \|\bP_k \cPrR f\|_{L^p L^r} \leq \frac{C_{p,r}}{(\nu+\nu')^{\frac1{p}-\frac18(1-\frac2{r})}} \frac{\Re^{7-\frac9{r}}}{\re^{\frac{13}2+\frac2{p}-\frac7{r}}} \ee^{\frac18 (1-\frac2{r})} \|\cPrR f_0\|_{L^2},
$$
where $C_{p,r}=C_0 \big[16(\frac1{p}-\frac18(1-\frac2{r}))\big]^{\frac1{p}-\frac18(1-\frac2{r})}$, which concludes the proof. $\blacksquare$

\subsubsection{Strichartz estimates when $\nu=\nu'$}
\label{Disp}

When $\nu=\nu'$, we have $L=\nu\Delta$ and System \eqref{systdisp} becomes:
\begin{equation}
\begin{cases}
\d_t f-(\nu \D-\frac{1}{\ee} \bP \cB) f=\Fe,\\
f_{|t=0}=f_0.
\end{cases}
\label{systdispb}
\end{equation}
We prove in this section the following Strichartz estimates, that are close to their counterparts from \cite{FCcompl}:
\begin{prop}
 \sl{For any $d\in \R$, $r\geq 2$, $q\geq 1$, $\theta\in[0,1]$ and $p\in[1, \frac{4}{\theta (1-\frac{2}{r})}]$, there exists a constant $C=C_{p,r,\theta}$ such that for any $f$ solving \eqref{systdispb} for initial data $f_0$ and external force $\Fe$ both with zero divergence and vorticity (that is in the kernel of $\bP_2$), then
 \begin{equation}
  \||D|^d f\|_{\tilde{L}_t^p\dot{B}_{r, q}^0} \leq \frac{C_{p,r,\theta}}{\nu^{\frac{1}{p}-\frac{\theta}{4}(1-\frac{2}{r})}} \ee^{\frac{\theta}{4}(1-\frac{2}{r})} \left( \|f_0\|_{\dot{B}_{2, q}^{\sigma_1}} +\|\Fe\|_{\tilde{L}_t^1 \dot{B}_{2, q}^{\sigma_1}} \right),
 \end{equation}
 where $\sigma_1= d+\frac32-\frac{3}{r}-\frac{2}{p}+\frac{\theta}{2} (1-\frac{2}{r})$.
\label{Estimdispnu}
 }
\end{prop}
\textbf{Proof: }As the result is close to Proposition 4 from \cite{FCcompl} we will refer to this article for details and only point out what is different: namely the eigenvalues of the hessian and the singularity of the phase near $\{(0,0)\}\times \R$.

As in \cite{FCcompl}, the classical non stationnary phase argument from \cite{CDGG, FCPAA} is enhanced thanks to the Riesz-Thorin theorem (as in \cite{Dutrifoy2}) together with the Littman theorem (as in \cite{KLT, LT, IMT}) that we recall here:
\begin{thm} (Littman \cite{Litt, Stein2})
 \sl{Assume that $\psi:\R^n \rightarrow \R$ is a smooth function compactly supported in $K$ and $\phi:\R^n \rightarrow \R$ is a smooth function such that for any $\xi \in K$, the Hessian $D^2 \phi (\xi)$ has at least $k$ nonzero eigenvalues. Then there exists a constant $A$ such that for any $\lambda\in \R$ and $x\in\R^n$,
 $$
 |\int_{\R^n} e^{ix\cdot \xi +i \lambda \phi(\xi)} \psi(\xi) d\xi| \leq A\sqrt{|x|^2+\lambda^2}^{-\frac{k}2} \leq A |\lambda|^{-\frac{k}2}.
 $$
 }
 \label{ThLitt}
\end{thm}
As in \cite{FCcompl}, we emphasize that the use of the Littman theorem only improves the coefficients: indeed, had we not used this result (but all the other arguments from \cite{FCPAA} enhanced with the Riesz-Thorin theorem and the method from the previous section), we would have ended with the following alternative (and very close) proposition:
\begin{prop}
 \sl{For any $d\in \R$, $r> 2$, $q\geq 1$, $\theta\in[0,1]$ and $p\in[1, \frac8{\theta (1-\frac{2}{r})}]$, there exists a constant $C=C_{p,r,\theta}$ such that for any $f$ solving \eqref{systdispb} for initial data $f_0$ and external force $\Fe$ both with zero divergence and vorticity (that in the kernel of $\bP_2$), then
 \begin{equation}
  \||D|^d f\|_{\tilde{L}_t^p\dot{B}_{r, q}^0} \leq \frac{C_{p,r,\theta}}{\nu^{\frac{1}{p}-\frac{\theta}8(1-\frac{2}{r})}} \ee^{\frac{\theta}{8}(1-\frac{2}{r})} \left( \|f_0\|_{\dot{B}_{2, q}^{\sigma_2}} +\|\Fe\|_{\tilde{L}_t^1 \dot{B}_{2, q}^{\sigma_2}} \right),
 \end{equation}
 where $\sigma_2= d+\frac32-\frac{3}{r}-\frac{2}{p}+\frac{\theta}{4} (1-\frac{2}{r})$.
\label{Estimdispnualt}
 }
\end{prop}
Let us continue the proof of Proposition \ref{Estimdispnu} and first assume that $\Fe=0$. As explained in \cite{FCPAA, FCcompl}, as $div f_0=0=\omega(f_0)$, we have $f_0=\bP f_0= \bP_3 f_0 +\bP_4 f_0$, and:
$$
f(t)=\mathcal{F}^{-1}\left(e^{-\nu t |\xi|^2 +i\frac{t}{\ee} \frac{|\xi_h|}{|\xi|}} \mathcal{P}_3(\xi,\ee) \hat{f_0}(\xi) +e^{-\nu t |\xi|^2 -i\frac{t}{\ee} \frac{|\xi_h|}{|\xi|}} \mathcal{P}_4(\xi,\ee) \hat{f_0}(\xi) \right).
$$
Thanks to the orthogonality properties when $\nu=\nu'$ (the projectors are of norm 1), we are reduced to the study of:
$$
f(t)=\mathcal{F}^{-1}\left(e^{-\nu t |\xi|^2 +i\frac{t}{\ee} \frac{|\xi_h|}{|\xi|}} \hat{f_0}(\xi) \right),
$$
If $\varphi$ is the truncation function introduced in the Littlewood-Paley decomposition supported in the annulus centered at zero of radii $\frac34$ and $\frac83$ (see appendices from \cite{FCPAA, FCcompl} or \cite{Dbook} for a complete presentation), we denote by $\varphi_1$ another smooth truncation function, with support in a slightly larger annulus than $\mbox{supp }\varphi$ (for instance the annulus centered at zero and of radii $\frac12$ and $3$) and equal to $1$ on $\mbox{supp }\varphi$.
\\

With the same set $\cA$ and the same steps as in the previous section, for any $j\in \Z$ and $r \geq 1$:
\begin{multline}
 \|\ddj f\|_{L^p L^r}= \sup_{\psi \in \cA} \int_0^\infty \int_{\R^3} \ddj f(t,x) \psi(t,x) dx dt\\
 =C \sup_{\psi \in \cA} \int_0^\infty \int_{\R^3} e^{-\nu t|\xi|^2+i\frac{t}{\ee}\frac{|\xi_h|}{|\xi|}} \widehat{\ddj f_0}(\xi) \varphi_1(2^{-j} \xi) \hat{\psi}(t,\xi) d\xi dt\\
 \leq C \sup_{\psi \in \cA} \|\ddj f_0\|_{L^2} \left[\int_0^\infty \int_0^\infty \|L_j(\frac{t-t'}{\ee})\psi(t,.)\|_{L^r} \|e^{\nu(t+t')\Delta} \varphi_1(2^{-j}D) \overline{\psi(t',.)}\|_{L^{\bar{r}}} dtdt'\right]^{\frac12},
\label{estimTT1}
 \end{multline}
with $L_j(\sigma)$ defined as follows:
\begin{equation}
 L_j(\sigma)g=\int_{\R^3} e^{ix\cdot \xi +i \sigma \frac{|\xi_h|}{|\xi|}} \varphi_1(2^{-j}|\xi|) \hat{g}(\xi) d\xi =K_j(\sigma)*g,
\label{defLj}
\end{equation}
where
\begin{equation}
K_j(\sigma)(x)=\int_{\R^3} e^{ix\cdot \xi +i \sigma \frac{|\xi_h|}{|\xi|}} \varphi_1(2^{-j}|\xi|) d\xi.
 \label{defKj}
\end{equation}
Thanks to the frequency truncation (remember that $\mbox{supp }\varphi_1 \subset \cC(0,\frac12,3)$) and the classical estimates for the heat kernel in this case (we refer for example to Lemma 2.3 from \cite{Dbook}) we easily get that:
\begin{equation}
 \|e^{\nu(t+t')\Delta} \varphi_1(2^{-j}D) \overline{\psi(t',.)}\|_{L^{\bar{r}}} \leq C'e^{-\frac{\nu}4 (t+t')2^{2j}} \|\psi(t')\|_{L^{\bar{r}}}.
\label{Estimphi}
 \end{equation}
To bound the other term we will, as in \cite{FCcompl}, successively bound its $L^2 \rightarrow L^2$ and $L^1 \rightarrow L^\infty$ norms, and conclude thanks to the Riesz-Thorin theorem. Thanks to the Plancherel formula there exists a constant $C$ (only depending on $\varphi_1$) such that :
\begin{equation}
 \|L_j(\sigma) g\|_{L^2} \leq C \|g\|_{L^2}.
 \label{EstimLjL2}
\end{equation}
Thanks to the Young estimates:
\begin{equation}
 \|L_j(\sigma)g\|_{L^\infty}\leq \|K_j(\sigma)\|_{L^\infty}\|g\|_{L^1},
 \label{EstimLjL1}
\end{equation}
and performing the change of variable $\xi=2^j \eta$, we get $K_j(\sigma)(x)=2^{3j} K_0(\sigma)(2^j x)$ and:
\begin{equation}
 \|K_j(\sigma)\|_{L^\infty}\leq 2^{3j}\|K_0(\sigma)\|_{L^\infty},
\label{Chgtvar}
 \end{equation}
And as in \cite{FCcompl, KLT, IMT, LT} we need bound the $L^\infty$-norm of $K_0$ (see \eqref{defKj}). First we rewrite it as follows, denoting $b(\xi)=\frac{|\xi_h|}{|\xi|}$:
\begin{equation}
 K_0(\sigma)(x)=\int_{\R^3} e^{ix\cdot \xi +i \sigma b(\xi)} \varphi_1(|\xi|) d\xi,
\end{equation}
Compared to the case of the Primitive system (see \cite{IMT, FCcompl} the problem is here that the derivatives of $b$ present singularities when $\xi_h=(0,0)$. Using the Littman theorem then requires frequency cut-offs, as done in Proposition 3.1 from \cite{LT}:
\begin{prop} \cite{LT},
 \sl{With the previous notations, there exists a constant $C>0$ such that for all $x\in\R^3$ and all $\sigma\neq 0$,
\begin{equation}
 \|K_0(\sigma)\|_{L^\infty}\leq C(1+|\sigma|)^{-\frac12} \leq C \min(1, |\sigma|^{-\frac12}).
\label{ThLitt2}
\end{equation}
 }
\end{prop}
Had we not the singularity problem, we would simply perform the frequency truncations (as in \cite{KLT, IMT, LT, FCcompl}) to split the integral into several parts, on which the Littman theorem can be applied. The decomposition corresponds to zones where we precisely know how many nonzero eigenvalues are featured by the Hessian $D^2 b$ which writes in our case:
\begin{equation}
D^2 b(\xi)=\frac1{|\xi_h|^3|\xi|^5} \left(
 \begin{array}{ccc}
 \vspace{0.2cm}
  \xi_3^2 \big(\xi_2^2 |\xi|^2 -3 \xi_1^2 |\xi_h|^2\big) & -\xi_1 \xi_2 \xi_3^2 \big(|\xi|^2 +3 |\xi_h|^2\big) & \xi_1 \xi_3 |\xi_h|^2 \big(3|\xi_h|^2 -|\xi|^2\big)\\
  \vspace{0.2cm}
  -\xi_1 \xi_2 \xi_3^2 \big(|\xi|^2 +3 |\xi_h|^2\big) & \xi_3^2 \big(\xi_1^2 |\xi|^2 -3 \xi_2^2 |\xi_h|^2\big) & \xi_2 \xi_3 |\xi_h|^2 \big(3|\xi_h|^2 -|\xi|^2\big)\\
  \xi_1 \xi_3 |\xi_h|^2 \big(3|\xi_h|^2 -|\xi|^2\big) & \xi_2 \xi_3 |\xi_h|^2 \big(3|\xi_h|^2 -|\xi|^2\big) & -|\xi_h|^4 \big(3|\xi_h|^2 -2|\xi|^2\big)
 \end{array} \right),
\end{equation}
To ease the understanding of what is done in \cite{KLT, IMT, LT} (where the proofs are done only reasoning with the determinant) we compute the eigenvalues of this matrix (as we did in \cite{FCcompl}):
\begin{equation}
 \{\frac{\xi_3^2}{|\xi_h| |\xi|^3}, \quad -\frac{|\xi_h|\pm \sqrt{|\xi|^2+3\xi_3^2}}{2|\xi|^3}\}.
\end{equation}
Notice that the behaviour of the eigenvalues is a little different: if $\xi_3=0$, then only one eigenvalue is nonzero (and is equal to $-\frac{1}{|\xi_h|^2}$) and if $\xi_3\neq 0$ then none of the eigenvalues is zero. The idea is basically to split the integral into:
\begin{multline}
 K_0(\sigma)(x)=\int_{\R^3} e^{ix\cdot \xi +i \sigma b(\xi)} \chi(|\xi_3|)\varphi_1(|\xi|) d\xi +\int_{\R^3} e^{ix\cdot \xi +i \sigma b(\xi)}\big(1-\chi(|\xi_3|)\big) \Big) \varphi_1(|\xi|) d\xi\\
 = K_{0,1}(\sigma)(x)+K_{0,2}(\sigma)(x),
\end{multline}
Thanks to the study of $D^2 b$, at least 1 eigenvalue is nonzero for the first term, and the three of them are non zero for the second one, and if we could directly apply Theorem \ref{ThLitt} we would end-up with:
$$
\|K_{0,1}(\sigma)\|_{L^\infty} \leq C \min(1, |\sigma|^{-\frac12}), \quad \|K_{0,2}(\sigma)\|_{L^\infty} \leq C \min(1, |\sigma|^{-\frac32}).
$$
But this is not so simple as function $b$ is not regular on $\mbox{supp }\varphi_1$. In \cite{LT} the authors overcome this difficulty with frequency truncations.
\begin{rem}
 \sl{We can complete Remark 9 from \cite{FCcompl}: in the case of the rotating fluids, we decompose the frequency space into three zones, on two of them at least two eigenvalues are nonzero, on the third zone, every eigenvalue are nonzero. In the case of the Primitive system, on one zone at least one eigen value is nonzero, at least two of them in the second zone and all of them in the third zone. For the Stratified Boussinesq system, singularities arise and on the first zone, at least one eigenvalue is non zero, all of them in the other zone. We also refer to \cite{GenTa} who revisit Strichartz estimates for the linearized systems of the three models thanks to restriction theory (without resorting to dispersive estimates).}
\end{rem}
Gathering \eqref{ThLitt2}, \eqref{EstimLjL1} and \eqref{Chgtvar} leads to
$$
\|L_j(\sigma)g\|_{L^\infty}\leq C_F 2^{3j}\min(1, |\sigma|^{-\frac12})\|g\|_{L^1},
$$
Doing as in \cite{FCcompl}, gathering the previous estimate with \eqref{EstimLjL2} and thanks to the  Riesz-Thorin theorem, we obtain that for all $r\in[2,\infty]$ and $\theta \in[0,1]$:
\begin{equation}
 \|L_j(\sigma)g\|_{L^r} \leq C_r \frac{2^{3j(1-\frac2{r})}}{|\sigma|^{\frac{\theta}2 (1-\frac2{r})}} \|g\|_{L^{\bar{r}}},
\end{equation}
so that, using also \eqref{Estimphi}, we can bound \eqref{estimTT1} and obtain:
$$
 \|\ddj f\|_{L^p L^r} \leq C_r 2^{\frac{3j}2 (1-\frac2{r})} \ee^{\frac{\theta}4 (1-\frac2{r})} \|\ddj f_0\|_{L^2} \sup_{\psi \in \cB} \left[\int_0^\infty \int_0^\infty \frac{h(t)h(t')}{|t-t'|^{\frac{\theta}2 (1-\frac2{r})}}  dtdt'\right]^{\frac12},
$$
with $h(t)=e^{-\frac{\nu}4 t 2^{2j}} \|\psi(t,.)\|_{L^{\bar{r}}}$. The last term is bounded exactly like in \cite{FCcompl} using the Hardy-Littlewood-Sobolev estimates (this is here that we need the condition $p \leq \frac{4}{\theta (1-\frac{2}{r})}$) and we finally obtain:
$$
 \|\ddj f\|_{L^p L^r} \leq \frac{C_{p,r,\theta}}{\nu^{\frac{1}{p}-\frac{\theta}{4}(1-\frac{2}{r})}} \ee^{\frac{\theta}4 (1-\frac{2}{r})} 2^{j(\frac32-\frac{3}{r} -\frac{2}{p} +\frac{\theta}2 (1-\frac{2}{r}))} \|\ddj f_0\|_{L^2}
$$
which leads to the result in the homogeneous case. The inhomogeneous case (i.-e. when $\Fe\neq 0$) easily follows applying the previous steps to the Duhamel term. $\blacksquare$

\begin{rem}
 \sl{\begin{enumerate}
      \item In the case $\nu=\nu'$, there is a difference between the power of $|t-t'|/\ee$ in the dispersive estimates produced with the Littman method compared to the non-stationnary phase method: respectively $-1/2$ and $-1/4$. This is also noticed in \cite{FCcompl, FCRF} in the case of the rotating fluids (respectively $-1$ and $-1/2$). In \cite{FCcompl} we also point out that for the primitive system both exponents are the same.
      \item When $\nu\neq \nu'$ or when anisotropic estimates are needed the non-stationnary phase method is preferred to the Littman method. We refer to \cite{FCStratif2} for more details.
     \end{enumerate}
}
\end{rem}

\subsubsection{Proof of Proposition \ref{Heatflow}}
\label{FinApp1}
We adapt here the proof of Lemma 2.3 from \cite{Dbook}. As in \cite{FCestimLp} we only present here what is new. Thanks to \eqref{PrR} we can write that:
$$
e^{t\D} u=\cF^{-1} \left(f_{\frac{r}2, 2R}(\xi) e^{-t|\xi|^2} \hat{u}(\xi) \right) =g(t,.) \star u,
$$
with
$$
g(t,x)=(2\pi)^{-3} \int_{\R^3} e^{i x\cdot \xi} f_{\frac{r}2, 2R}(\xi) e^{-t|\xi|^2} d\xi (2\pi)^{-3} \int_{\R^3} e^{i x\cdot \xi} \chi(\frac{|\xi|}{2\Re})\Big(1-\chi(\frac{|\xi_h|}{\re})\Big) e^{-t|\xi|^2} d\xi,
$$
so that $\|e^{t\D} u\|_{L^p} \leq \|g(t,.)\|_{L^1} \|u\|_{L^p}$ and in order to bound the $L_x^1$-norm of $g$ we use the same method as in \cite{Dbook} but will only perform twice the integrations by parts (we are in $\R^3$):
$$
 g(t,x)=(2\pi)^{-3} (1+|x|^2)^{-2} \int_{\R^3} e^{i x\cdot \xi} (I_d-\D_{\xi})^2 \left(\chi(\frac{|\xi|}{2\Re})\Big(1-\chi(\frac{|\xi_h|}{\re})\Big) e^{-t|\xi|^2} \right) d\xi.
$$
The rest of the proof is then a matter of computing the derivatives of $f_{\frac{r}2, 2R}(\xi) e^{-t|\xi|^2}$ using the fact that:
$$
\begin{cases}
\vspace{0.2cm}
 \D(fg)=(\D f)g +2\n f\cdot \n g + f\D g,\\
 \n(e^{-t|\xi|^2})= -2t e^{-t|\xi|^2} \xi,\quad \mbox{and}\quad \D(e^{-t|\xi|^2})= -2t (3-2t |\xi|^2)e^{-t|\xi|^2}.
\end{cases}
$$
To simplify the notations, let us put $K=f_{\frac{r}2, 2R}$. We have:
\begin{multline}
 (I_d-\D)^2(K(\xi) e^{-t|\xi|^2})= \Bigg[(I_d-\D)^2 K(\xi)+ 2t \Big(3 K(\xi)-10 \D K(\xi)+2\xi\cdot \n K(\xi)\Big)\\
 +4t^2 \Big(15 K(\xi)+20 \xi\cdot \n K(\xi)+4\xi\cdot(\xi\cdot \n)\n K(\xi)\Big)\\
 +8t^3 \Big(-10 |\xi|^2 K(\xi)-4|\xi|^2 \xi\cdot \n K(\xi)\Big) +16 t^4 |\xi|^4 K(\xi)\Bigg] e^{-t|\xi|^2}.
\end{multline}
We can get rid of the powers of $t$ with estimates such as $t^\aa |\xi|^{2\aa} e^{-t|\xi|^2} \leq C' e^{-\frac{t}2 |\xi|^2} \leq C' e^{-\frac{t}2 r^2}$ and the proof follows from the fact that there exists a constant $C>0$ such that for all $\xi \in \cC_{\frac{r}2, 2R}$, we have
$$
\begin{cases}
\vspace{0.2cm}
 \|\n K\| \leq C(\frac1{R}+\frac1{r}),\\
 \vspace{0.2cm}
 |\D K| \leq C(\frac1{R^2}+\frac1{Rr}+\frac1{r^2}),\\
 \vspace{0.2cm}
 \|(\xi \cdot \n) \n K\| \leq C(\frac1{R}+\frac1{r}+\frac{R}{r^2}),\\
|\D^2 K| \leq C (\frac1{R^4}+\frac1{R^3 r}+\frac1{R^2 r^2}+\frac1{R r^3}+\frac1{r^4}),
\end{cases}
$$
so that, integrating on $\cC_{\frac{r}2, 2R}$, we finally obtain:
$$
|g(t,x)|\leq \frac{C}{(1+|x|^2)^2} \frac{R^3}{r^4} e^{-\frac{t}2 r^2} \quad \mbox{and} \quad
\|g(t,.)\|_{L^1} \leq C \frac{R^3}{r^4} e^{-\frac{t}2 r^2},
$$
which concludes the proof. $\blacksquare$

\section{Appendix 2}

For $0<\aa<R$, and $\bb\geq 0$, let us define, for any $x\in \R$,
$$
f_\aa(x)=\frac{\aa x}{(x^2+\aa^2)^\frac32},
$$
and
\begin{equation}
  I_{\aa, \bb}^R(\sigma) \overset{def}{=} \int_0^{\sqrt{R^2-\aa^2}} \frac{dx}{1+\sigma(f_\aa(x)-\bb)^2},
 \label{DefIab}
\end{equation}
It is obvious that $\IbR \leq R$ but we wish to bound this integral when $\sigma$ is large. The aim of this section is to prove the following proposition, that is crucial to obtain the Strichartz estimates from Proposition \ref{Estimdispnu2}:

\begin{prop}
 \sl{There exists a constant $C_0>0$ such that for any $\aa>0$, $R\geq \frac2{\sqrt{3}} \aa$,
 \begin{equation}
  \sup_{\bb \in \R_+} I_{\aa, \bb}^R(\sigma) \leq C_0 \frac{R^7}{\aa^\frac{11}2} \min(1,\sigma^{-\frac14}).
 \end{equation}
 Moreover, the exponent $-\frac14$ is optimal in the sense that there exist $c_0,\sigma_0>0$ such that for any $R\geq \frac{\sqrt{3}}{\sqrt{2}} \aa$ and $\sigma\geq\sigma_0$,
 $$
 \sup_{\bb \in \R_+} I_{\aa, \bb}^R(\sigma) \geq c_0\sigma^{-\frac14} \aa^{\frac32}.
 $$
 }
 \label{PropIab}
\end{prop}

\subsection{Reduction of the problem}

If $\aa>0$, it is immediate that for any $\lambda>0$ and $x\in \R$, $f_\aa(\lambda x)=\frac1{\lambda} f_{\frac{\aa}{\lambda}} (x)$, so that in particular, we get that for any $x\in \R$, $f_\aa(x)=\frac1{\aa} f_1(\frac{x}{\aa})$. Performing the change of variable $x=\aa y$, we obtain:
$$
I_{\aa, \bb}^R(\sigma)=\aa I_{1,\aa \bb}^\frac{R}{\aa}(\frac{\sigma}{\aa^2}),
$$
so that we are reduced, given $\bb\in \R$ and $R>1$ to prove that
\begin{prop}
 \sl{There exists a constant $C_0>0$ such that for any $R\geq \frac2{\sqrt{3}}$,
 \begin{equation}
  \sup_{\bb \in \R_+} \IbR =\sup_{\bb \in \R_+}  \int_0^{\sqrt{R^2-1}} \frac{dx}{1+\sigma(f_1(x)-\bb)^2} \leq C_0 R^7 \min(1,\sigma^{-\frac14}).
 \end{equation}
 Moreover, the exponent $-\frac14$ is optimal in the sense that there exist $c_0,\sigma_0>0$ such that for any $R\geq \frac{\sqrt{3}}{\sqrt{2}}$ and $\sigma\geq\sigma_0$,
 $$
 \sup_{\bb \in \R_+} \IbR \geq \IbRd \geq c_0 \sigma^{-\frac14}.
 $$
 }
 \label{PropIab2}
\end{prop}

\subsection{First properties of $f_1$ and $\IbR$}

In what follows we will need to study not only $f_1$ but also its first three derivatives. We easily obtain that for all $x\geq 0$:
$$
\begin{cases}
\vspace{0.3cm}
f_1(x)=\frac{x}{(x^2+1)^\frac32}, \quad f_1'(x)= \frac{1-2x^2}{(x^2+1)^\frac52}, \quad f_1'' (x)=\frac{3x(2x^2-3)}{(x^2+1)^\frac72},\\
f_1^{(3)} = \frac{-3}{(x^2+1)^\frac92}(8x^4-24 x^2+3), \quad f_1^{(4)} = \frac{15x}{(x^2+1)^\frac{11}2}(8x^4-40 x^2+15).
\end{cases}
$$
An elementary study gives us the following variations:
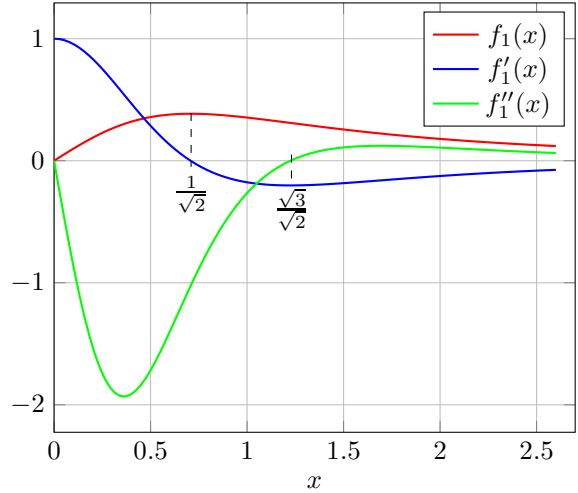
\begin{figure}[H]
		\centering
		\begin{tikzpicture}
			\begin{axis}
				[xmin=0,xmax=2.7,
				xlabel=\(x\),
				samples=100,
				grid=both
				]
				\addplot[domain=0:2.6,color=red,thick]{x/((x^2+1)^(1.5))};
				\addplot[domain=0:2.6,color=blue,thick]{(1-2*x^2)/((x^2+1)^(2.5))};
				\addplot[domain=0:2.6,color=green,thick]{3*x*(2*x^2-3)/((x^2+1)^(3.5))};
				\draw [dashed] (0.71,0.39) -- (0.71,-0.05) node[below] {$\frac1{\sqrt{2}}$};
				\draw [dashed] (1.23,0.05) -- (1.23,-0.15) node[below] {$\frac{\sqrt{3}}{\sqrt{2}}$};
				\legend{\(f_1(x)\), \(f_1'(x)\), \(f_1''(x)\)};
			\end{axis}
		\end{tikzpicture}
		\caption{Variations of $f_1$ and derivatives}
	\end{figure}
More precisely, $f_1$ is bijective from $]0,\frac1{\sqrt{2}}[$ to $]0, \frac2{3\sqrt{3}}[$, and from $]\frac1{\sqrt{2}}, \infty[$ to $]0, \frac2{3\sqrt{3}}[$. Hence for any $y\in ]0, \frac2{3\sqrt{3}}[$ the equation $f_1(x)=y$ admits exactly two solutions denoted as:
\begin{equation}
 0<z_1(y)<\frac1{\sqrt{2}} <z_2(y).
\end{equation}
\begin{figure}[H]
		\centering
		\begin{tikzpicture}
			\begin{axis}
				[xmin=0,xmax=3.4,
				ymin=-0.06, ymax=0.42,
				xlabel=\(x\),
				samples=100,
				grid=both
				]
				\addplot[domain=0:3.4,color=red,thick]{x/((x^2+1)^(1.5))};
				\draw [dashed,thick] (0.71,0.39) -- (0.71,0) node[below] {$\frac1{\sqrt{2}}$};
				\draw [dashed] (0.26,0.22) -- (0.26,0) node[below] {$z_1(y)$};
				\draw [dashed] (1.7,0.22) -- (1.7,0) node[below] {$z_2(y)$};
				\draw [dashed] (1.7,0.22) -- (0.26,0.22) node[left] {$y$};
				\legend{\(f_1(x)\), \(f_1'(x)\), \(f_1''(x)\)};
			\end{axis}
		\end{tikzpicture}
		\caption{Solutions of $f_1(x)=y$}
\end{figure}
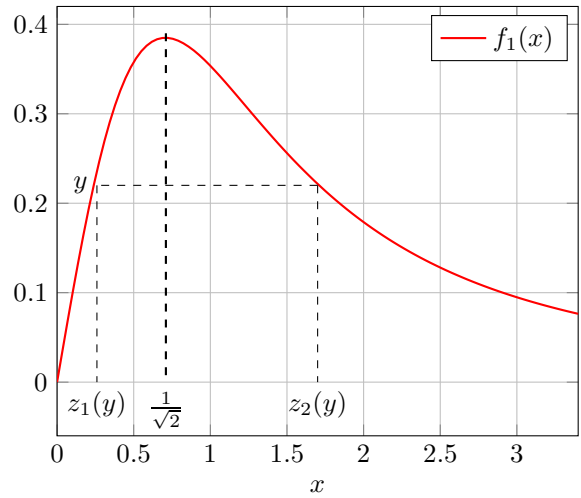

\begin{rem}
 \sl{
 \begin{enumerate}
  \item Both $0$ and $\frac2{3\sqrt{3}}$ have only one antecedent by $f_1$, respectively $0$ and $\frac1{\sqrt{2}}$.
  \item We will abundantly use in what follows that $z_1,z_2$ are the respective inverse functions of $f_{1/]0, \frac2{3\sqrt{3}}[}$, and $f_{1/]\frac2{3\sqrt{3}}, \infty[}$.
 \end{enumerate}
}
\end{rem}
\begin{prop}
 \sl{The previous functions can be explicitely expressed as follows. For any $y\in ]0, \frac2{3\sqrt{3}}[$, we have
 \begin{equation}
  \begin{cases}
  \vspace{0.2cm}
   z_1(y)=\sqrt{\frac2{y\sqrt{3}} \cos \left[\frac13 \arccos(-\frac{3 \sqrt{3}}2 y) +\frac{4\pi}3\right] -1},\\
   z_2(y)=\sqrt{\frac2{y\sqrt{3}} \cos \left[\frac13 \arccos(-\frac{3 \sqrt{3}}2 y)\right] -1}
  \end{cases}
 \end{equation}
 }
 \label{z12expr}
\end{prop}
\textbf{Proof:} for any $x>0$ and $y\in ]0, \frac2{3\sqrt{3}}[$, $f_1(x)=y \Leftrightarrow \left(\frac{x}{y}\right)^2=(x^2+1)^3 \Leftrightarrow$ the new variable $X\overset{def}{=}x^2+1$ satisfies:
\begin{equation}
 X^3+pX+q=0,
 \label{Eqdeg3}
\end{equation}
where $p=-\frac1{y^2}=-q$, which is exactly the reduced form featured in the statement of the Cardan formulae. We then follow the well-known formulae and first compute the discriminant:
$$
\D =4p^3+27q^2= -\frac4{y^6}+\frac{27}{y^4}.
$$
As $y\in ]0, \frac2{3\sqrt{3}}[$, $\D<0$ and \eqref{Eqdeg3} admits exacly three solutions given by the following expressions: if $\phi$ is such that $\cos (3\phi) =\frac32 \frac{q}{p}\sqrt{-\frac3{p}}$, then
\begin{equation}
\begin{cases}
\vspace{0.2cm}
 X_1(y)=\sqrt{-\frac{4p}3} \cos(\phi) =\frac2{y\sqrt{3}} \cos\left[\frac13 \arccos(-\frac{3 \sqrt{3}}2 y) \right],\\
 \vspace{0.2cm}
 X_2(y)=\sqrt{-\frac{4p}3} \cos(\phi +\frac{2\pi}3) =\frac2{y\sqrt{3}} \cos\left[\frac13 \arccos(-\frac{3 \sqrt{3}}2 y) +\frac{2\pi}3\right],\\
 X_3(y)=\sqrt{-\frac{4p}3} \cos(\phi +\frac{4\pi}3) =\frac2{y\sqrt{3}} \cos\left[\frac13 \arccos(-\frac{3 \sqrt{3}}2 y) +\frac{4\pi}3\right].
\end{cases}
\end{equation}
We recall that for any given $y\in ]0, \frac2{3\sqrt{3}}[$ the only solutions $0<z_1(y)<\frac1{\sqrt{2}} <z_2(y)$ of the equation $f_1(x)=y$ satisfy
$$
\{1+z_1(y)^2, 1+z_2(y)^2\} \subset \{X_1(y),X_2(y),X_3(y)\},
$$
and thanks to the fact that $\phi=\frac13 \arccos(-\frac{3 \sqrt{3}}2 y)$, we get that for any $y\in ]0, \frac2{3\sqrt{3}}[$,
$$
 \cos(\phi) \in \Big]\frac12,\frac{\sqrt{3}}2\Big[,\quad \cos(\phi +\frac{2\pi}3) \in \Big]-\frac{\sqrt{3}}2,-1\Big[, \quad \cos(\phi +\frac{4\pi}3) \in \Big]0,\frac12\Big[,
$$
which immediately implies that $X_2<0<X_3<X_1$ so that we are sure that
$$
(1+z_1(y)^2, 1+z_2(y)^2) =(X_3(y),X_1(y)),
$$
which concludes the proof. $\blacksquare$
\\
We can immediately state the following elementary properties:
\begin{prop}
 \sl{With the previous notations,
 \begin{enumerate}
  \item For any $\bb>\frac2{3\sqrt{3}}$,
  $$
  \IbR \leq \frac{R}{1+\sigma (\bb-\frac2{3\sqrt{3}})^2} \leq R \min\left(1, \sigma^{-1} (\bb-\frac2{3\sqrt{3}})^{-2}\right),
  $$
  \item There exists $C_0>0$ such that for any $R>1$ and $\sigma>0$, $I_{1,0}^R(\sigma)\leq C_0 \sigma^{-\frac12} R^3$.
  \item For any $\bb,\bb_0\geq 0$,
  \begin{equation}
   |\IbR-\IbRo| \leq \left(\sigma^\frac12|\bb-\bb_0| +\sigma |\bb-\bb_0|^2\right) \min\left(\IbR, \IbRo\right).
  \end{equation}
In particular, for any $\bb$ we have $|\bb-\bb_0|\leq \sigma^{-\frac12} \Rightarrow \IbR \leq 3 \IbRo$.
 \end{enumerate}
 }
 \label{Propelem}
\end{prop}
\begin{rem}
 \sl{Point 1 is interesting when $\beta\geq \frac2{3\sqrt{3}}+ k$ (for some $k>0$) but useless when $\bb$ goes to $\frac2{3\sqrt{3}}$, Point 3 will be crucial.}
\end{rem}
\textbf{Proof:} The first point is immediate as, thanks to the bounds of $f_1$, for any $\bb>\frac2{3\sqrt{3}}$ and any $x\geq 0$, $|f_1(x)-\bb|=\bb -f_1(x) \geq \bb -\frac2{3\sqrt{3}}$. For the second point, we simply remark that for any $x\in[0,\sqrt{R^2-1}]$, $f_1(x) \geq x R^{-3}$ so that thanks to the change of variable $z=\sigma^\frac12 x R^{-3}$,
$$
I_{1,0}^R(\sigma) \leq \int_0^{\sqrt{R^2-1}} \frac{dx}{1+\sigma x^2 R^{-6}} \leq \sigma^{-\frac12} R^3 \int_0^\infty \frac{dz}{1+z^2}.
$$
To prove the last point, we begin with
$$
 \IbR -\IbRo =2\sigma (\bb-\bb_0)\int_0^{\sqrt{R^2-1}} \frac{f_1(x)-\frac12(\bb+\bb_0)}{(1+\sigma(f_1(x)-\bb)^2)(1+\sigma(f_1(x)-\bb_0)^2)} dx.
$$
Noticing that $f_1(x)-\frac12(\bb+\bb_0) =f_1(x)-\bb +\frac12(\bb-\bb_0)$ and using the fact that
\begin{equation}
 |f_1(x)-\bb|=\sigma^{-\frac12} (\sigma^\frac12 |f_1(x)-\bb|) \leq \frac{\sigma^{-\frac12}}2 (1+\sigma (f_1(x)-\bb)^2),
\label{Youngsigma}
\end{equation}
we get that
\begin{multline}
 |\IbR -\IbRo| \leq 2\sigma |\bb-\bb_0| \Bigg(\frac{\sigma^{-\frac12}}2 \int_0^{\sqrt{R^2-1}} \frac{dx}{1+\sigma(f_1(x)-\bb_0)^2} \\
 +\frac{|\bb-\bb_0|}2 \int_0^{\sqrt{R^2-1}} \frac{dx}{(1+\sigma(f_1(x)-\bb)^2)(1+\sigma(f_1(x)-\bb_0)^2)}\Bigg)\\
 \leq \sigma^\frac12|\bb-\bb_0| \IbRo +\sigma |\bb-\bb_0|^2 \int_0^{\sqrt{R^2-1}} \frac{dx}{(1+\sigma(f_1(x)-\bb)^2)(1+\sigma(f_1(x)-\bb_0)^2)}\\
 \leq \sigma^\frac12|\bb-\bb_0| \IbRo +\sigma |\bb-\bb_0|^2 \min(\IbR, \IbRo).
\end{multline}
As we also have $f_1(x)-\frac12(\bb+\bb_0) =f_1(x)-\bb_0-\frac12(\bb-\bb_0)$, we similarly obtain that:
$$
|\IbR -\IbRo| \leq \sigma^\frac12|\bb-\bb_0| \IbR +\sigma |\bb-\bb_0|^2 \min(\IbR, \IbRo),
$$
which leads to the result. The last result is an immediate consequence of the second point. $\blacksquare$

\subsection{Study when $\bb>0$ is small}
In this section we will need the following asymptotic expansions near zero.
\begin{prop}
 \sl{When $y>0$ is small we have:
 \begin{equation}
  \begin{cases}
  \vspace{0.2cm}
z_1(y)= y\left(1+\frac32 y^2+o(y^2)\right),\\
 \vspace{0.1cm}
z_2(y)=\frac1{y^\frac12}\left(1-\frac34 y-\frac{15}{32}y^2-\frac{77}{128}y^3+o(y^3)\right),\\
 \vspace{0.2cm}
f_1'(z_1(y))=1-\frac92 y^2-\frac{33}8 y^4+o(y^4),\\
f_1'(z_2(y))=-2 y^\frac32\left(1-\frac32 y -\frac{27}{32}y^2+o(y^2)\right).\\
  \end{cases}
 \end{equation}
\label{DL0}
 }
\end{prop}
We will then prove the following result:
\begin{prop}
 \sl{There exist $\delta_0>0$ and $C_{\delta_0}>0$ such that for any $\bb\in[0,\delta_0]$, any $R\geq \frac2{\sqrt{3}}$ and $\sigma>0$,
 $$
 \IbR \leq C_{\delta_0} R^7 \sigma^{-\frac12}.
 $$
 }
 \label{estimbetapetit}
\end{prop}
\textbf{Proof:} For $\bb>0$ ($\delta_0$ will be specified later) let us split $\IbR$ as follows
\begin{equation}
\IbR=\int_0^\frac1{\sqrt{2}} \frac{dx}{1+\sigma(f_1(x)-\bb)^2} +\int_\frac1{\sqrt{2}}^{\sqrt{R^2-1}} \frac{dx}{1+\sigma(f_1(x)-\bb)^2} \overset{def}{=} \Jb +\Kb.
\label{Split}
\end{equation}
The methods will be similar for both integrals but as $z_2(y)$ goes to infinity when $y>0$ goes to zero, the second one will require more discussion. Let us begin with splitting $\Jb$ into three parts as follows: for $\frac12\leq m<1<M \leq 2$ (that will be precised later) we define:
\begin{multline}
 \Jb= J_1 +J_2+ J_3\\
 \overset{def}{=} \int_0^{z_1(m\bb)} \frac{dx}{1+\sigma(f_1(x)-\bb)^2} +\int_{z_1(m\bb)}^{z_1(M\bb)} \frac{dx}{1+\sigma(f_1(x)-\bb)^2} +\int_{z_1(M\bb)}^\frac1{\sqrt{2}} \frac{dx}{1+\sigma(f_1(x)-\bb)^2}.
 \label{DecompoJ}
\end{multline}
As $f_1$ is strictly increasing on $[0,\frac1{\sqrt{2}}]$, and by the definition of $z_1$, for all $x\in [0,z_1(m\bb)]$, we have $0\leq f_1(x) \leq m\bb$, so that $\bb\geq \frac1{m} f_1(x)$ and (also using that $x\leq \frac1{\sqrt{2}}$):
$$
\bb-f_1(x) \geq \frac{1-m}m f_1(x) \geq \frac{1-m}m \frac{2\sqrt{2}}{3\sqrt{3}} x \geq C_0 (1-m) x\geq 0,
$$
so that, thanks to the change of variable $z= C_0 \sigma^\frac12 (x-m) x$ (and constantly denoting as $C_0$ the constant, that may vary from line to line),
\begin{equation}
 J_1 \leq \int_0^{z_1(m\bb)} \frac{dx}{1+\sigma(C_0 (1-m) x)^2}\leq C_0 \frac{\sigma^{-\frac12}}{1-m} \int_0^\infty \frac{dz}{1+y^2} \leq \frac{C_0}{1-m} \sigma^{-\frac12}.
 \label{estimb0J1}
\end{equation}
Similarly, for any $x\in [z_1(M\bb), \frac1{\sqrt{2}}]$, $\bb\leq \frac1{M} f_1(x)$ and we obtain that:
\begin{equation}
 J_3 \leq \frac{C_0}{M-1} \sigma^{-\frac12}.
 \label{estimb0J3}
\end{equation}
All the work consists in correctly estimating $J_2$ and to do this we introduce (by definition, $\bb=f_1(\zib)$):
\begin{equation}
\Jbb \overset{def}{=} \int_{z_1(m\bb)}^{z_1(M\bb)} \frac{dx}{1+\sigma\big[f_1'(\zib)\big(x-\zib\big)\big]^2}.
 \label{defJ2b}
\end{equation}
Now we can write that:
\begin{multline}
 J_2-\Jbb=-\sigma \int_{z_1(m\bb)}^{z_1(M\bb)} \displaystyle{\frac{f_1(x)-f_1(\zib)-f_1'(\zib)(x-\zib)}{\left(1+\sigma\big(f_1(x)-\bb\big)^2 \right) \left(1+\sigma\big[f_1'(\zib)\big(x-\zib\big)\big]^2\right)}}\\
 \times \left[f_1(x)-f_1(\zib)+f_1'(\zib)\big(x-\zib\big) \right] dx.
\end{multline}
Thanks to the Taylor series expansion with integral remainder, and using that $f_1''$ is bounded by $2$, we have:
\begin{multline}
 |f_1(x)-f_1(\zib)-f_1'(\zib)\big(x-\zib\big)|=\Big|(x-\zib)^2 \int_0^1 (1-u) f_1''\big((1-u)\zib+ux\big) du\Big|\\
 \leq |x-\zib|^2.
 \label{Taylor2}
\end{multline}
Injecting in the previous lines, and using once more \eqref{Youngsigma}, we have:
\begin{multline}
 |J_2-\Jbb| \leq \sigma \int_{z_1(m\bb)}^{z_1(M\bb)} \frac{|x-\zib|^2 |f_1(x)-f_1(\zib)|+|f_1'(\zib)|\cdot|x-\zib|^3}{\left(1+\sigma(f_1(x)-\bb)^2 \right) \left(1+\sigma |f_1'(\zib)|^2 |x-\zib|^2\right)}dx\\
 \leq \sigma \int_{z_1(m\bb)}^{z_1(M\bb)} \left[\frac{\frac12\sigma^{-\frac12}}{\sigma |f_1'(\zib)|^2} +\frac{|x-\zib|}{\left(1+\sigma(f_1(x)-\bb)^2 \right) C_0\sigma |f_1'(\zib)|} \right]dx\\
 \leq \frac{\sigma^{-\frac12}}{2|f_1'(\zib)|^2} (z_1(M\bb)-z_1(m\bb))+\frac{\max\Big(z_1(M\bb)-z_1(\bb), z_1(\bb)-z_1(m\bb)\Big)}{|f_1'(\zib)|} J_2.
 \label{estimJJA}
\end{multline}
We need to estimate the differences between $z_1(M\bb),z_1(m\bb)$ and $z_1(\bb)$, but as $m<1<M$ will be intended to be very close to $1$, we won't use the asymptotic expansions from Proposition \ref{DL0} and simply write instead that as $z_1$ is the inverse function of $f_1$ (restricted to $]0,\frac1{\sqrt{2}}$):
$$
z_1(M\bb)-z_1(\bb)=\int_\bb^{M\bb} z_1'(t)dt =\int_\bb^{M\bb} \frac{dt}{f_1'(z_1(t))}.
$$
Thanks to the fact that $f_1'$ decreases on $[0,\frac{\sqrt{3}}{\sqrt{2}}]$ and satisfies $f_1'(0)=1$, and as $\frac12\leq m<1<M \leq 2$, we have $[m\bb, M\bb]\subset[0,2\bb]\subset[0,2\delta_0]$ so that for all $t\in[0,2\delta_0]$ (we recall that $z_1$ increases as $f_1$ does on $[0,\frac1{\sqrt{2}}]$),
\begin{equation}
f_1'(z_1(t))\geq \fdo>0.
\label{Minofp}
\end{equation}
Plugging this into the previous integral gives:
\begin{equation}
 z_1(M\bb)-z_1(\bb) \leq \frac{(M-1)\bb}{\fdo}.
\label{zM}
 \end{equation}
Similarly, we obtain that:
\begin{equation}
z_1(\bb)-z_1(m\bb) \leq \frac{(1-m)\bb}{\fdo}, \quad \mbox{and} \quad z_1(M\bb)-z_1(m\bb) \leq \frac{(M-m)\bb}{\fdo}.
\label{zmM}
\end{equation}
Returning to \eqref{estimJJA}, we obtain that (using \eqref{Minofp}, we also have $f_1'(z_1(\beta))\geq \fdo$):
\begin{equation}
 J_2 \leq \Jbb +\frac{\sigma^{-\frac12}}2 \frac{(M-m)\bb}{\fdo^3} +\frac{\beta}{\fdo^2}\max(M-1,1-m) J_2.
\label{estimJ2mM}
 \end{equation}
Considering $(m,M)=(1-k,1+k)$ (with $k>0$ small), as $\beta \leq \delta_0\leq 1$, the previous estimates turns into:
$$
 J_2 \leq \Jbb +\sigma^{-\frac12} \frac{k}{\fdo^3} +\frac{k}{\fdo^2} J_2 \leq \Jbb +\frac{\sigma^{-\frac12}}{2\fdo} +\frac12 J_2,
$$
when we choose $k= \frac12\min(1, \fdo^{-2})>0$. This entails that:
\begin{equation}
 J_2 \leq 2\Jbb +\frac{\sigma^{-\frac12}}{\fdo}
 \label{estimJJA2}
\end{equation}
To conclude, we need to estimate $\Jbb$ (see \eqref{defJ2b}) which is easy thanks to the change of variable $z=\sigma^{-\frac12} f_1'(\zib)(x-\zib)$:
\begin{multline}
 \Jbb = \frac{\sigma^{-\frac12}}{f_1'(\zib)} \int_{\sigma^{\frac12} f_1'(\zib) (z_1(m\bb)-\zib)}^{\sigma^{\frac12} f_1'(\zib) (z_1(M\bb)-\zib)} \frac{dz}{1+z^2}
 \leq \frac{\sigma^{-\frac12}}{\fdo} \int_{-\infty}^\infty \frac{dz}{1+z^2}\\
 \leq 2\pi\frac{\sigma^{-\frac12}}{\fdo}.
 \label{estimJbar2}
\end{multline}
Plugging this into \eqref{estimJJA2}, we obtain that:
$$
J_2\leq (4\pi+1)\frac{\sigma^{-\frac12}}{\fdo}.
$$
Gathering this with \eqref{estimb0J1}, \eqref{estimb0J3} and recalling that $k= \frac12\min(1, \fdo^{-2})>0$, we get that
\begin{equation}
 \Jb \leq \left(4C_0 \max(1,\fdo^2) +\frac{4\pi+1}{\fdo}\right) \sigma^{-\frac12}.
 \label{estimJJA3}
\end{equation}
\begin{rem}
 \sl{We emphasize that the constant is harmless as when $\delta_0$ goes to zero, it goes to $4C_0+4\pi+1$.}
\end{rem}
Now we continue with the second term $\Kb$: the methods will be similar except that we now write $\beta=f_1(z_2(\beta))$ where $z_2$ is the inverse function of $f_1$ restricted to $[\frac1{\sqrt{2}}, \infty[$ (the function is now decreasing as well as $z_2$).

We will also introduce $\frac12\leq m<1<M \leq 2$ but we will need to discuss the position of $\sqrt{R^2-1}$ relatively to the points:
$$
\frac1{\sqrt{2}} <z_2(2\bb) < z_2(\bb) < z_2(\frac{\bb}2),
$$
and will distinguish the following three cases:
\begin{enumerate}
 \item $\sqrt{R^2-1} <z_2(2\bb) < z_2(\bb) \quad \Longleftrightarrow \quad \bb <\frac12 f_1(\sqrt{R^2-1})$,
 \item $z_2(2\bb) \leq \sqrt{R^2-1} <z_2(\frac{\bb}2)  \quad \Longleftrightarrow \quad \frac12  f_1(\sqrt{R^2-1}) \leq \bb < 2 f_1(\sqrt{R^2-1})$,
 \item $z_2(2\bb) < z_2(\bb) \leq \sqrt{R^2-1}  \quad \Longleftrightarrow \quad 2 f_1(\sqrt{R^2-1}) \leq \bb$.
\end{enumerate}
The first case is easy: as $f_1$ decreases, for all $x\in[\frac1{\sqrt{2}}, \sqrt{R^2-1}]$ we have $f_1(x) \geq f_1(\sqrt{R^2-1}) >2\bb$ and $\bb<\frac12 f_1(x)$ so that
$$
|\bb-f_1(x)| =f_1(x)-\bb > \frac12 f_1(x) >0.
$$
In addition as $x\leq \sqrt{R^2-1}$, $f_1(x)\geq \frac{x}{R^3}$ so that using the change of variable $z=\frac{\sigma^{\frac12}x}{2R^3}$, we get:
\begin{equation}
 \Kb \leq \int_{\frac1{\sqrt{2}}}^{\sqrt{R^2-1}} \frac{dx}{1+\frac{\sigma}4 \frac{x^2}{R^6}} \leq 2R^3 \sigma^{-\frac12} \int_0^\infty \frac{dy}{1+y^2} \leq \pi R^3 \sigma^{-\frac12}.
\end{equation}
Let us jump to the third case: as $z_2(\bb) \leq \sqrt{R^2-1}$ we can reproduce the arguments used for $\Jb$ and split the integral in three according to:
\begin{multline}
 \Kb =K_1 +K_2+ K_3\\
 \overset{def}{=} \int_{\frac1{\sqrt{2}}}^{z_2(M\bb)} \frac{dx}{1+\sigma(f_1(x)-\bb)^2} +\int_{z_2(M\bb)}^{z_2(m\bb)} \frac{dx}{1+\sigma(f_1(x)-\bb)^2} +\int_{z_2(m\bb)}^{\sqrt{R^2-1}} \frac{dx}{1+\sigma(f_1(x)-\bb)^2}.
 \label{DecompoK}
\end{multline}
For $K_1$ and $K_3$ the methods are the same as for $\Jb$ except that, as $x\in[\frac1{\sqrt{2}}, \sqrt{R^2-1}]$, we can only write that $f_1(x) \geq \frac{x}{R^3}$ which leads to:
\begin{equation}
 \begin{cases}
 \vspace{0.2cm}
  K_1 \leq \frac{\pi}{M-1} R^3 \sigma^{-\frac12},\\
  K_3 \leq \frac{\pi}{2(1-m)} R^3 \sigma^{-\frac12}.
 \end{cases}
\label{K1K3}
\end{equation}
For $K_2$, introducing (now, by definition, $\bb=f_2(\ziib)$):
\begin{equation}
\Kbb \overset{def}{=} \int_{z_2(M\bb)}^{z_2(m\bb)} \frac{dx}{1+\sigma[f_1'(\ziib)(x-\ziib)]^2},
 \label{defK2b}
\end{equation}
we similarly obtain that:
\begin{multline}
 |K_2-\Kbb| \leq \frac{\sigma^{-\frac12}}{2|f_1'(\ziib)|^2} (z_2(m\bb)-z_2(M\bb))\\
 +\frac{\max\Big(z_2(m\bb)-\ziib, \ziib-z_2(M\bb)\Big)}{|f_1'(\ziib)|} K_2.
 \label{estimKKA}
\end{multline}
Also similarly we have (we recall that now, $f_1$ and $z_1$ are decreasing):
\begin{equation}
 z_2(\bb)-z_1(M\bb)=\int_{M\bb}^\bb z_1'(t)dt =\int_{M\bb}^\bb \frac{dt}{f_1'(z_1(t))} =\int_\bb^{M\bb} \frac{dt}{-f_1'(z_1(t))} =\int_\bb^{M\bb} \frac{dt}{|f_1'(z_1(t))|}.
\label{estimzM2}
\end{equation}
We emphasize that now $f_1'$ is negative and decrases on $[\frac1{\sqrt{2}}, \frac{\sqrt{3}}{\sqrt{2}}]$, increases on $[\frac{\sqrt{3}}{\sqrt{2}}, \infty]$, and is small, this is here that we will need the expansions from Proposition \ref{DL0}. Noticing that:
$$
z_2(2\beta) \geq \frac{\sqrt{3}}{\sqrt{2}}\quad \Longleftrightarrow \quad f_1(z_2(2\beta)) \leq f_1(\frac{\sqrt{3}}{\sqrt{2}}) \quad  \Longleftrightarrow \quad \beta \leq \frac{\sqrt{3}}{5\sqrt{5}},
$$
so that we can reduce to the case where $|f_1'|$ is decreasing on $[z_2(2\bb), z_2(\frac{\bb}2)]$ asking that $\delta_0\leq \frac{\sqrt{3}}{5\sqrt{5}}$. Now we can state that for all $t\in[m\bb, M\bb]\subset[\frac{\bb}2,2\bb]$,
\begin{equation}
 \frac1{|f_1'(z_2(2\bb))|} \leq \frac1{|f_1'(z_2(t))|} \leq \frac1{|f_1'(z_2(\frac{\bb}2))|}.
 \label{Minofp2}
\end{equation}
Thanks to Proposition \ref{DL0}, we know that
$$
\frac{|f_1'(z_2(\frac{\bb}2))|}{\bb^\frac32} \underset{\bb \rightarrow 0}{\longrightarrow} \frac1{\sqrt{2}},
$$
and we will ask in addition that $\delta_0 \in]0, \frac{\sqrt{3}}{5\sqrt{5}}]$ is so small that for all $\beta\in]0,\delta_0]$ (this is here that the size of $\delta_0$ is specified):
\begin{equation}
 |f_1'(z_2(\frac{\bb}2))| \geq \frac12 \bb^\frac32.
 \label{Minofp3}
\end{equation}
Now we can go back to \eqref{estimzM2} and obtain that:
\begin{equation}
 \begin{cases}
 \vspace{0.2cm}
  z_2(\bb)-z_1(M\bb) \leq \frac{2(M-1)}{\bb^\frac12},\\
  \vspace{0.2cm}
  z_2(m\bb)-z_1(\bb) \leq \frac{2(1-m)}{\bb^\frac12},\\
  z_2(m\bb)-z_1(M\bb) \leq \frac{2(M-m)}{\bb^\frac12}.
 \end{cases}
\label{estimzM2b}
\end{equation}
As $|f_1'(\ziib)| \geq \frac12 \bb^\frac32$ (thanks to \eqref{Minofp2} and \eqref{Minofp3}), \eqref{estimzM2b} entails that (also considering $(m,M)=(1-k,1+k)$)
$$
 K_2 \leq \Kbb +4\frac{\sigma^{-\frac12}}{\bb^\frac72} (M-m) +\frac4{\bb^2} \max(M-1,1-m) K_2 \leq \Kbb +\frac{8 k}{\bb^\frac72}\sigma^{-\frac12} +\frac{4k}{\bb^2}  K_2.
$$
Choosing $k=\frac{\bb^2}8$, we similarly obtain that:
\begin{equation}
K_2 \leq 2(\Kbb +\frac{\sigma^{-\frac12}}{\bb^\frac32}).
\label{estimKKA2}
\end{equation}
The auxiliary term $\Kbb$ (see \eqref{defK2b}) can be bounded with the same methods as its counterpart $\Jbb$ but using \eqref{Minofp2} and \eqref{Minofp3}:
$$
 \Kbb \leq \int_{z_2(M\bb)}^{z_2(m\bb)} \frac{dx}{1+\sigma \Big(\frac{\bb^\frac32}2 (x-\ziib)\Big)^2} =\frac2{\bb^\frac32} \sigma^{-\frac12} \int_{\sigma^{\frac12} \frac{\bb^\frac32}2 (z_2(M\bb)-\ziib)}^{\sigma^{\frac12} \frac{\bb^\frac32}2 (z_2(m\bb)-\ziib)} \frac{dz}{1+z^2} \leq \frac{2\pi}{\bb^\frac32}\sigma^{-\frac12}.
$$
Gathering this estimate with \eqref{estimKKA2} and \eqref{K1K3} (with $k=\frac{\bb^2}8$), we get:
$$
\Kb \leq \frac{12\pi}{\bb^2} R^3 \sigma^{-\frac12} +2\frac{2\pi+1}{\bb^\frac32} \sigma^{-\frac12}.
$$
Finally, using that $\bb \geq 2f_1(\sqrt{\R^2-1})= 2\frac{\sqrt{\R^2-1}}{R^3} \geq \frac1{R^2}$ when $R\geq \frac2{\sqrt{3}}$, we end-up with:
\begin{equation}
\Kb\leq 2(8\pi+1)R^7 \sigma^{-\frac12},
 \label{estimKKA3}
\end{equation}
To finish we focus on the last (second) case that we also separate into two subcases:
\begin{itemize}
 \item $z_2(2\bb) < \ziib \leq \sqrt{R^2-1} <z_2(\frac{\bb}2)  \quad \Longleftrightarrow \quad f_1(\sqrt{R^2-1}) \leq \bb < 2 f_1(\sqrt{R^2-1})$,
 \item $z_2(2\bb) \leq \sqrt{R^2-1} < \ziib <z_2(\frac{\bb}2)  \quad \Longleftrightarrow \quad \frac12 f_1(\sqrt{R^2-1}) \leq \bb < f_1(\sqrt{R^2-1})$.
\end{itemize}
In the first subcase, we have $z_2(\frac{\bb}2) \leq z_2\big(\frac12 f_1(\sqrt{R^2-1})\big) =z_2(\frac{\sqrt{R^2-1}}{2R^3}) \sim R\sqrt{2}\geq\sqrt{R^2-1}$ so that we can write:
$$
\Kb\leq \int_\frac1{\sqrt{2}}^{z_2(\frac{\sqrt{R^2-1}}{2R^3})} \frac{dx}{1+\sigma(f_1(x)-\bb)^2},
$$
and we can reproduce the method for the previous third case (which only slightly changes the constants but not the conditions on $\delta_0$) and obtain once more \eqref{estimKKA3}. The second subcase is treated identically but this time replacing the upper bound of the integral by $z_2(\frac{\sqrt{R^2-1}}{4R^3})\sim 2R$, which completes the second case.

Now combining \eqref{estimKKA3} with \eqref{estimJJA3}, we have proved that when $\delta_0 \in]0, \frac{\sqrt{3}}{5\sqrt{5}}]$ is so small that for all $\beta\in]0,\delta_0]$:
$$
 |f_1'(z_2(\frac{\bb}2))| \geq \frac12 \bb^\frac32.
$$
there exists a constant $C_{\delta_0}$ (that has a finite limit when $\delta_0$ goes to zero) such that for any $\beta \in ]0,\delta_0]$, $\IbR \leq C_{\delta_0} R^7 \sigma^{-\frac12}$. To conclude the proof we just recall that thanks to Point 2 from Proposition \ref{Propelem}, we can replace $C_{\delta_0}$ by $\max(C_0,C_{\delta_0})$ so that the previous estimates is also true for $\bb=0$. $\blacksquare$

\subsection{Study when $\bb=\frac2{3\sqrt{3}}$}
In this section we prove the following result.
\begin{prop}
 \sl{There exists $C_0>0$ such that for any $R\geq 1$ and $\sigma>0$,
 $$
 \IbRd =\int_0^{\sqrt{R^2-1}} \frac{dx}{1+\sigma \big(f_1(x)-f_1(\frac1{\sqrt{2}})\big)^2} \leq C_0 R^3 \sigma^{-\frac14}.
 $$
 Moreover, there exist $c_0, \sigma_0>0$ such that for any $R\geq \frac{\sqrt{3}}{\sqrt{2}}$ and $\sigma\geq \sigma_0$,
 $$
 \IbRd \geq c_0 \sigma^{-\frac14}.
 $$
 }
 \label{estimbetadeux}
\end{prop}
\textbf{Proof:} The reason why we do not obtain $\sigma^{-\frac12}$ in Proposition \ref{PropIab2} but $\sigma^{-\frac14}$ is that in the present case $f_1'(\frac1{\sqrt{2}})=0$ (but $f_1''(\frac1{\sqrt{2}})\neq 0$), which forces us to push one rank further the Taylor series expansion with integral remainder according to:
\begin{multline}
 f_1(x)-f_1(\frac1{\sqrt{2}})\\
 =\frac12 f_1''(\frac1{\sqrt{2}})(x-\frac1{\sqrt{2}})^2 +\frac12 (x-\frac1{\sqrt{2}})^3 \int_0^1 (1-u)^2 f_1'''\Big((1-u)\frac1{\sqrt{2}}+ux\Big) du.
 \label{Taylor3}
\end{multline}
We will then cut the integral as follows (for some $\frac12 <m<1$ close to $1$ in a way precised later) and adapt the previous method using \eqref{Taylor3}:
\begin{multline}
  \IbRd =I_1 +I_2 + I_3 \overset{def}{=} \int_0^{z_1(\frac{2m}{3\sqrt{3}})} \frac{dx}{1+\sigma \big(f_1(x)-f_1(\frac1{\sqrt{2}})\big)^2}\\
  +\int_{z_1(\frac{2m}{3\sqrt{3}})}^{z_2(\frac{2m}{3\sqrt{3}})} \frac{dx}{1+\sigma \big(f_1(x)-f_1(\frac1{\sqrt{2}})\big)^2} +\int_{z_2(\frac{2m}{3\sqrt{3}})}^{\sqrt{R^2-1}} \frac{dx}{1+\sigma \big(f_1(x)-f_1(\frac1{\sqrt{2}})\big)^2}.
\end{multline}
\begin{rem}
 \sl{As $m$ is close to $1$, $z_1(\frac{2m}{3\sqrt{3}})\sim \frac1{\sqrt{2}} \sim z_2(\frac{2m}{3\sqrt{3}})$.}
\end{rem}
For all $x\in [0, z_1(\frac{2m}{3\sqrt{3}})] \cup [z_2(\frac{2m}{3\sqrt{3}}), \sqrt{R^2-1}]$ we have $f_1(x)\leq \frac{2m}{3\sqrt{3}}$ so that:
$$
\frac2{3\sqrt{3}} -f_1(x) \geq \frac{1-m}{m} f_1(x) \geq (1-m)\frac{x}{R^3},
$$
and reproducing the arguments from the previous section we obtain that
\begin{equation}
I_1+I_3 \leq \int_0^{\sqrt{R^2-1}} \frac{dx}{1+\sigma (1-m)^2\frac{x^2}{R^6}} \leq \frac{\pi}{1-m} R^3 \sigma^{-\frac12}.
 \label{estimI1I3}
\end{equation}
The last piece is also decomposed as follows:
\begin{equation}
 I_2= J+K \overset{def}{=} \int_{z_1(\frac{2m}{3\sqrt{3}})}^{\frac1{\sqrt{2}}} \frac{dx}{1+\sigma \big(f_1(x)-f_1(\frac1{\sqrt{2}})\big)^2} +\int_{\frac1{\sqrt{2}}}^{z_2(\frac{2m}{3\sqrt{3}})} \frac{dx}{1+\sigma \big(f_1(x)-f_1(\frac1{\sqrt{2}})\big)^2}
\label{decpI2}
\end{equation}
We begin with $J$ and introduce:
\begin{equation}
\bJ \overset{def}{=} \int_{z_1(\frac{2m}{3\sqrt{3}})}^{\frac1{\sqrt{2}}} \frac{dx}{1+\sigma \big[\frac12 f_1''(\frac1{\sqrt{2}})(x-\frac1{\sqrt{2}})^2 \big]^2}.
 \label{defJb}
\end{equation}
Doing as in \eqref{estimJJA}, and using in \eqref{Taylor3} that $\|f_1'''\|_{L^\infty(\R_+)}\leq 9$, we obtain:
\begin{equation}
 |J-\bJ| \leq \frac32 \sigma \int_{z_1(\frac{2m}{3\sqrt{3}})}^{\frac1{\sqrt{2}}} \frac{|f_1(x)-f_1(\frac1{\sqrt{2}})|\cdot |x-\frac1{\sqrt{2}}|^3 +\frac12 |f_1''(\frac1{\sqrt{2}})|\cdot|x-\frac1{\sqrt{2}}|^5}{\left(1+\sigma \big(f_1(x)-f_1(\frac1{\sqrt{2}})\big)^2\right) \left(1+\frac{\sigma}4 |f_1''(\frac1{\sqrt{2}})|^2 |x-\frac1{\sqrt{2}}|^4\right)}dx.
\end{equation}
Using twice the Young estimates (with coefficients $(2,2)$ and $(4,\frac43)$) we can write:
$$
|f_1(x)-f_1(\frac1{\sqrt{2}})| \leq \frac{\sigma^{-\frac12}}2 \Big(1+\sigma \big(f_1(x)-f_1(\frac1{\sqrt{2}})\big)^2\Big),
$$
and
\begin{multline}
 |x-\frac1{\sqrt{2}}|^3= \left(\frac{\sigma}4 |f_1''(\frac1{\sqrt{2}})|^2\right)^{-\frac34} \left(\frac{\sigma}4 |f_1''(\frac1{\sqrt{2}})|^2 |x-\frac1{\sqrt{2}}|^4\right)^\frac34\\
 \leq \frac34 \left(\frac{\sigma}4 |f_1''(\frac1{\sqrt{2}})|^2\right)^{-\frac34} \left(1+\frac{\sigma}4 |f_1''(\frac1{\sqrt{2}})|^2 |x-\frac1{\sqrt{2}}|^4\right),
\end{multline}
which entails that
\begin{multline}
 |J-\bJ| \leq \frac32 \sigma \int_{z_1(\frac{2m}{3\sqrt{3}})}^{\frac1{\sqrt{2}}} \left(\frac3{2\sqrt{2}} \sigma^{-\frac54} \frac1{|f_1''(\frac1{\sqrt{2}})|^\frac32} + \frac2{\sigma} \frac1{|f_1''(\frac1{\sqrt{2}})|} \frac{|x-\frac1{\sqrt{2}}|}{1+\sigma \big(f_1(x)-f_1(\frac1{\sqrt{2}})\big)^2} \right) dx\\
 \leq \left(\frac1{\sqrt{2}}- z_1(\frac{2m}{3\sqrt{3}})\right) \left(\frac9{4\sqrt{2}} \frac{\sigma^{-\frac14}}{|f_1''(\frac1{\sqrt{2}})|^\frac32} +\frac3{|f_1''(\frac1{\sqrt{2}})|} J \right).
  \label{estimJJb}
 \end{multline}
 As in the previous part we need to be careful when estimating the following term (performing in the last line the change of variable $t=\frac2{3\sqrt{3}}-\eta$):
 \begin{multline}
  \frac1{\sqrt{2}}- z_1(\frac{2m}{3\sqrt{3}}) =z_1(\frac2{3\sqrt{3}})- z_1(\frac{2m}{3\sqrt{3}}) =\int_{\frac{2m}{3\sqrt{3}}}^{\frac2{3\sqrt{3}}} z_1'(t) dt =\int_{\frac{2m}{3\sqrt{3}}}^{\frac2{3\sqrt{3}}} \frac{dt}{f_1'(z_1(t))}\\
  =\int_0^{\frac2{3\sqrt{3}}(1-m)} \frac{d\eta}{f_1'(z_1(\frac2{3\sqrt{3}}-\eta))}.
  \label{estimdiff}
 \end{multline}
We recall that $f_1$ is strictly increasing on $]0, \frac1{\sqrt{2}}[$ but $f_1'(z_1(\frac2{3\sqrt{3}}-\eta))$ is small and we will need the following asymptotic expansions.

\begin{prop}
 \sl{When $\eta>0$ is small we have:
 \begin{equation}
  \begin{cases}
  \vspace{0.2cm}
z_1(\frac2{3\sqrt{3}}-\eta)= \frac1{\sqrt{2}}- \frac{3^\frac54}{2\sqrt{2}} \eta^\frac12 -\frac{25\sqrt{3}}{8\sqrt{2}}\eta +o(\eta),\\
 \vspace{0.2cm}
z_2(\frac2{3\sqrt{3}}-\eta)=\frac1{\sqrt{2}}+ \frac{3^\frac54}{2\sqrt{2}} \eta^\frac12 +\frac{7\sqrt{3}}{8\sqrt{2}}\eta +o(\eta),\\
 \vspace{0.2cm}
f_1'(z_1(\frac2{3\sqrt{3}}-\eta))\sim \frac{4\sqrt{2}}{3^\frac54} \eta^\frac12,\\
f_1'(z_2(\frac2{3\sqrt{3}}-\eta))\sim -\frac{4\sqrt{2}}{3^\frac54} \eta^\frac12.\\
  \end{cases}
 \end{equation}
\label{DL2}
 }
\end{prop}
Using the last two results, there exists $\eta_0>0$ such that for any $\eta\in]0,\eta_0]$, we have for $k\in\{1,2\}$
\begin{equation}
\frac{3^\frac54}8 \eta^{-\frac12} \leq \frac1{|f_1'(z_k(\frac2{3\sqrt{3}}-\eta))|} \leq \frac{3^\frac54}4 \eta^{-\frac12}.
\label{estimf1prime}
\end{equation}
If we ask that $0<(1-m)\leq \eta_0 \frac{3\sqrt{3}}2$, then
\begin{equation}
  \frac1{\sqrt{2}}- z_1(\frac{2m}{3\sqrt{3}}) \leq  \frac{3^\frac54}4 \int_0^{\frac2{3\sqrt{3}}(1-m)} \eta^{-\frac12} d\eta \leq \frac{\sqrt{3}}{\sqrt{2}} \sqrt{1-m}.
 \end{equation}
Plugging this into \eqref{estimJJb} we obtain:
\begin{equation}
 |J -\bJ| \leq\frac{9\sqrt{3}}{8} \frac{\sqrt{1-m}}{|f_1''(\frac1{\sqrt{2}})|^\frac32} \sigma^{-\frac14} +\frac{3\sqrt{3}}{\sqrt{2}} \frac{\sqrt{1-m}}{|f_1''(\frac1{\sqrt{2}})|} J.
  \label{estimJJbbis}
 \end{equation}
Choosing $m$ such that
\begin{equation}
 1-m =\min\left(\eta_0 \frac{3\sqrt{3}}2, \frac{|f_1''(\frac1{\sqrt{2}})|^2}{54}\right),
 \label{choixm}
\end{equation}
we obtain that
\begin{equation}
 J\leq 2\bigg(\bJ +\frac3{8\sqrt{2}} \frac{\sigma^{-\frac14}}{|f_1''(\frac1{\sqrt{2}})|^\frac12}\bigg).
 \label{estimJJb2}
\end{equation}
Now we need to precise $\bJ$, performing the change of variable $z=\big(\frac{\sigma}4\big)^\frac14 |f_1''(\frac1{\sqrt{2}})|^\frac12 (x-\frac1{\sqrt{2}})$, we get:
\begin{equation}
 \bJ= \big(\frac4{\sigma}\big)^\frac14 \frac1{|f_1''(\frac1{\sqrt{2}})|^\frac12} \int_{\big(\frac{\sigma}4\big)^\frac14 |f_1''(\frac1{\sqrt{2}})|^\frac12 (z_1(\frac{2m}{3\sqrt{3}})-\frac1{\sqrt{2}})}^0 \frac{dz}{1+z^4} \leq C_0 \sigma^{-\frac14}.
 \label{estimJbar}
\end{equation}
Going back to \eqref{estimJJb2}, we finally obtain that there exists some constant $C_0$ such that if \eqref{choixm} is fulfilled,
$$
J \leq C_0 \sigma^{-\frac14}.
$$
The second term $K$ from \eqref{decpI2} is dealt the same way as $I$ and gathering \eqref{decpI2}, \eqref{estimI1I3}, \eqref{choixm}, we obtain that $I_2\leq C_0 R^3 \sigma^{-\frac14}$ which gives the right inequality from Proposition \ref{estimbetadeux}.
\\

We now turn to the reverse inequality from Proposition \ref{estimbetadeux}: as $\IbRd \geq I_2 \geq J$ we can go back to \eqref{estimJJbbis} and write that (thanks to \eqref{estimJbar}):
\begin{multline}
 (1 +\frac{3\sqrt{3}}{\sqrt{2}} \frac{\sqrt{1-m}}{|f_1''(\frac1{\sqrt{2}})|}) J \geq \bJ -\frac{9\sqrt{3}}{8} \frac{\sqrt{1-m}}{|f_1''(\frac1{\sqrt{2}})|^\frac32} \sigma^{-\frac14}\\
 \geq \frac{\sqrt{2}}{|f_1''(\frac1{\sqrt{2}})|^\frac12} \left( \int_{\big(\frac{\sigma}4\big)^\frac14 |f_1''(\frac1{\sqrt{2}})|^\frac12 (z_1(\frac{2m}{3\sqrt{3}})-\frac1{\sqrt{2}})}^0 \frac{dz}{1+z^4} -\frac{9\sqrt{3}}{8\sqrt{2}} \frac{\sqrt{1-m}}{|f_1''(\frac1{\sqrt{2}})|} \right) \sigma^{-\frac14}.
 \end{multline}
If we now choose $m=m_0<1$ sufficiently close to $1$ such that
\begin{equation}
 \sqrt{1-m_0} \leq \min\left(\sqrt{\eta_0 \frac{3\sqrt{3}}2},\frac{4\sqrt{2}}{9\sqrt{3}} \int_{-1}^0 \frac{dz}{1+z^4}\right).
\end{equation}
We recall that $z_1(\frac{2m_0}{3\sqrt{3}})<\frac1{\sqrt{2}}$, so if $\sigma>0$ is so large that $\big(\frac{\sigma}4\big)^\frac14 |f_1''(\frac1{\sqrt{2}})|^\frac12 (z_1(\frac{2m_0}{3\sqrt{3}})-\frac1{\sqrt{2}}) \leq -1$ (that is $\sigma\geq \sigma_0=4 |f_1''(\frac1{\sqrt{2}})|^{-2} |z_1(\frac{2m_0}{3\sqrt{3}})-\frac1{\sqrt{2}}|^{-4}$) then:
\begin{equation}
 J \geq (1 +\frac{3\sqrt{3}}{\sqrt{2}} \frac{\sqrt{1-m_0}}{|f_1''(\frac1{\sqrt{2}})|})^{-1} \frac{\sqrt{2}}{|f_1''(\frac1{\sqrt{2}})|^\frac12} \left(\frac12 \int_{-1}^0 \frac{dz}{1+z^4}\right) \sigma^{-\frac14}.
 \end{equation}
Note that the condition $R\geq \frac{\sqrt{3}}{\sqrt{2}}$ in the statement of Proposition \ref{estimbetadeux} is needed only to ensure that $\sqrt{R^2-1} \geq \frac1{\sqrt{2}}$ so that $\IbRd \geq I_2 \geq J$. $\blacksquare$

\subsection{Study for $\bb\geq \frac2{3\sqrt{3}}$}
The aim of this section is to prove the following result (much more useful than the first point from Proposition \ref{Propelem}):

\begin{prop}
 \sl{For any $\bb \geq\frac2{3\sqrt{3}}$, we have $\IbR \leq 4 \IbRd \leq 4C_0 R^3 \sigma^{-\frac14}$.}
 \label{estimbetagrand}
\end{prop}
\textbf{Proof:} following the steps of the proof of Proposition \ref{Propelem}, we have for any $\bb>\frac2{3\sqrt{3}}$:
\begin{multline}
 |\IbR -\IbRd| \leq 2\sigma |\bb-\frac2{3\sqrt{3}}| \int_0^{\sqrt{R^2-1}} \frac{|f_1(x)-\beta|+\frac12|\beta-\frac2{3\sqrt{3}}|}{\sigma(f_1(x)-\bb)^2(1+\sigma(f_1(x)-\frac2{3\sqrt{3}})^2)} dx\\
 \leq 2|\bb-\frac2{3\sqrt{3}}| \int_0^{\sqrt{R^2-1}} \left(\frac1{|f_1(x)-\beta|} +\frac12 \frac{|\bb-\frac2{3\sqrt{3}}|}{|f_1(x)-\beta|^2} \right) \frac{dx}{1+\sigma(f_1(x)-\frac2{3\sqrt{3}})^2}.
\end{multline}
As $\bb>\frac2{3\sqrt{3}}$, for any $x\in[0, \sqrt{R^2-1}]$, we have (as both terms below are nonnegative):
$$
|f_1(x)-\beta|=(\beta-\frac2{3\sqrt{3}}) +(\frac2{3\sqrt{3}}-f_1(x)) \geq
\beta-\frac2{3\sqrt{3}} =|\beta-\frac2{3\sqrt{3}}|,
$$
so that
\begin{multline}
 |\IbR -\IbRd|\\
 \leq 2|\bb-\frac2{3\sqrt{3}}| \int_0^{\sqrt{R^2-1}} \frac32 \frac1{|\bb-\frac2{3\sqrt{3}}|} \frac{dx}{1+\sigma(f_1(x)-\frac2{3\sqrt{3}})^2} = 3 \IbRd,
\end{multline}
which entails the result as the inequality remains true when $\bb=\frac2{3\sqrt{3}}$. $\blacksquare$

\subsection{Study when $\bb=\frac2{3\sqrt{3}}-\eta$}
This part is devoted to the proof of the following result:

\begin{prop}
 \sl{
 There exists $C_0>0$ such that for any $\beta \in [\frac2{3\sqrt{3}}-\eta_0, \frac2{3\sqrt{3}}[$ (with $\eta_0$ defined in \eqref{estimf1prime}):
 $$
 \IbR \leq C_0 \left(\frac{\sigma^{-\frac12}}{|\beta-\frac2{3\sqrt{3}}|^\frac12}+R^3 \sigma^{-\frac14}\right).
 $$
 }
 \label{estimbetaproche}
\end{prop}
\textbf{Proof:} let us begin by splitting the integral as in \eqref{Split}. Next, writing $\bb=\frac2{3\sqrt{3}}-\eta$, we decompose $\Jb$ a little differently: for $\frac12\leq m<1<M \leq 2$ (precised later) we define:
\begin{multline}
 \Jb= J_1 +J_2+ J_3 \overset{def}{=} \int_0^{z_1(\frac2{3\sqrt{3}}-M\eta)} \frac{dx}{1+\sigma(f_1(x)-\bb)^2}\\
 +\int_{z_1(\frac2{3\sqrt{3}}-M\eta)}^{z_1(\frac2{3\sqrt{3}}-m\eta)} \frac{dx}{1+\sigma(f_1(x)-\bb)^2} +\int_{z_1(\frac2{3\sqrt{3}}-m\eta)}^\frac1{\sqrt{2}} \frac{dx}{1+\sigma(f_1(x)-\bb)^2}.
\end{multline}
Improving what is done before we introduce:
$$
\Jub \overset{def}{=} \int_0^{z_1(\frac2{3\sqrt{3}}-M\eta)} \frac{dx}{1+\sigma\big(f_1(x)-\frac2{3\sqrt{3}}\big)^2}.
$$
For all $x\in[0,z_1(\frac2{3\sqrt{3}}-M\eta)]$, $f_1(x) \leq \frac2{3\sqrt{3}}-M\eta=\beta+(1-M)\eta$, so
\begin{equation}
 |f_1(x)-\beta| \geq (M-1)\eta,
\end{equation}
and we obtain that:
\begin{multline}
 |J_1-\Jub| \leq 2\sigma \eta \int_0^{z_1(\frac2{3\sqrt{3}}-M\eta)} \frac{|f_1(x)-(\frac2{3\sqrt{3}}-\eta)|+\frac{\eta}2}{\Big(1+\sigma\big(f_1(x)-(\frac2{3\sqrt{3}}-\eta)\big)^2\Big)\Big(1+\sigma\big(f_1(x)-\frac2{3\sqrt{3}}\big)^2\Big)} dx\\
 \leq 2 \eta \int_0^{z_1(\frac2{3\sqrt{3}}-M\eta)} \left(\frac1{|f_1(x)-\bb|} +\frac{\eta}2\frac1{|f_1(x)-\bb|^2}\right) \frac{dx}{1+\sigma\big(f_1(x)-\frac2{3\sqrt{3}}\big)^2}\\
 \leq \left(\frac2{M-1} +\frac1{(M-1)^2}\right) \Jub.
\end{multline}
As $\Jub \leq \IbRd$ and considering that $0<M-1$ will be chosen small:
\begin{equation}
 J_1\leq \left(1+\frac2{M-1} +\frac1{(M-1)^2}\right) \Jub \leq \frac4{(M-1)^2} \IbRd \leq \frac4{(M-1)^2} C_0 R^3 \sigma^{-\frac14}.
 \label{estimJ1bis}
\end{equation}
Similarly, we get that (with $0<1-m$ small):
\begin{equation}
 J_3\leq \frac4{(1-m)^2} C_0 R^3 \sigma^{-\frac14}.
  \label{estimJ3bis}
\end{equation}
To deal with the last piece, we introduce:
$$
\Jbb \overset{def}{=} \int_{z_1(\frac2{3\sqrt{3}}-M\eta)}^{z_1(\frac2{3\sqrt{3}}-m\eta)} \frac{dx}{1+\sigma[f_1'(\zib)(x-\zib)]^2}.
$$
Reproducing the arguments that lead to \eqref{estimJJA} we obtain that (with $\bb=\frac2{3\sqrt{3}}-\eta$):
\begin{multline}
 |J_2-\Jbb| \leq \frac{\sigma^{-\frac12}}{2|f_1'(\zib)|^2} |z_1(\frac2{3\sqrt{3}}-m\eta)-z_1(\frac2{3\sqrt{3}}-M\eta)|\\
 +\frac{\max\Big(|z_1(\frac2{3\sqrt{3}}-m\eta)-z_1(\bb)|, |z_1(\bb)-z_1(\frac2{3\sqrt{3}}-M\eta)|\Big)}{|f_1'(\zib)|} J_2.
 \label{estimJJbb}
\end{multline}
As in \eqref{estimdiff}, and when $\eta\leq \eta_0$ where $\eta_0$ has been introduced in \eqref{estimf1prime}
\begin{multline}
  z_1(\frac2{3\sqrt{3}}-m\eta)-z_1(\bb) =\int_{\frac2{3\sqrt{3}}-\eta}^{\frac2{3\sqrt{3}}-m\eta} \frac{dt}{f_1'(z_1(t))} =\int_{m\eta}^{\eta} \frac{du}{f_1'(z_1(\frac2{3\sqrt{3}}-u))}\\
  \leq \int_{m\eta}^{\eta} \frac{3^\frac54}4 u^{-\frac12}du \leq \frac{3^\frac54}2 (1-\sqrt{m}) \eta^\frac12 \leq \frac{3^\frac54}2 (1-m) \eta^\frac12.
 \end{multline}
Similarly,
$$
\begin{cases}
\vspace{0.1cm}
 z_1(\bb)-z_1(\frac2{3\sqrt{3}}-M\eta) \leq \frac{3^\frac54}2 (M-1) \eta^\frac12,\\
 z_1(\frac2{3\sqrt{3}}-m\eta)-z_1(\frac2{3\sqrt{3}}-M\eta) \leq \frac{3^\frac54}2 (M-m) \eta^\frac12,
\end{cases}
$$
so that plugging into \eqref{estimJJbb}, using the bound from \eqref{estimf1prime} and choosing $(m,M)=(1-k,1+k)$,
$$
|J_2-\Jbb| \leq \frac{\sigma^{-\frac12}}{\eta^\frac12} \frac{3^\frac{15}4}{32}k +\frac{3^\frac52}8k J_2.
$$
Performing the change of variable $z=\sigma\frac12 |f_1'(\zib)|(x-\zib)$, we easily bound $\Jbb\leq C_0 \sigma^{-\frac12} \eta^{-\frac12}$ so that choosing $k\leq 4\cdot 3^{-\frac52}$, we finally obtain that
$$
J_2 \leq C_0 \sigma^{-\frac12} \eta^{-\frac12}
$$
Gathering it with \eqref{estimJ1bis} and \eqref{estimJ3bis} ($m,M$ now precised with the choice of $k$), we have proved that:
$$
\Jb \leq C_0 (\sigma^{-\frac12} \eta^{-\frac12}+R^3 \sigma^{-\frac14}).
$$
The case of $\Kb$ is similar and the proof of Proposition \ref{estimbetaproche} is finished. $\blacksquare$

\subsection{Study when $\beta \in [\frac2{3\sqrt{3}}-\eta_0, \frac2{3\sqrt{3}}]$}
We prove the following result:

\begin{prop}
 \sl{
 With $\eta_0$ defind in \eqref{estimf1prime}, there exists $C_0>0$ such that for any $\beta \in [\frac2{3\sqrt{3}}-\eta_0, \frac2{3\sqrt{3}}]$:
 $$
 \IbR \leq C_0 R^3 \sigma^{-\frac14}.
 $$
 }
 \label{estimbetaproche2}
\end{prop}
\textbf{Proof:} the difference with the previous subsection is that we have to overcome the singularity when $\beta=\frac2{3\sqrt{3}}$. This is simply done coupling the previous proposition with Proposition \ref{Propelem} (Point 3), for any $\beta \in [\frac2{3\sqrt{3}}-\eta_0, \frac2{3\sqrt{3}}]$:
\begin{itemize}
 \item either $|\beta-\frac2{3\sqrt{3}}| \leq \sigma^{-\frac12}$ and $\IbR \leq 3\IbRd \leq C_0 R^3 \sigma^{-\frac14}$,
 \item or $|\beta-\frac2{3\sqrt{3}}| \geq \sigma^{-\frac12}$, and from Proposition \ref{estimbetaproche} we also have $\IbR \leq C_0 R^3 \sigma^{-\frac14}$. $\blacksquare$
\end{itemize}

\subsection{Study when $\beta \in [\delta,\frac2{3\sqrt{3}}-\delta]$}
We prove the following result:
\begin{prop}
 \sl{For any $\delta\in]0,\frac1{3\sqrt{3}}[$, there exists $D_{\delta}>0$ (going to infinity when $\delta$ goes to zero) such that for any $\beta \in [\delta,\frac2{3\sqrt{3}}-\delta]$,
 $$
 \IbR \leq D_{\delta} R^3 \sigma^{-\frac12}.
 $$
 }
 \label{estimbetaintermediaire}
\end{prop}
\textbf{Proof:} as before let us begin by splitting the integral as in \eqref{Split}, and boths parts $\Jb$ and $\Kb$ are also decomposed as in \eqref{DecompoJ} and \eqref{DecompoK} for some $\frac12<m<1<M$ that we choose as before according to $(m,M)=(1-k,1+k)$ with $k\in]0,1[$ small. When $\beta \in [\delta,\frac2{3\sqrt{3}}-\delta]$, we have
$$
\frac{\delta}2 \leq m\delta \leq m\bb \leq \bb \leq M\bb \leq M(\frac2{3\sqrt{3}}-\delta)=(1+k)(\frac2{3\sqrt{3}}-\delta),
$$
and
\begin{equation}
 (1+k)(\frac2{3\sqrt{3}}-\delta) \leq \frac2{3\sqrt{3}}-\frac{\delta}2 \Longleftrightarrow k\leq \frac{\delta}2 \frac1{\frac2{3\sqrt{3}}-\delta}.
 \label{condkdelta}
\end{equation}
In the decomposition \eqref{DecompoJ}, \eqref{estimb0J1} and \eqref{estimb0J3} are still true as well as \eqref{estimJJA}. As $t \mapsto |f_1'(z_1(t))|$ is continuous on the compact set $[\frac{\delta}2,\frac2{3\sqrt{3}}-\frac{\delta}2]$ (and always positive) it is bounded from below by some positive constant $c_{\delta}>0$. And thanks to \eqref{condkdelta}, for all $x\in [m\bb, M\bb] \subset [\frac{\delta}2,\frac2{3\sqrt{3}}-\frac{\delta}2]$, $|f_1(z_1(x))|\geq c_{\delta}$, and \eqref{zM} turns into (the other two terms are transformed similarly):
$$
z_1(M\bb)-z_1(\bb) \leq \frac{(M-1)\bb}{c_{\delta}},
$$
and \eqref{estimJ2mM} turns into ($\bb \leq \frac2{3\sqrt{3}} \leq 1$):
\begin{equation}
 J_2 \leq \Jbb +\sigma^{-\frac12} \frac{\bb k}{c_{\delta}^3} +\frac{\beta k}{c_{\delta}^2} J_2 \leq \Jbb +\sigma^{-\frac12} \frac{k}{c_{\delta}^3} +\frac{k}{c_{\delta}^2} J_2 \leq \Jbb +\frac{\sigma^{-\frac12}}{2c_{\delta}} +\frac12 J_2,
\end{equation}
if $k$ satisfies:
$$
k=\min\left(\frac12,\frac{\delta}2 \frac1{\frac2{3\sqrt{3}}-\delta}, \frac{c_{\delta}^2}2\right).
$$
Similarly \eqref{estimJbar2} turns into $\Jbb\leq \frac{2\pi}{c_{\delta}} \sigma^{-\frac12}$, so that we finally obtain (replacing $m,M$ for the choice of $k$ in \eqref{estimb0J1} and \eqref{estimb0J3}) that there exists some constant $D_{\delta}>0$ (still going to infinity when $\delta$ goes to zero) such that
$$
\Jb \leq D_{\delta} \sigma^{-\frac12}.
$$
The other integral $\Kb$ is treated similarly but due to the bounds, it is estimated as follows:
$$
\Jb \leq D_{\delta} R^3 \sigma^{-\frac12},
$$
which ends the proof. $\blacksquare$

\subsection{Conclusion}
We are now able to prove Proposition \ref{PropIab2} (which implies Proposition \ref{PropIab} thanks to the scaling argument). Fixing $\delta=\min(\delta_0, \eta_0)>0$ (the parameters from Propositions \ref{estimbetapetit} and \ref{estimbetaproche2}), we apply Proposition \ref{estimbetaintermediaire} to this choice for $\delta_0$. Using Proposition \ref{estimbetagrand} we then have obtained that there exists a constant $C_0$ such that for any $\sigma>0$ and $R\geq \frac1{\sqrt{3}}$
$$
\sup_{\bb\geq 0} \IbR \leq C_0 R^7 (\sigma^{-\frac12}+ \sigma^{-\frac14}).
$$
When $\sigma \geq 1$, this turns into (still denoting as $C_0$ the constant)
$$
\sup_{\bb\geq 0} \IbR \leq C_0 R^7 \sigma^{-\frac14} =C_0 R^7 \min\big(1,\sigma^{-\frac14}\big).
$$
When $\sigma \leq 1$, from the definition of $\IbR$ it is obvious that
$$
\IbR \leq R \leq C_0 R^7 =C_0 R^7 \min\big(1,\sigma^{-\frac14}\big),
$$
which concludes the proof, the bound from below being done in the proof of Proposition \ref{estimbetadeux}. $\blacksquare$
\begin{rem}
 \sl{Another way of proving Proposition \ref{PropIab2} would consist in performing in $\Jb$ the change of variable $z=\sigma^{\frac12} (f_1(x)-\beta) \Longleftrightarrow x=z_1(\beta+\sigma^{-\frac12}z)$ so that:
 $$
 \Jb =\sigma^{-\frac12} \int_{-\beta\sigma^\frac12}^{(\frac2{3\sqrt{3}}-\beta)\sigma^\frac12} \frac1{f_1'(z_1(\beta+\sigma^{-\frac12}z))} \frac{dz}{1+z^2}.
 $$
 Thanks to this we simply see that for any fixed  $\beta \in ]0,\frac2{3\sqrt{3}}[$,
 $$
 \sigma^\frac12 \Jb \underset{\sigma \rightarrow \infty}{\longrightarrow} \frac{\pi}{f_1'(\zib)},
 $$
 which is of course problematic when $\bb=\frac2{3\sqrt{3}}$. But if we want to obtain precise bounds for any $\sigma>0$ more work is needed (having in mind that the function is not integrable on $\R$). Much more work is necessary for $\Kb$ as large parameters $R,\sigma$ are mixed with $\bb$ (in the previous limit everything is fixed except $\sigma$). On top of that, it is still needed to treat separately the case $\bb=\frac2{3\sqrt{3}}$.
 }
\end{rem}


\end{document}